\documentclass[10pt]{article}
\usepackage{amsmath}
\usepackage{amssymb}
\usepackage{amsthm}
\usepackage{amsbsy}
\usepackage{graphicx}
\usepackage{subfigure}
\usepackage{verbatim}
\usepackage{enumerate}

\font\elevenss=cmss11

\font\eightss=cmss8

\font\sixss=cmss8 at 6pt

\newfam\ssfam
\textfont\ssfam=\elevenss \scriptfont\ssfam=\eightss
 \scriptscriptfont\ssfam=\sixss
\def\ss{\fam\ssfam \elevenss}%

\theoremstyle{plain}
\newtheorem{thm}{Theorem}[section]
\newtheorem{lem}[thm]{Lemma}
\newtheorem{pr}[thm]{Proposition}
\newtheorem{cor}[thm]{Corollary}
\newtheorem{defn}[thm]{Definition}

\newtheorem{conj}[thm]{Conjecture}
\newtheorem{property}{Property}[section]
\newtheorem{example}[thm]{Example}
\theoremstyle{remark}

\newtheorem*{unremark}{Remark}

\newcommand{\Em}[1]{\textbf{#1}}

\def\ee{\epsilon}

\def\vv{{\bf v}}
\def\uu{{\bf u}}
\def\xx{{\bf x}}
\def\yy{{\bf y}}
\def\zz{{\bf z}}

\def\ZZ{{\bf Z}}

\def\ww{{\bf w}}
\def\rr{{\bf r}}

\def\bb{{\bf b}}

\def\tor3{{\mathbb T}^3}

\def\cK{{\bf K}}

\def\CC{{\cal C}}

\def\sing{{\cal V}}

\def\one{{\bf 1}}
\def\zero{{\bf 0}}

\def\yhat{\hat{\yy}}

\def\nbd{{\cal N}}
\def\B{{\cal B}}
\def\Cox{\hfill \Box}
\def\C{{\mathbb C}}
\def\Z{{\mathbb Z}}

\def\R{{\mathbb R}}
\def\CP{{\mathbb C}{\mathbb P}}

\def\Real{{\rm Re\,}}

\def\Log{{\rm ReLog\,}}

\def\amoeba{{\mbox{\ss amoeba}}}

\def\pt{\tilde{p}}

\def\Det{{\rm Det \,}}
\def\Per{{\rm Per \,}}
\def\trace{{\rm Tr \,}}
\def\Cap{{\rm Cap \,}}

\def\homog{\mbox{\ss hom}}

\def\disp{\displaystyle}

\def\SS{{\cal S}}

\def\P{{\mathbb P}}
\def\E{{\mathbb E}}

\def\|{{\, | \, }}
\def\LP{{\cal L}\hspace{-3pt}-\hspace{-3pt}{\cal P}}
\def\LPI{{\cal L}\hspace{-3pt}-\hspace{-3pt}{\cal P}I}
\def\blambda{{\boldsymbol{\lambda}}}
\def\halfplane{{\cal H}}
\def\mcovers{\triangleright}
\def\covers{\; \cdot \kern-5pt >}
\def\diag{{\rm diag}\,}
\def\lap{{\cal L}}
\def\hadamard{\bullet}

\makeatletter \def\romenumi{ \def\theenumi{\roman{enumi}}
\def\p@enumi{\theenumi} \def\labelenumi{(\@roman\c@enumi)}} 
\makeatother

\begin{document}
\renewcommand{\thepage}{\roman{page}}

\begin{titlepage}
\begin{center}
{\large \bf Hyperbolicity and stable polynomials in combinatorics
and probability} \\
\end{center}
\vspace{5ex}
\begin{flushright}
Robin Pemantle \footnote{Supported in part by National Science Foundation 
grant \# DMS 0905937}$^,$\footnote{University of Pennsylvania, 
Department of Mathematics, 209 S. 33rd Street, Philadelphia,
PA 19104 USA, pemantle@math.upenn.edu}
\end{flushright}

\vfill

\noindent{\bf ABSTRACT:} These lectures survey the theory of 
hyperbolic and stable polynomials, from their origins 
in the theory of linear PDE's to their present uses
in combinatorics and probability theory.
\hfill \\[1ex]

\vfill



\end{titlepage}

\tableofcontents

\pagebreak

\renewcommand{\thepage}{\arabic{page}}
\setcounter{page}{1}

\setcounter{equation}{0}
\section{Introduction} \label{sec:intro} 

These lectures concern the development and uses of two properties, 
\Em{hyperbolicity} and \Em{stability}.  Hyperbolicity has been
used chiefly in geometric and analytic contexts, concerning 
wave-like partial differential equations, lacunas, and generalized 
Fourier transforms.  Stability has its origins in control theory
but has surfaced more recently in combinatorics and probability, where
its algebraic properties, such as closure under various operations,
are paramount.  We begin with the definitions.
\begin{defn}[hyperbolicity] \label{def:hyperbolicity}
A homogeneous complex $d$-variable polynomial $p(z_1 , \ldots , z_d)$ 
is said to be hyperbolic in direction $\xx \in \R^d$ if and only if 
\begin{equation} \label{eq:def 1}
p(\xx + i \yy) \neq 0 \mbox{ for all } \yy \in \R^d  \, .
\end{equation}
A non-homogeneous $d$-variable polynomial $q$ with leading 
homogeneous part $p$ is said to be hyperbolic if and only if 
$p(\xx) \neq 0$ and there is some $t_0 > 0$ such that
\begin{equation} \label{eq:def hyperbolic}
q(i t \xx + \yy) \neq 0 \mbox{ for all } \yy \in \R^d  \mbox{ and real }
   t < t_0 \, .
\end{equation}
A polynomial is said to be hyperbolic if it is hyperbolic in direction
$\xi$ for some $\xi \in \R^d$.
\end{defn}

\begin{defn}[stability] \label{def:stable}
The complex polynomial $q$ is said to be \Em{stable} 
if and only if 
\begin{equation} \label{eq:def stable}
\Im \{ z_j \} > 0 \mbox{ for all } j = 1 , \ldots , d
   \mbox{ implies } q(z_1 , \ldots , z_d) \neq 0 \, .
\end{equation}
A \Em{real stable} polynomial is one that is both real and stable. 
\end{defn}

The following relation holds between these two definitions.
\begin{pr} \label{pr:orthant}
A real homogeneous polynomial $p$ is stable if and only if it is 
hyperbolic in direction $\xi$ for every $\xi$ in the positive orthant.
\end{pr}
Hyperbolicity is a real geometric property, that is, it is 
invariant under the general linear group.  It is not invariant 
under complex linear transformations, rather it distinguishes
the real subspace of $\C^d$.  Stability is invariant under yet
fewer maps than is hyperbolicity, due to the fact that it 
distinguishes positive from negative.  Stability is invariant 
under coordinatewise transformations of the upper half-plane, 
such as inversion $z \mapsto -1/z$ and dilation 
$z \mapsto \lambda z$, from which one can get a 
surprising amount of mileage.  

The notion of hyperbolicity for polynomials has been around since its 
introduction by L.\ G{\aa}rding sixty years 
ago~\cite{garding-hyperbolic}.\footnote{G{\aa}rding credits many
of the ideas to I.\ G.\ Petrovsky in his seminal paper~\cite{petrovsky}.}
The notion of stability for real or complex polynomials in several variables
goes back at least to multivariate versions of Hurwitz stability in
the 1980's, but appears to have become established as tools for
combinatorialists only in the last five
years~\cite{COSW,branden07,borcea-branden-liggett,gurvits-vdW}.  
These two related notions have now been applied in several 
seemingly unrelated fields.  

G{\aa}rding's  original purpose was to find the right condition to
guarantee analytic stability of a PDE as it evolves in time from an 
original condition in space; the notion of stability here is that
small perturbations of the initial conditions should produce small
perturbations of the solution at time $t$.  That G{\aa}rding found 
the right definition is undeniable because his result 
(Theorem~\ref{th:garding} below) is an equivalence.  
His proof uses properties of hyperbolic functions 
and their cones that are derived directly from the definitions.
Specfically, in order to construct solutions to boundary value
problems via the \Em{Riesz kernel}, an inverse Fourier transform
is computed on a linear space $\{ \xx + i \yy : \yy \in \R^d \}$;
hyperbolicity in the direction $\xx$ is used to ensure that
this space does not intersect the set $\{ q = 0 \}$ where 
$q^{-a}$ is singular.  G{\aa}rding developed some properties
of hyperbolic polynomials in this paper and later devoted a
separate paper to the further development of some
inequalities~\cite{garding-inequalities}.

Twenty years later, with Atiyah and Bott~\cite{ABG}, this work
was extended considerably.  They were able to compute inverse 
Fourier transforms more explicitly for a number of homogeneous
hyperbolic functions.  To do this, they developed semi-continuity 
properties of hyperbolic functions.  These properties enable
the construction of certain vector fields and deformations, which
in turn are used to deform the plane $\{ \xx + i \yy : \yy \in \R^d \}$ 
into a cone on which the inverse Fourier transforms are directly
integrable.  Direct integrability provides useful estimates and also
enables explicit computation via topological dimension reduction 
theorems of Leray and Petrovsky.  My introduction to this
exceptional body of work came in the early 2000's when I needed
to apply these deformations to multivariate analytic 
combinatorics; refining the relatively crude techniques 
in~\cite{PW1,PW2} led to~\cite{BP-cones}.

In control theory, a continuous-time (respectively discrete-time)
system whose transfer function is rational will be stable if 
its poles all lie in the open left half-plane (respectively 
the open unit disk).  Accordingly, we say a univariate polynomial 
is \Em{Hurwitz stable} (respectively \Em{Schur stable})
if its zeros are all in the open left half-plane (respectively 
the open unit disk).  
Multivariate generalizations of Hurwitz and Schur stability 
have arisen in a variety of problems.  The development of 
these generalizations splits into two veins.  The review 
paper~\cite{sokal-review} surveys a number of ``hard'' results 
on zero-free regions, meaning results that give good estimates 
on zero-free regions that depend on specific parameters of the 
models or graphs from which the polynomial is formed.  
Most relevant to the applications in Sections~\ref{sec:negdep} 
and~\ref{sec:determinants} are ``soft'' results, which give 
simple domains free of zeros (such as products of half planes), 
valid for all graphs or for wide classes of graphs.  

The development of this vein of multivariate stable 
function theory took place mainly in the area of 
statistical physics and combinatorial extensions thereof.
Multivariate Schur stability arises 
in the celebrated Lee--Yang Theorem for ferromagnetic Ising 
models~\cite{lee-yang1,lee-yang2}, where it implies the absence 
of phase transitions at nonzero magnetic field.  Likewise,
multivariate Hurwitz stability arises in electrical circuit
theory~\cite{fettweis85,COSW,sokal-lectures} and matroid 
generalizations thereof~\cite{COSW,branden07,wagner-wei},
and has applications to combinatorial 
enumeration~\cite{wagner06,wagner-rayleigh}.  It also
arises in generalizations of the Lee--Yang Theorem~\cite{lieb-sokal}
and in the Heilmann-Lieb Theorem for matching polynomials
(also known as monomer-dimer models)~\cite{heilmann-lieb,COSW}.
It should be noted that Hurwitz stability~\cite{hurwitz},
differs from the notion in Definition~\ref{def:stable} in
that the zero-free region is the right half-plane rather than 
the upper half-plane.  For general complex polynomials 
the two notions are equivalent under a linear change of
variables, but for real polynomials the two notions are
very different.
Except for brief portions of Sections~\ref{sec:stability 1}
and~\ref{sec:stability 2}, our concern will be with
stability as defined in Definition~\ref{def:stable}

Multivariate stability behaves nicely under certain transformations.  
The closure properties were investigated further by Borcea and
Br\"and\'en~\cite{borcea-branden-LYPS1}.  These closure
properties turn out to have powerful implications for systems
of negatively dependent random variables, as worked out in 2009
by Borcea, Br\"and\'en and Liggett~\cite{borcea-branden-liggett}.
My second brush with this subject was when I needed to apply some
of their results to determinantal point processes~\cite{PP-rayleigh}.

These two encounters with hyperbolic/stable polynomials in seemingly
dissimilar areas of mathematics prompted my search for a unified
understanding.  As pointed out at the end of the introduction
of~\cite{borcea-branden-liggett}, such a viewpoint was taken
25 years ago by Gian-Carlo Rota.
\begin{quote}
``The one contribution of mine that I hope will be remembered has
consisted in just pointing out that all sorts of problems of 
combinatorics can be viewed as problems of the locations of
zeros of certain polynomials...''
\end{quote}
to which we would add ``and also problems in differential equations,
number theory, probability, and perhaps many more areas.''  

Presently, the greatest interest in the subject of hyperbolic/stable
polynomials seems to be its potential for unlocking some combinatorial
conjectues concerning determinants such as the Bessis--Moussa--Villani 
conjecture, Lieb's ``permanent on top'' conjecture and extensions of 
the van der Waerden conjecture, proved in 1979, but still an active
subject (see, e.g.,~\cite{gurvits-vdW-06}).  The distances between
areas of mathematics in which hyperbolicity plays a role are such
that it was difficult to find a primary subject classification for
these lectures.  It also poses a unique challenge for me, as I am
an expert in none of these, save for applications to analytic 
combinatorics.  I therefore apologize in advance for any historical
inaccuracies perpetrated here, or for idiosyncratic viewpoints
arising from my personal history with the subject.  Here follows
an attempt to lay out its two main pillars, hyperbolicity and stability,
and to follow their progress up to the present day.

\pagebreak

\begin{center}
{\Huge Part I: Hyperbolicity}
\end{center}
\vspace{1in}

\setcounter{equation}{0}
\section{Origins, definitions and properties} 
\label{sec:hyperbolicity} 

\subsection{Relation to the propagation of wave-like equations}

G{\aa}rding's objective was to prove stability results for wave-like
partial differential equations.  Endow the space of smooth complex
valued functions of $d$ real variables, denoted $C^\infty (\R^d)$,
with the topology of uniform convergence on compact sets of all
the partial derivatives.  For any polynomial $q$ in $d$ variables,
let $D[q] = q(\partial / \partial \xx)$ denote the linear differential
operator obtained by formally substituting $(\partial / \partial x_i)$
for $x_i$ in $q$.  If we consider $\R^d$ as spacetime with the positive 
time direction given by a vector $\xi$, and denote by $H_\xi$ the
hyperplane orthogonal to $\xi$, then the boundary value equation
\begin{equation} \label{eq:diffeq}
D[q] (f) = 0, \hspace{0.5in} f = g \mbox{ on } H_\xi 
\end{equation}
on the halfspace $\{ \xx : \xx \cdot \xi \geq 0 \}$ may be 
interpreted as the evolution of $g$ under the equation $D[q] (f) = 0$ 
with initial condition $g$.  One would expect $m-1$ further initial 
conditions, such as $(\partial^j / \partial x_\xi^j) f = g_j$ for 
$1 \leq j \leq m-1$, where $(\partial / \partial x_\xi)$ denotes 
a derivative in direction $\xi$ and $m$ is the degree of $q$.

We consider~\eqref{eq:diffeq} to be \Em{stable under perturbations}
(in time direction $\xi$) if 
for all sequences $\{ f_n \}$ of functions in $C^\infty (\R^d)$
satisfying $D[q] (f_n) = 0$, convergence to zero of the restrictions
of $f_n$ to $H_\xi$ implies convergence to zero of $f_n$ on
the whole space.  Intuitively, arbitrary small change in the 
initial conditions cannot lead to macroscopic changes in the
evolution at a later time.  The first of G{\aa}rding's results
can be stated as follows. 

\begin{thm}[\protect{\cite[Theorem~III]{garding-hyperbolic}}]
\label{th:garding}
The differential equation~\eqref{eq:diffeq} is stable in direction $\xi$
if and only if $q$ is hyperbolic in direction $\xi$.
\end{thm}

It should be noted that the
definition of hyperbolicity in direction $\xx$ in G{\aa}rding's
1951 paper is that $q(t \xx + i \yy) \neq 0$ when $\yy \in \R^d$ 
and $t > t_0$.  This does not appear to me to be equivalent to 
Definition~\ref{def:hyperbolicity}, which appears in all later
work including~\cite{ABG}.  The two definitions specialize to
the same thing when $q$ is homogeneous.  Although the primary 
focus of these notes is on later uses outside of PDE's, I will
give indications of the proof.  Not only does this satisfy curiosity
and a historical sense, but it promotes the goal of using some
of the geometric understanding when dealing with contemporary,
more algebraic problems.  

To keep things simple, assume $q$ is homogeneous.  Examples
will be discussed shortly, but for now, a good mental picture
is to keep in mind the example $d=3$ and $q(x,y,z) = 
z^2 - x^2 - y^2$.  Even though $\R^d$ and its dual are isomorphic, 
I find it useful to classify a vector as belonging to the dual space
$(\R^d)^*$ rather than $\R^d$ if it plays the role of a linear 
functional on the original space.  For example the frequency 
$\xi$ of a wave $\xx \mapsto \exp (i \xi \cdot \xx)$ is thought 
of as living in $(\R^d)^*$.

\begin{figure}[ht]
\centering
\includegraphics[scale=.30]{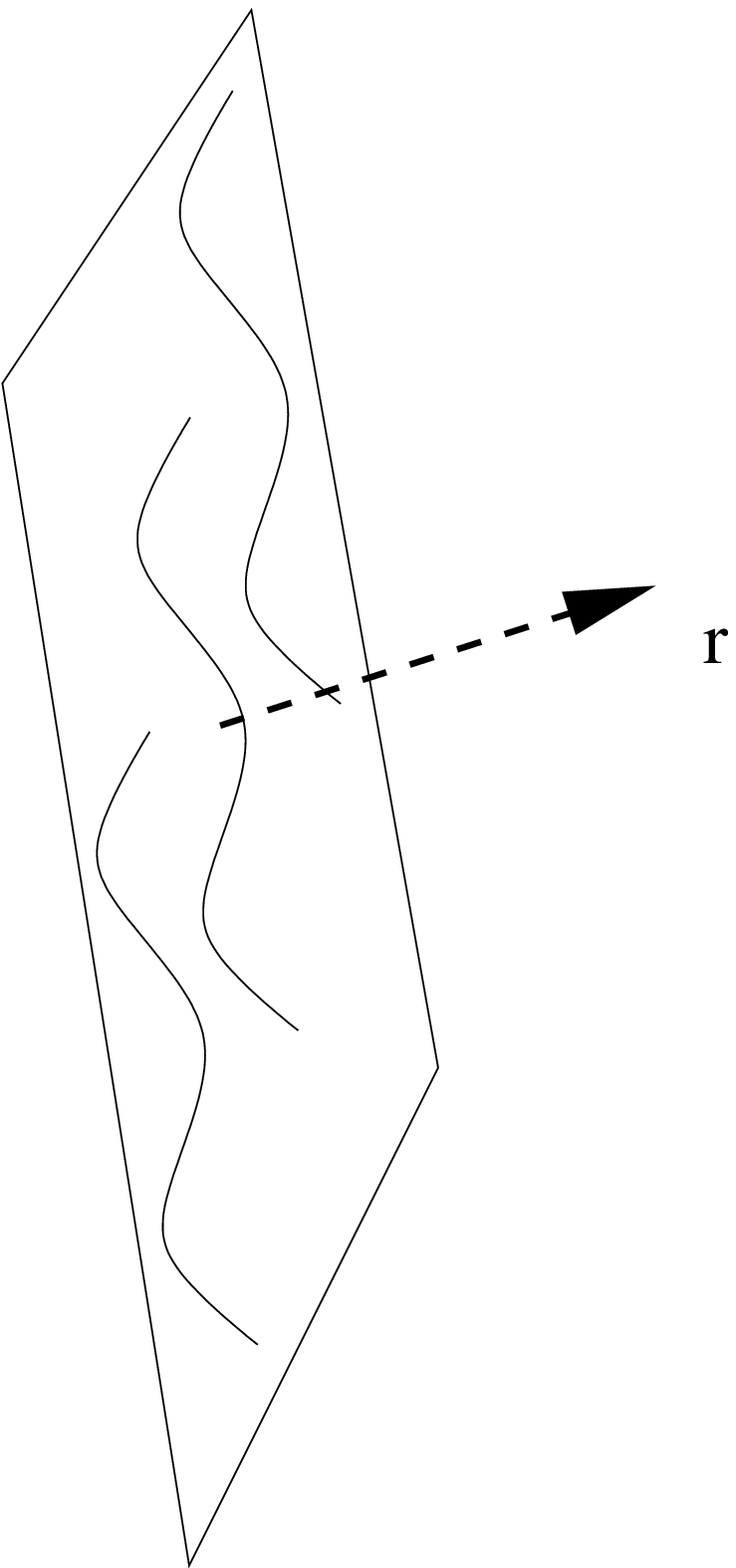}
\caption{$f_\rr$ when $\rr$ is real}
\label{fig:wave}
\end{figure}

\noindent{\sc Forward direction:}  For any $\xi \in (\R^d)^*$, 
let $f_\xi$ denote the function $\xx \mapsto e^{i \xi \cdot \xx}$
mapping $\R^d$ to $\C$.  Of course $(\R^d)^*$ is isomorphic to $\R^d$
but we call it $(\R^d)^*$ to remind ourselves that $\xi$ runs
over frequency space.  When $\rr$ is real, the function $f_\rr$ 
is a sinusoidal wave and is bounded, as in Figure~\ref{fig:wave}.
Let us assume without loss of generality 
that $\xi = (0 , \ldots , 0 , 1)$.  If $q(\xi) = 0$ then $f_\xi$ 
is a solution to $D[q] (f) = 0$.  If $q$ is not hyperbolic in 
direction $\xi$ then there is some $(r_1 , \ldots , r_{d-1})$ 
such that the $d$ solutions to $q(r_1 , \ldots , r_{d-1} , t) = 0$
include at least one conjugate pair of values that are not real.
Denoting by $r_d$ the one with the negative imaginary part, we
see that $f_\rr (\xx)$ grows exponentially as $x_d$ increases.
By homogeneity, $f_{\lambda \rr}$ grows at $\lambda$ times the
exponential rate.  Choose $c_\lambda$ so that $c_\lambda f_{\lambda \rr}
(0 , \ldots , 0 , 1) = 1$.  Then $f_{n \rr} \to 0$ on the
hyperplane $H_\xi$ in the topology of uniform convergence
(due to periodicity, no compact set restriction is required), 
while $f_{n \rr} \equiv 1$ at $(0 , \ldots , 0 , 1)$, proving 
the contrapositive, namely that lack of hyperbolicity implies 
lack of stability.
$\Cox$

\noindent{\sc Backward direction, handwaving argument:}  
For fixed $(r_1 , \ldots , r_{d-1})$, the set $\{ f_\rr : 
\rr \in (\R^d)^* \}$ of solutions to $D[q] f = 0$ spans 
a vector space $V(r_1 , \ldots r_{d-1})$ of dimension $k$.  
Since we are waving our hands, we have assumed \Em{strong} 
hyperbolicity, namely that the $d$ values of $r_d$ are distinct.
Any integral $\int c (r_1 , \ldots , r_{d-1}) d r_1 \cdots d r_{d-1}$
is also a solution, where $c(\cdot)$ is a section of $V(\cdot)$.
The handwaving part is that this gives all solutions of $D[q] f = 0$.
Assuming this, we write a generic solution $f$ as such an integral.
Each $f_\xi$ evolves unitarily, so more handwaving along the
lines of a Parseval relation shows that if $f$ and its first
$k$ derivatives are small at time zero then the time~$t$ value
is small as well.  
$\Cox$

Removing all the handwaving takes considerable work, beyond our 
scope here.  The key is the construction of the Riesz kernel,
which is interesting enough to merit a brief digression.
Suppose that $q$ is $m$-homogeneous and hyperbolic in direction 
$\xx$ and let $\alpha$ is any complex number.  Define the Riesz 
kernel $Q_\alpha$ by
\begin{equation} \label{eq:riesz kernel}
Q_\alpha (\yy) := (2 \pi)^{-d} \int_{\zz \in \xx + i \R^d}
   q(\zz)^{-\alpha} \exp (\yy \cdot \zz) \, d\zz \, .
\end{equation}

\begin{figure}[ht]
\centering
\includegraphics[scale=.30]{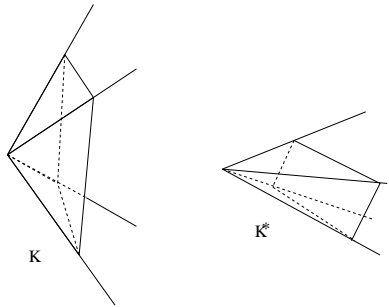}
\caption{The cone $K$ and its dual $K^*$.}
\label{fig:dual}
\end{figure}

\begin{figure}[ht]
\centering
\includegraphics[scale=.30]{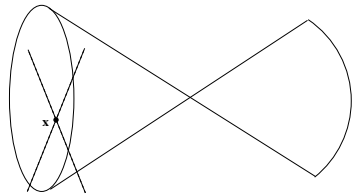}
\caption{Chain of integration in the definition of $Q_\alpha$}
\label{fig:imfiber}
\end{figure}

We have not yet defined cones of hyperbolicity, but let us
nonetheless try to visualize why the Riesz kernel is well
defined and supported on the dual cone.  To explain the 
terminology, the convex dual of a cone $K \subseteq \R^d$
is the cone $K^* \subseteq (\R^d)^*$ consisting of vectors $\yy$
such that $\ww \cdot \yy \geq 0$ for all $\ww \in K$.  Duality
maps fat cones to skinny cones and vice versa (see Figure~\ref{fig:dual}).
The connected components of the open set $\{ \ww \in \R^d : 
q(\ww) \neq 0 \}$ are cones.  If convex they have duals.
Let $C$ be the component containing $\xx$.  We will see later
that hyperbolicity guarantees this is one of the convex ones.
Hyperbolicity in direction $\xx$ guarantees that $q$ is 
nonvanishing on the domain of integration $\xx + i \R^d$.  
This is illustrated in Figure~\ref{fig:imfiber}, where the
dashed lines signify that the contour of integration is
varying in imaginary directions only and therefore does not
hit the surface $\{ q = 0 \}$.
The magnitude of the exponential term is constant, 
whence the integral is convergent at infinity when 
$m \, \Real \{ \alpha \} > 1 + d$.  Because the chain
of integration avoids $\{ q = 0 \}$, integrability is assured
as is a well defined branch of $q^{-\alpha}$.  It follows easily 
from Cauchy's Theorem that the Riesz kernel is independent of the
choice of $\xx$ within the component $C$.   Furthermore, if 
$\yy \in (\R^d)^*$ is a vector for which some $\xx \in C$ has
$\xx \cdot \yy < 0$, then it is clear that the integral defining
$Q_\alpha (\yy)$ can be deformed so as to be arbitrarily small,
hence is equal to zero.  In other words, $Q_\alpha$ is supported
on the set of $\yy$ for which there is no such $\xx$, which is
precisely the dual cone $C^*$.
$\Cox$

To complete the sketchy proof of the backward direction,
convolving with the Riesz kernel gives an operator $I_q$ 
inverted by $D[q]$.  This constuction is stable under perturbations.  
If $D[q] f = 0$, one shows that $f$ is equal to $(I - I_q D[q]) \tilde{f}$, 
where $\tilde{f}$ is defined stably in terms of $f$ and 
its first $d-1$ normal derivatives on $H_\xi$.  Small
changes in the initial conditions thus give rise to small 
changes in $\tilde{f}$, which give rise to small changes 
in $f$, proving the theorem.

The Riesz kernels with parameters $\alpha$ and $\alpha - 1$
are related by a differential identity as long as $\alpha \neq 1$.
When $\alpha = 1$, the Riesz kernel is called the \Em{fundamental 
solution} by reason of a result of G{\aa}rding.  This is explained
more fully in Theorem~\ref{th:delta} below, but the short
statement, found for instance in~\cite[Theorem~2.2]{guler},
is that the solution of 
$$D[q] f (\xx) = \delta (\xx)$$
exists, is unique, and is supported on the dual cone, $C^*$. 

\subsection{Homogeneous hyperbolic polynomials}

Many of the properties of hyperbolic polynomials are easier 
to state, prove and understand in the homogeneous case.  
Throughout this section, therefore, we deal only with 
homogeneous polynomials.  We will use $p$ rather than $q$ 
as a visual cue when speaking about homogeneous polynomials.  
Hyperbolicity is easily seen to be equivalent to the equivalent 
to the following ``real root'' property.
\begin{pr} \label{pr:homogeneous hyp}
The homogeneous polynomial $p$ is hyperbolic in direction $\xx$ if
and only if for any $\yy \in \R^d$ the univariate polynomial 
$t \mapsto p(\yy + t \xx)$ has only real roots.
\end{pr}

\noindent{\sc Proof:}
Because $p$ is homogeneous, when $\lambda \neq 0$,
we have $p(\lambda \zz) = 0$ if and only if $p(\zz) = 0$.  With
$\lambda = i s$ for some real nonzero value of $s$
we see that $p(\xx + i \yy) \neq 0$ for all $\yy \in \R^d$ 
(the definition of hyperbolicity) is equivalent
to $p(\yy + i s \xx) \neq 0$ for all $\yy \in \R^d$ and all
nonzero real $s$.  This is equivalent to $p(\yy + t \xx + i s \xx) \neq 0$
for all $\yy \in \R^d$ and real $s,t$ with $s \neq 0$.
Writing $z = t + i s$, this is equivalent to $p(\yy + z \xx) \neq 0$ 
when $z$ is not real.  
$\Cox$

Denote by $t_k (\xx , \yy)$ the roots of $t \mapsto p(\yy + t \xx)$.  
Then
\begin{equation} \label{eq:canonical}
p(\yy + t \xx) = p(\xx) \prod_{k=1}^m [t - t_k (\xx , \yy)] \, .
\end{equation} 
By Proposition~\ref{pr:homogeneous hyp}, when $p$ is hyperbolic in 
direction $\xi$, all values $t_k (\xx , \yy)$ are real, hence,
setting $t=0$, we see that $p/p(\xx)$ is a real polynomial.  We
see that little generality is lost in restricting our discussion
of homogeneous polynomials to those with real coefficients.

The two nontrivial examples of hyperbolic polynomials given
in G{\aa}rding's original paper, namely Lorentzian quadratics
and the determinant, provide much of the intuition as to
the meaning of hyperbolicity.  We now discuss these.  
A preliminary observation is that, trivially, all real 
homogeneous polynomials of degree one are hyperbolic in 
all directions in which they do not vanish.  Next, observe that 
hyperbolicity of the homogeneous polynomial $p$ in direction $\xi$ 
is preserved when $p$ and $\xi$ are transformed by the same invertible
real linear map.  This allows us to classify all nondegenerate 
quadratics: it suffices to consider all polynomials of the form 
$\sum_{j=1}^d \pm x_j^2$, hyperbolicity evidently being
determined by signature.  A necessary and sufficient condition 
for hyperbolicity is Lorentzian signature, that is precisely
one sign different from the others.

\begin{example}[Lorentzian quadratics] \label{eg:lorentz}
Let $p(\xx) := x_1^2 - \sum_{j=2}^d x_j^2$.  Real vectors $\xi$
may be classified as time-like, light-like or space-like according
to whether $p(\xi)$ is respectively positive, zero or negative.
The time-like vectors form the open convex cone $\{ \xx : (x_2/x_1)^2 +
\cdots + (x_d/x_1)^2 < 1 \}$.  Fix any time-like vector $\xi$.
If $\eta = \lambda \xi$ then the line $t \mapsto \eta + t \xi$ contains
a doubled root of $p$ at the origin.  For any other real $\eta$,
the line $\eta + t \xi$ intersects the hyperplane $H_\xi$ orthogonal 
to $\xi$ at some point other than the origin and $p$ takes a negative
value there.  On the other hand, the quadratic $p(\eta + t \xi)$ 
has positive leading term, so goes to $+\infty$ at $t = \pm \infty$.
Hence the line intersects the zero set of $p$ twice.  The degree
of $p$ is two, hence for any $\eta$, the polynomial 
$t \mapsto p(\eta + t \xi)$ has all real roots.  Hence $p$ is 
hyperbolic.  Replacing $p$ by $-p$ does not affect hyperbolicity.
Hence homogeneous quadratics with signature~1 or $d-1$ are
hyperbolic in timelike directions. 
\begin{figure}[ht]
\centering
\includegraphics[scale=.40]{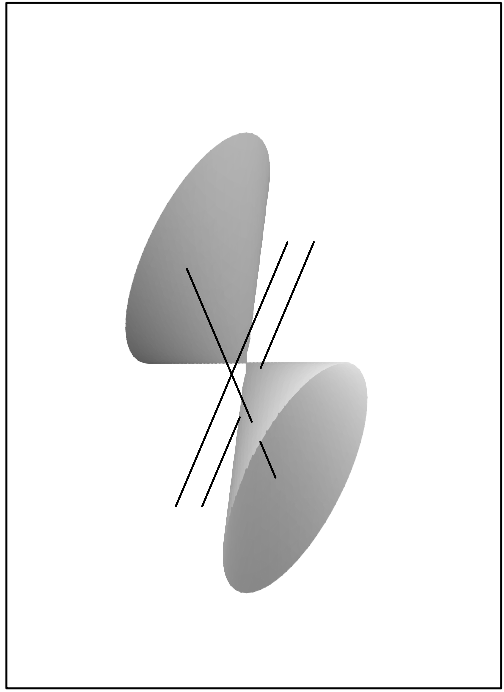}
\caption{This Lorentzian quadratic is hyperbolic in the $x$-direction
but not the $y$-direction.}
\label{fig:lorentz}
\end{figure}
\end{example}
For the converse, elliptic quadratics have lines in every direction 
with no real roots, hence are obviously not hyperbolic.  If $d > 4$ 
and there are at least two positive and two negative directions,
the any line on which $p$ has a minimum may be translated in a 
positive direction not on the line so that the minimum becomes 
positive, and similarly any line on which $p$ has a maximum may be
translated so that the maximum becomes negative.  We conclude there
are no directions of hyperbolicity.

The other classical example is as follows.
\begin{example}[Determinants of Hermitian matrices] 
\label{eg:hermitian}
An $n \times n$ matrix $M$ is Hermitian if $M_{ij} = \overline{M_{ji}}$.
The space of Hermitian matrices of size $n$ is parametrized by $n^2$ 
real parameters, these being the $\binom{n}{2}$ arbitrary
complex numbers $\{ M_{ij} : 1 \leq i < j \leq n \}$ and the
$n$ real numbers $\{ M_{ii} : 1 \leq i \leq n \}$.  Let $p$
be the homogeneous polynomial on $\R^{n^2}$ defined by the 
determinant as a function of these $n^2$ real numbers.  
Hyperbolicity of $p$ in the direction of any positive
definite matrix $A$ is equivalent to the well known fact
that the zeros of $M + zA$ are all real if $M$ is Hermitian
and $A$ is positive definite.  To see this, let $z = a + bi$
with $b > 0$ and write 
$$\Det (M + zA) = \Det (bA) \Det \left [ 
   (bA)^{-1/2} (a A + M) (bA)^{-1/2} + i I \right ] \, ,$$
where $(bA)^{1/2}$ is the Hermitian positive definite square root
of $bA$.  If this quantity is equal to zero then $-i$ is
an eigenvalue of $(bA)^{-1/2} (a A + M) (bA)^{-1/2}$
which is impossible because this matrix is Hermitian.
\end{example}

There are several ways to construct new hyperbolic functions
from old ones.  Note that the first does not require homogeneity.
\begin{pr}[products and polarization] \label{pr:new hyp}
~\\[-0.2in]
\begin{enumerate}[(i)]
\item Let $q_1$ and $q_2$ be hyperbolic with respect to $\xx$.
Then $q_1 q_2$ is also hyperbolic with respect to $\xx$.
\item Let $p$ be homogeneous of degree $m$ and hyperbolic with 
respect to $\xx$.  Let $p_0 , \ldots , p_m$ be the coefficients
of $p(\yy + t \xx)$ as a polynomial in $t$:
\begin{equation} \label{eq:expand}
p(\yy + t \xx) = \sum_{k=0}^m p_k (\yy) t^k \, .
\end{equation}
Then $p_k$ is hyperbolic with respect to $\xx$ for all $0 \leq k \leq m$.
\end{enumerate}
\end{pr}

\noindent{\sc Proof:} The fact for products is immediate from the 
definition.  For polarization, we first observe that hyperbolicity
of $p$ in direction $\xx$ implies hyperbolicity of the directional
derivative $D[\xx] [p]$ in direction $\xx$.  To see this, observe 
that $f(t)$ having only real roots implies $f'(t)$ has only real roots.  
Setting $f (t) = p (\yy + t \xx)$, we see that $f'(t)$, namely 
$t \mapsto D_\xx (p) (\yy + t \xx)$, has only real roots for any 
$\yy \in \R^d$, which is hyperbolicity of $D_\xx (p)$ in direction $\xx$.    

Now consider the expansion~\eqref{eq:expand} of $p(\yy + t \xx)$
into homogeneous parts.  The parts are given by
$$p_k (\yy) = \frac{1}{k!} \, \left. \left ( \frac{d}{dt} 
    \right )^k \right |_{t=0}
   p(\yy + t \xx) = \left ( D_\xx \right )^k (p) (\yy) \, .$$ 
Hyperbolicity in direction $\xx$ is stable under $D_\xx$, which 
shows that $p_k$ is hyperbolic in direction $\xx$.
$\Cox$

\subsection{Cones of hyperbolicity for homogeneous polynomials}

Hyperbolic polynomials have associated with them certain convex cones.
In the homogeneous case the definition is simple and is
contained in the following proposition.  This was proved 
first in~\cite{garding-hyperbolic} and reproduced many times;
the proof below follows~\cite{guler}.
Let $p$ be a (complex) homogeneous hyperbolic polynomial, hyperbolic
in direction $\xi$.  Dividing by a real multiple of $p(\xi)$ we
may assume, cf.~\eqref{eq:canonical}, that $p$ is real and $p(\xi) = 1$.
\begin{pr} \label{pr:convex}
Let $K(p,\xi)$ denote the connected component 
of the set $\R^d \setminus \{ \xx : p(\xx) = 0 \}$ that contains $\xi$.
\begin{enumerate}[(i)]
\item $p$ is hyperbolic in direction $\xx$ for every $\xx \in K(p,\xi)$. 
\item The set $K(p , \xi)$ is an open convex cone; we call this a
\Em{cone of hyperbolicity} for $p$.
\item $K(p , \xi)$ is equal to the set $K$ of vectors $\xx$ for which
all roots of $t \mapsto p(\xx + t \xi)$ are real and negative,
which (by hyperbolicity of $p$ in direction $\xi$) is the same as 
the set of vectors $\xx$ for which no root of this polynomial is real 
and nonnegative.
\end{enumerate}
\end{pr}

\noindent{\sc Proof:} By continuity of the roots of a polynomial
with respect to its coefficients, the set $K$ is open.  Also 
$\xi \in K$ because $p(\xi + t \xi) = (1+t)^m p(\xi)$ where $m$
is the degree of the homogeneous polynomial $p$.  If $\xx$ is
in the closure $ \overline{K}$ then, again by continuity of 
the roots, $p(\xx + t \xi) \neq 0$ for $t > 0$; in particular,
$p(\xx) \neq 0$ implies $\xx \in K$, which means that 
$K$ is closed in $\R^d \setminus \{ p = 0 \}$.  Also
$K$ is connected: by definition of $K$, if $\xx \in K$
then $\xx + t \xi \in K$ for $t > 0$, and by homogeneity $s \xx 
+ t \xi \in K$ for $s,t > 0$, and sending $s$ to zero proves
that $K$ is star-convex at $\xi$ (hence connected).
Being connected and both closed and open in $\R^d \setminus \{ p=0 \}$,
$K$ is a component of $\R^d \setminus \{ p=0 \}$ and is thus
equal to $K(p , \xi)$.  

Next we check that for $\vv \in K$, $\xx$ real, and $s$ and $t$ complex,
\begin{equation} \label{eq:imag}
p(\xx + t \vv + s \xi) \neq 0 \mbox{ whenever } 
   \Im (s) \leq 0 , \Im (t) \leq 0 , 
   \Im (s) + \Im (t) < 0 \, .
\end{equation}
First, suppose $\Im (s) < 0$.  If $u \geq 1$ then the polynomial 
\begin{equation} \label{eq:poly}
t \mapsto p \left [ \frac{\xx}{u} + t \vv + \left ( \frac{s}{u} 
   + i \left ( \frac{1}{u} - 1 \right ) \right ) \xi \right ] 
\end{equation}
is nonvanishing for real $t$ because this is $p$ evaluated at the sum 
of a real vector and a complex multiple of $\xi$ with nonzero imaginary 
part.  As $u \to \infty$, the number of roots of~\eqref{eq:poly} 
in the lower half-plane remains constant.
The limit at $u = \infty$ is $t \mapsto p(t \vv - i \xi)$.  This
has roots in the upper half-plane because its roots are $-i$
divided by the roots of $t \mapsto p(\vv + t \xi)$, the
latter of which we have seen to be negative real.  We conclude
that for $u=1$ there are no roots in the lower half-plane,
in other words,
\begin{equation} \label{eq:one way}
p (\xx + t \vv + s \xi) \neq 0 \mbox{ whenever }
   \Im (s) < 0 , \Im (t) \leq 0 \, .
\end{equation} 
To see that this also holds for $\Im (s) = 0$ and $\Im (t) < 0$,
note that $K$ is open, so $\vv - \ee \xi \in K$ for some $\ee > 0$,
whence 
$$p(\xx + t \vv + s \xi) = p(\xx + t (\vv - \ee \xi) + (s + \ee t) \xi)$$
which is nonzero by~\eqref{eq:one way}.  This completes the verification
of~\eqref{eq:imag}.

Setting $s=0$ in~\eqref{eq:imag} shows that any root of
$t \mapsto p(\xx + t \vv)$ satisfies $\Im (t) \geq 0$.  
Because $p$ is real, complex conjugation may be applied to the
entire argument, showing that also $\Im (t) \leq 0$, and
hence that all roots of $t \mapsto p(\xx + t \vv)$ are real.
The vector $\vv$ was chosen arbitrarily in $K$.  We conclude
that $p$ is hyperbolic with respect to every $\vv \in K$.
It follows that $K$ is star convex with respect to $\vv$ as
well, and since $\vv \in K$ is arbitrary, that $K$ is convex.
$\Cox$

Because the cones of hyperbolicity are characterized as
components of the nonzero set in real space, it follows
that the cones of hyperbolicity of $pq$ are pairwise
intersections of the cones of hyperbolicity of $p$ with 
the cones of hyperbolicity of $q$.  

\begin{figure}[ht]
\centering
\includegraphics[scale=.40]{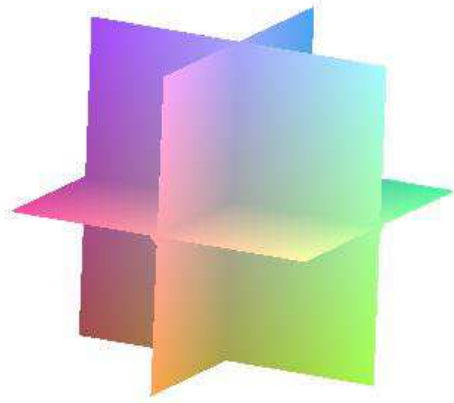}
\caption{The cones of hyperbolicity for $\prod_{j=1}^d z_j$ are the orthants}
\label{fig:planes}
\end{figure}

\begin{example}[coordinate planes] \label{eg:coords}
Each coordinate function $z_j$ is homogeneous of degree~1, 
therefore hyperbolic in every direction not contained in
the plane $\{ z_j = 0 \}$.  The cones of hyperbolicity
are the two half spaces bounded by this plane.  It follows 
that the product $\prod_{j=1}^d z_j$ is hyperbolic in every
direction in which no coordinate vanishes, and that the cones 
of hyperbolicity are the orthants.   This is also obvious from
Figure~\ref{fig:planes}.
\end{example}

Returning briefly to the origins of hyperbolicity theory in the
properties of wave-like partial differential equations, the 
following result was proved in~\cite[Theorem~4.1]{ABG}.
It says roughly that cones of hyperbolicity are propagation 
cones for solutions to wave-like equations.  Given any convex
cone $K \subseteq \R^d$ over the origin, let $K^*$ denote
the dual cone of vectors $\vv^* \in (\R^d)^*$ such that
$\langle \vv^* , \xx \rangle \leq 0$ for all $\xx \in K$.
Recall the linear partial differential operator $D[p]$
defined by substituting $\partial / \partial x_i$ for $x_i$ 
in $p$.  Let $E = E(p,\xx)$ be a solution to
\begin{equation} \label{eq:delta}
D[p] E = \delta_\xx
\end{equation}
where $\delta_\xx$ is a delta function at $\xx$.

\begin{thm} \label{th:delta}
Suppose the homogeneous polynomial $p$ is hyperbolic in
direction $\xi$.  Then there exists a unique distribution $E$
solving~\eqref{eq:delta} and its support is the cone $K(p , \xi)^*$
dual to the cone of hyperbolicity of $p$ containing $\xi$.
In fact $E$ is given as in inverse Fourier transform:
$$E(p,\xi) = (2 \pi)^{-d} \int_{\R^d} p_- (\yy)^{-1} 
   e^{i \langle \xx , \yy \rangle} \, d\yy$$
where $p_- (y)^{-1} := \lim_{t \downarrow 0} p(\yy - i t \nu)^{-1}$
as distributions and $\nu$ is any element of $K (p , \xi)$.
$\Cox$
\end{thm}

Although there is no time here to explore the properties
of hyperbolic polynomials as barrier functions, I will
mention this example briefly, as it has proved to be of
some importance to convex programming.

\begin{example}[self-concordant barrier functions]
Let $Q \subseteq \R^d$ be an open convex set.  A function
$F : Q \to \R$ is a \Em{logarithmically homogeneous 
self-concordant barrier function} for $Q$ if it is 
smooth, convex, and satisfies several properties:
\begin{eqnarray*}
F(\xx) \to \infty \mbox{ as } \xx \to \partial Q && \\[1ex]
\left | D^3 F(\xx) [u, u, u] \right | & \leq & 
   2 \left ( D^2 F (\xx) [u,u] \right )^{3/2} \\
|\left | DF (\xx) [u] \right |^2 & \leq & \theta D^2 F(\xx) [u,u] \\
F(t\xx) & = & F(\xx) - \theta \log t  \, .
\end{eqnarray*}
The interior point method for convex programming problems 
(see, e.g.~\cite{nesterov-nemirovskii}) is based on
finding such a function for a given region, $Q$.  
There is a universal construction but it is not always
useful for computations.  Also, properties beyond those
satisfied by the universal construction are required for
long term stability of the interior point method.  Such
properties will depend on the region $Q$.  When $Q$ is the 
cone of hyperbolicity for a homogeneous polynomial $p$,
it turns out that the function $F(\xx) := - \log p(\xx)$
is a logarithmically homogeneous self-concordant barrier
function with a number of other useful properties.  These
are detailed in~\cite{guler}.  
\end{example}

\subsubsection*{Hyperbolicity without homogeneity}

Delving into the relation between hyperbolicity
in the homogeneous and non-homogeneous cases, 
(definitions both given in Definition~\ref{def:hyperbolicity}, we
begin by looking at the homogeneous parts of a polynomial, $f$.
Sorting the terms of $f$ by total degree, the highest degree 
part will be denoted $LT (f)$ for ``leading term(s)''.
The lowest degree is called the \Em{localization} and
will be discussed further in the next section.  If $f$
is hyperbolic, then both the leading term and the 
localization are hyperbolic as well (respectively
\cite[Lemma~3.20]{ABG} and~\cite[Lemma~3.42]{ABG}).
By definition, the cones of hyperbolicity of a non-homogeneous
polynomial $f$ are the cones of hyperbolicity of $LT(f)$.
When $f$ is hyperbolic, these cones still characterize
the directions of stable propagation of PDE's and
supports of solutions to these.

Because hyperbolicity is easier to understand in the
homogeneous case, G{\aa}rding looked for a converse
to these and came up with a criterion for homogeneous
polynomials called strong hyperbolicity.  This is most
naturally stated using the equivalent definition of 
hyperbolicity from Proposition~\ref{pr:homogeneous hyp}:
the roots of $p(\yy + t \xx)$ must not only be real but
also distinct (except when $\yy$ is a multiple of $\xx$,
when all the zeros coincide perforce).  We then have:
\begin{thm}
If $A$ is strongly hyperbolic then $f$ is hyperbolic for
any $f$ such that $LT(f) = A$.
\end{thm}
This is proved as~(3.9) of~\cite{ABG} and in fact is already
in~\cite{garding-hyperbolic} after the statement and proof 
of Lemma~2.5.  Necessary and sufficient conditions for
$f$ to be hyperbolic are given 
in~\cite[Theorem~12.4.6]{hormander-partial-differential2};
these were conjectured by G{\aa}rding and proved first by
Svensson~\cite{svensson}.

\pagebreak
\setcounter{equation}{0}
\section{Semi-continuity and Morse deformations} \label{sec:morse} 

\subsection{Localization} \label{ss:local}

Given a function $f$ analytic on a neighborhood of the origin,
its order of vanishing is the least total degree of a nonvanishing
term in its Taylor series.  The sum of all such terms is
called the homogeneous part of $f$ at the origin and denoted 
$\homog (f)$.  For any $\xx$ we let
$$\homog (f,\xx) := \homog (f (\xx + \cdot))$$ 
denote the homogeneous part of $f$ at $\xx$.  When $\xx \neq 0$,
Atiyah {\em et al.}~\cite{ABG} define this by taking the $t=0$ term of
$t^m f(t^{-1} \xx + \yy)$ as a function of $\yy$, where
$m - \deg (f)$, and they refer to this as the \Em{localization} 
of $f$ at $\xx$.

\begin{pr} \label{pr:semi 1} 
Let $p$ be any hyperbolic homogeneous polynomial, and let $m$ be
its degree.  Fix $\xx$ with $p(\xx) = 0$ and let $\pt := \homog (p , \xx)$
denote the leading homogeneous part of $p$ at $\xx$.  If
$p$ is hyperbolic in direction $\uu$ then $\pt$ is also
hyperbolic in direction $\uu$.  Consequently, if $B$ is
any cone of hyperbolicity for $p$ then there is some
cone of hyperbolicity for $\pt$ containing $B$.
\end{pr}
\noindent{\sc Proof:} This follows from the conclusion~(3.45)
of~\cite[Lemma~3.42]{ABG}.
Because the development there is long and complicated, we
give here a short, self-contained proof, provided by 
J. Borcea~\cite[Proposition~2.8]{BP-cones}.  If $Q$ is a 
polynomial whose degree at zero is $k$, we may recover its 
leading homogeneous part $\homog(Q)$ by
$$\homog (Q) (\yy) = \lim_{\lambda \to \infty}
   \lambda^k Q (\lambda^{-1} \yy) \, .$$
The limit is uniform as $\yy$ varies over compact sets.
Indeed, monomials of degree $k$ are invariant under the scaling
on the right-hand side, while monomials of degree $k+j$ scale
by $\lambda^{-j}$, uniformly over compact sets.

Apply this with $Q (\cdot) = p(\xx + \cdot)$ and $\yy + t \uu$
in place of $\yy$ to see that for fixed $\xx, \yy$ and $\uu$,
$$\pt (\yy + t \uu) = \lim_{\lambda \to \infty}
   \lambda^{k} p(\xx + \lambda^{-1}(\yy + t \uu))$$
uniformly as $t$ varies over compact sub-intervals of $\R$.
Because $p$ is hyperbolic in direction $\uu$, for any fixed
$\lambda$, all the zeros of this polynomial in $t$ are real.
Hurwitz' theorem on the continuity of zeros~\cite[Corollary~2.6]{conway}
says that a limit, uniformly on bounded intervals, of polynomials
having all real zeros will either have all real zeros or vanish
identically.  The limit $\pt (\yy + t \uu)$ has degree $k \geq 1$;
it does not vanish identically and therefore it has all real zeros.
This shows $\pt$ to be hyperbolic in direction $\uu$.
$\Cox$

\begin{defn}[family of cones in the homogeneous case]
\label{def:KAB}
Let $p$ be a hyperbolic homogeneous polynomial and let $B$
be a cone of hyperbolicity for $p$.  If $p(\xx) = 0$, define
$$\cK^{p,B} (\xx)$$
to be the cone of hyperbolicity of $\homog (p,\xx)$ containing
$B$, whose existence we have just proved.  If $p(\xx) \neq 0$
we define $\cK^{p,B} (\xx)$ to be all of $\R^d$.
\end{defn}

\begin{example}[strongly hyperbolic functions] \label{eg:strong}
As mentioned in passing in Section~\ref{sec:hyperbolicity}, 
a homogeneous polynomial $p$ is said to be strongly
hyperbolic if the roots of $p(\yy + t \xx)$, in addition
to being real, are distinct.  Equivalently, the projective
variety $\sing := \{ \xx : q(\xx) = 0 \}$ is smooth.
In this case, for each nonzero $\xx$, the localization 
$\homog (p , \xx)$ has degree one.  The cones of hyperbolicity
$K^{p , B} (\xx)$ are halfspaces whose common tangent hyperplane 
is tangent to $B$ at $\xx$; the convexity of $B$ implies
that the tangent hyperplane is a support hyperplane to $B$
at $\xx$, and we see that indeed $K^{p,B} (\xx)$ contains $B$.
\end{example}

Suppose $\xx_n \to \xx$.  It is not in general true that
$\homog (f , \xx_n) \to \homog (f , \xx)$.  However, it
is true that 
$$\homog (f , \xx_n) = \homog (\homog (f , \xx) , \xx_n)$$
for $n$ sufficiently large.  This implies that if $p$ is
homogeneous and hyperbolic with cone of hyperbolicity $B$,
then $K^{p,B} (\xx)$ is semi-continuous in $\xx$:
\begin{equation} \label{eq:semi 1}
K^{p,B} (\xx) \subseteq \liminf K^{p,B} (\xx_n) 
   \mbox{ as } \xx_n \to \xx \, .
\end{equation}
This is proved in~\cite[Lemma~5.9]{ABG}.  We will want a
version of this valid for polynomials that are not
necessarily homogeneous.  If $q$ is hyperbolic
but not homogeneous, the cone $K^{q,B} (\xx)$ has not
been defined.  To do this, at least for some points $\xx$,
we need the notion of the amoeba of $q$.

\subsection{Amoeba boundaries} \label{ss:amoeba}

\begin{figure}[ht]
\centering
\includegraphics[scale=.20]{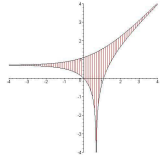}
\caption{The amoeba of the polynoimal $2 - x - y$}
\label{fig:amoeba}
\end{figure}

\begin{defn}[amoeba] \label{def:amoeba}
The \Em{amoeba} of any polynomial $q$ is defined
to be the image of the zero set $\sing := \{ \zz : q = 0 \}$ of
$q$ under the coordinatewise log modulus map $\Log (\zz) := 
(\log |z_1| , \ldots , \log |z_d|)$.  We denote this image
by $\amoeba (q)$.  
\end{defn}
The connected components of $\R^d \setminus \amoeba (q)$
are convex sets and are in one to one correspondence with the
Laurent series expansions of $F$, with each expansion converging
on precisely the set $\{ \exp (\xx + i \yy) : \yy \in B \}$
for some component $B$ of $\R^d \setminus \amoeba (q)$.  This
is well known and is presented, for instance, in Chapter~6
of~\cite{GKZ}.  For those who have not seen an amoeba before,
one of the simplest nontrivial amoebas is shown in 
Figure~\ref{fig:amoeba}.

The following result, proved 
in~\cite[Proposition~2.12]{BP-cones}, defines a family of
cones via localizations on the boundary of an amoeba.
\begin{pr} \label{pr:amoeba hyp}
Let $q$ be any polynomial and let $B$ be a component of the complement
of $\amoeba (q)$.  Denote the boundary of $B$ by $\partial B$.
Fix $\xx \in \partial B$ and let $f = q \circ \exp$
so that $f$ vanishes at some point $\xx + i \yy$.  Let $f_\yy := 
\homog (f , \xx + i \yy)$.  Then each $f_\yy$ is hyperbolic (meaning
hyperbolic in at least one direction) and one of its cones of 
hyperbolicity contains $B$.  We denote this cone
by $K^{q,B} (\yy)$.  (The point $\xx$ is considered fixed and
is suppressed from the notation.) 
$\Cox$
\end{pr}

One may extend semi-continuity for localizations beyond the
homogeneous case, to the families in Proposition~\ref{pr:amoeba hyp}.
\begin{thm}[\protect{\cite[Corollary~2.15]{BP-cones}}]
\label{th:semi 2}
If $B$ is a component of $\amoeba(q)^c$ and $\xx \in \partial B$,
then the family of cones $K^{q,B} (\yy)$ defined in 
Proposition~\ref{pr:amoeba hyp} satisfies
\begin{equation} \label{eq:semi 2}
K^{q,B} (\yy) \subseteq \liminf K^{q,B} (\yy_n) \mbox{ whenever }
   \yy_n \to \yy \, .
\end{equation} 
$\Cox$
\end{thm}

\begin{figure}[ht]
\centering
\includegraphics[scale=.40]{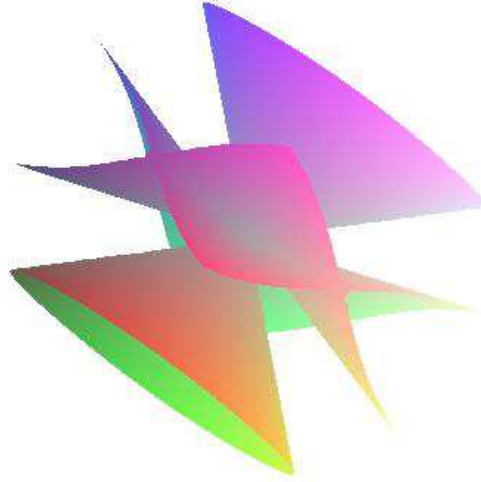}
\caption{Another hyperbolic homogeneous polynomial}
\label{fig:fortress}
\end{figure}

The fact that localizations of {\em any} polynomial to 
points on the amoeba boundary are hyperbolic allows us
to give numerous examples fo hyperbolic polynomials 
beyond the classic ones: planes, quadrics, and the determinant
function.  For example, the homogeneous polynomial variety
shown in Figure~\ref{fig:fortress} is a localization of
the famous so-called fortress generating function denominator.

\begin{figure}[ht]
\centering
\includegraphics[scale=.40]{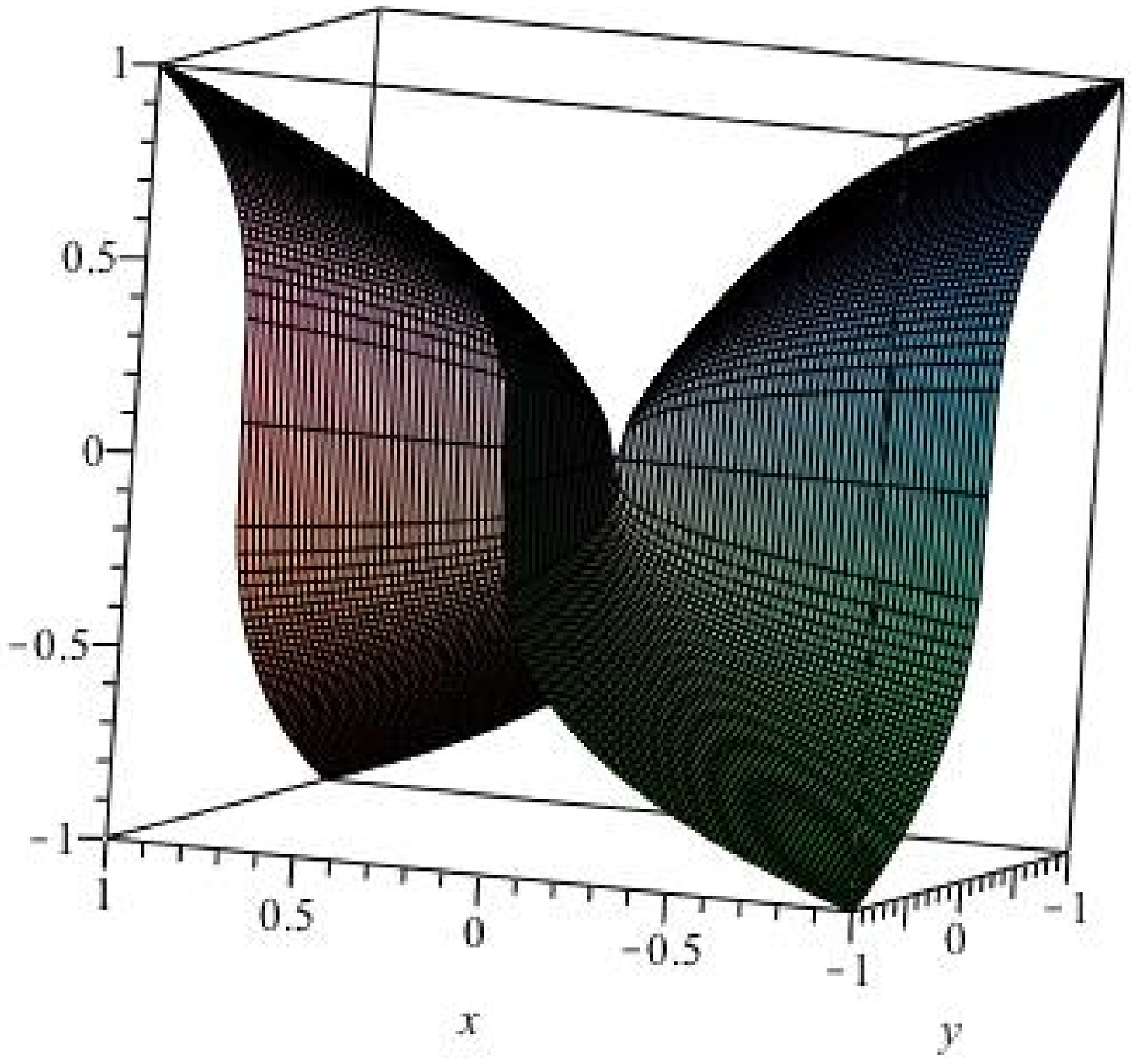}
\caption{The sheared cubic $xy+z^3$}
\label{fig:twisted}
\end{figure}

It is instructive to see what role hyperbolicity plays in 
Theorem~\ref{th:semi 2} by considering the counterexample
$f(x,y,z) = xy + z^3$ as in Figure~\ref{fig:twisted}.  
As $\yy$ varies over a neighborhood $U$ of the origin, 
is it possible to choose open convex cones 
$\{ K(\yy) : \yy \in U \}$ over $\yy$ in such a way that 
$K(\yy)$ varies semi-continuously with $\yy$ and each $K(\yy)$ 
is, locally, a subset of $\{ \yy : f (\yy) \neq 0 \}$?  
The points $\yy = (0,y,0)$ are forced to choose whether $K (\yy)$ 
contains points with positive $\xx$ components or negative $\xx$
components.  One of these violates semi-continuity with cones
$K(0,y,\ee)$ as $\ee \downarrow 0$ while the other violates
semi-continuity with cones $K(0,y,-\ee)$ as $-\ee \uparrow 0$.
The key to avoiding this in Theorem~\ref{th:semi 2} is that
$f_\yy$ must be hyperbolic, and this cannot happen with hyperbolic
functions.

\subsection{Morse deformations}

Suppose that $K (\yy)$ is any family of convex cones in $\R^d$ varying
semi-conintuously with $\yy$ in the sense of~\eqref{eq:semi 2}.
Let $\rr$ be any vector such that each $K(\yy)$ contains a
vector $\vv$ with $\rr \cdot \vv > 0$.  
\begin{pr} \label{pr:VF}
Over any compact set $S$ a section $\{ \vv (\yy) : \yy \in S \}$
may be chosen continuously with $\vv (\yy) > 0$ for all $\yy \in S$.
\end{pr}

\noindent{\sc Proof:} Given $\yy$, we may choose $\vv$ in the 
interior of $K(\yy)$ with $\vv \cdot \rr > 0$, whence by
semi-continuity, $\vv \in K(\yy')$ for all $\yy'$ in some
neighborhood of $\yy$ which we denote $\nbd (\yy , \vv)$. 
We may cover the compact set $S$ with finitely many of these
neighborhoods, $\{ \nbd (\yy^{(k)} , \vv^{(k)}) : 1 \leq k \leq m \}$.  
Let $\{ \psi^{(k)} : 1 \leq k \leq m \}$ be a partition of unity 
subordinate to this cover and for any $\yy \in S$ define
$$\vv (\yy) := \sum_{k=1}^m \psi^{(k)} (\yy) \vv^{(k)} \, .$$
It is clear that $\rr \cdot \vv (\yy) > 0$ for all $\yy \in S$.
For each $\yy \in S$, if $\psi^{(k)} (\yy) > 0$ then
$\vv^{(k)} \in K(\yy)$.  By convexity of $K(\yy)$, it follows
that $\vv (\yy) \in K(\yy)$.
$\Cox$

Two applications of this are as follows.  Suppose $p$ is a homogeneous
hyperbolic function and let $B$ be a cone of hyperbolicity for $p$.
Due to homogeneity, $K^{p,B} (\lambda \yy) = K^{p,B} (\yy)$, whence
the family may be described as $\{ K^{p,B} (\yhat) : \yhat \in S^{d-1} \}$.
Suppose $\rr$ is such that each $K (\yy)$ contains a $\vv$ with 
$\rr \cdot \vv > 0$.
Because $S^{d-1}$ is compact, an application of Proposition~\ref{pr:VF}
yields a continuous section $\{ \vv (\yy) : \yy \in S^{d-1} \}$ with
$\vv(\yy) \in K^{p,B} (\yy)$ and $\rr \cdot \vv (\yy) > 0$.  Extending
to all of $\R^d \setminus \{ \zero \}$ by $\vv (\lambda \yy) 
:= \lambda \vv (\yy)$, we arrive at:
\begin{cor} \label{cor:VF 1}
There is a vector field $\vv (\cdot)$ on $\R^d$ which is 1-homogeneous,
vanishes only at $\zero$, is a section of $K^{p,B} (\cdot)$, and
for some $\ee > 0$ satisfies $\rr \cdot \vv (\yy) \geq \ee |\yy|$. 
$\Cox$
\end{cor}

\begin{example} \label{eg:strong hyp}
If the homogeneous polynomial $q$ is strongly hyperbolic then
its zero set is the cone over a smooth projective hypersurface
and each $K(\yy)$ is a halfspace (see Example~\ref{eg:strong}).
The condition that every $K(\yy)$ contain a vector $\vv$ with
$\rr \cdot \vv > 0$ is the same as requiring that $\rr$ not
be the outward normal to the bounding hyperplane of $K(\yy)$.
In other words, $\rr$ may not be on the boundary of the 
dual cone to $B$.
\end{example}

The second application is when $B$ is a component of $\amoeba (q)^c$.
In this case the family $\{ f_\yy \}$ is periodic with period
$2 \pi$ in each coordinate.  Accordingly, the cones 
$K(\yy) = K^{q,B} (\yy)$ from Proposition~\ref{pr:amoeba hyp}
are indexed by a compact set, namely the torus $\R^d / (2 \pi \Z)^d$. 
Applying Proposition~\ref{pr:VF} then yields:

\begin{cor} \label{cor:VF 2}
If each cone $K^{q,B} (\yy)$ in Proposition~\ref{pr:amoeba hyp}
contains a vector $\vv$ with $\rr \cdot \vv > 0$ then there is
a continuously varying section $\{ \vv (\yy) : \yy \in \R^d \}$
with $\vv (\yy) \in K^{q,B} (\yy)$ and $\rr \cdot \vv (\yy) > 0$
for all $\yy$.
$\Cox$
\end{cor}

With these vector fields in hand we have enough to carry out
the programs of~\cite{ABG} and~\cite{BP-cones}.  I will briefly
describe the former, then go into a little more detail on the 
latter.  Let $q$ be a homogeneous polynomial, strongly hyperbolic 
in the direction $\xx$.  We wish to compute its inverse Fourier
transform, which we have called the Riesz kernel, $Q_\alpha$.
We recall that its value at $\yy$ is given by
$$ Q_\alpha (\yy) := (2 \pi)^{-d} \int_{\zz \in \xx + i \R^d}
   q(\zz)^{-\alpha} \exp (\yy \cdot \zz) \, d\zz \, .  
   \eqno{\protect{\eqref{eq:riesz kernel}}} $$
We set $\rr = - \yy$ and assume $\rr$ not to be on the boundary
of $B^*$, the dual cone of the cone $B$ of hyperbolicity of $q$ 
that contains $\xx$.  As we have seen in Example~\ref{eg:strong hyp},
this guarantees the existence of the 1-homogeneous vector field
$\vv (\cdot)$.  The sets $K(\yy)$ are convex and each contains
both $\xx$ and $\vv (\yy)$.  Therefore, each contains the line
segment joining $\xx$ to $\vv (\yy)$.  Define a homotopy by
$$H(\uu , t) := i \uu + (1-t) \xx + t \vv (\uu) \, .$$
This deforms the domain of integration $\xx + i \R^d$ 
in~\eqref{eq:riesz kernel} to the cone $\CC := \{ i \uu + \vv (\uu) :
\uu \in \R^d \}$ while avoiding the set $\sing$ where $q$ vanishes.
Cauchy's theorem implies that deforming the contour does not change 
the integral and therefore that
$$Q_\alpha (\yy) = \int_{\CC} q(\zz)^{-\alpha} \exp (\yy \cdot \zz) 
   \, d\zz \, .$$
Omitting pages of detail and skipping to the punch line, this
representation allows us to factor the integral.  Integrating 
radially reduces the integral to an integral over the
\Em{leray cycle}, which is a $(d-1)$-dimensional homology class 
in the complement of $\{ q=0 \}$ in $\CP^{d-1}$.  The dimension
reduction allows a number of explicit computations which are carried
out in~\cite[Section~7]{ABG}. 

\subsection{Asymptotics of Taylor coefficients}

Let $F = P/Q = \sum_\rr a_\rr \ZZ^\rr$ be a Laurent series
expansion converging on the component $B$ of $\amoeba (Q)^c$.
The coefficients $\{ a_\rr \}$ of $F$ may be recovered via
Cauchy's formula
\begin{equation} \label{eq:cauchy}
a_\rr = (2 \pi i)^{-d} \int_T \ZZ^{-\rr} F(\ZZ) \frac{d\ZZ}{\ZZ}
\end{equation}
where $T$ is the torus $T(\xx) = \{ \exp (\xx + i \yy) : 
\yy \in \R^d / (2 \pi \Z)^d \}$ for some $\xx \in B$
and $d\ZZ / \ZZ := dz_1 \wedge \cdots \wedge dz_d / (z_1 \cdots z_d)$.
The motivation for the evaluation or asymptotic estimation of 
the coefficients $\{ a_\rr \}$ comes from analytic combinatorics,
where the primary object of study is the array $\{ a_\rr \}$,
which counts something of interest.  One constructs the generating
function $F (\ZZ) := \sum_\rr a_\rr \ZZ^\rr$ and hopes
to identify a closed form representation of $F$.  If $\{ a_\rr \}$
satisfy a recursion, then (depending on boundary values) this
will usually succeed.  A number of examples of combinatorial
interest are surveyed in~\cite{PW9}.  In most of the cases
surveyed there, the generating function $F$ is rational.  
In~\cite{raichev-wilson-safonov}, based on a result 
of~\cite{safonov}, it is shown how to embed any algebraic 
function of $d$ variables as a diagonal of a rational function 
of $(d+1)$ variables.  Thus the evaluation of Taylor coefficients
of rational functions solves the enumeration problem for any
array of numbers whose generating function is rational or
algebraic.

Exact evaluation of the Cauchy integral~\eqref{eq:cauchy} 
is not easy but some methods are known for evaluating it
asymptotically.  In the case where the pole variety $\sing :=
\{ \zz : Q(\zz) = 0 \}$ is smooth, a formula was given 
in~\cite{PW1}; a coordinate free version is given in~\cite{BBBP}. 
Normal self-intersection in $\sing$ can also be dealt with~\cite{PW2}.
A number of examples from combinatorics and statistical physics
have generating functions whose pole variety has a singularity
with nontrivial monodromy.  In the remainder of this section I will
explain how hyperbolicity and the resulting deformations allow
explicit asymptotic evaluation of some of these generating functions.

\begin{figure}[ht]
\centering
\vspace{0.6in} \includegraphics[scale=0.64, bb=110 -20 210 120]{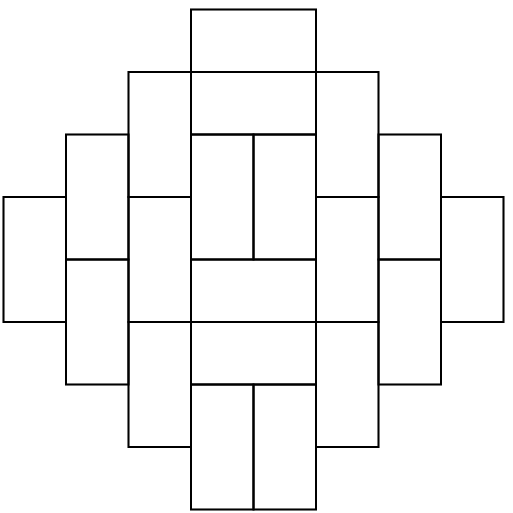} 
\hspace{0.4in}
\vspace{0.2in} \includegraphics[scale=0.20, bb=010 90 130 190]{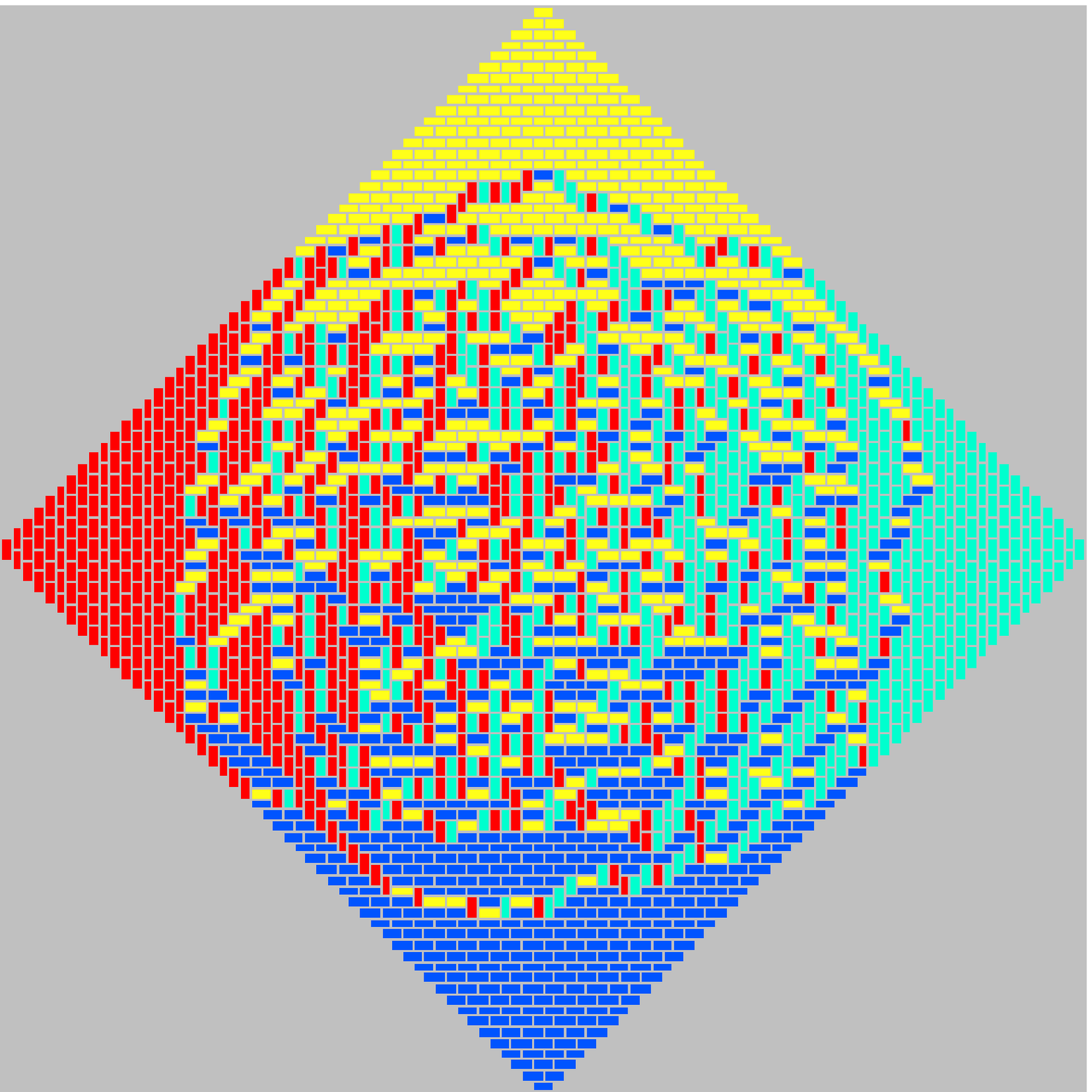}
\label{fig:aztec}
\caption{An Aztec diamond and a random tiling of another
Aztec diamond}
\end{figure}

So as to keep the conversation more concrete, we consider as a running
example the generating function for the probability of a Northgoing
diamond in a uniform random tiling of the Aztec Diamond.  This
example is taken from~\cite[Section~4]{BP-cones}.  We have
$$F(X,Y,Z) = \frac{Z/2}{(1 - (X + X^{-1} + Y + Y^{-1}) Z + Z^2) (1-YZ)}
   = \sum_{r,s,t} a_{r,s,t} X^r Y^s Z^t$$
where $a_{rst}$ is the probability that the domino covering 
the square $(r,s)$ in the order $t$ Aztec Diamond is oriented
in a Northgoing direction.  Figure~\ref{fig:aztec} illustrates
the Aztec Diamond shape of order~4 and the macroscopic features of
a random tiling by dominoes of a larger Aztec Diamond (order~47).  

\begin{figure}[ht]
\centering
\includegraphics[scale=0.5]{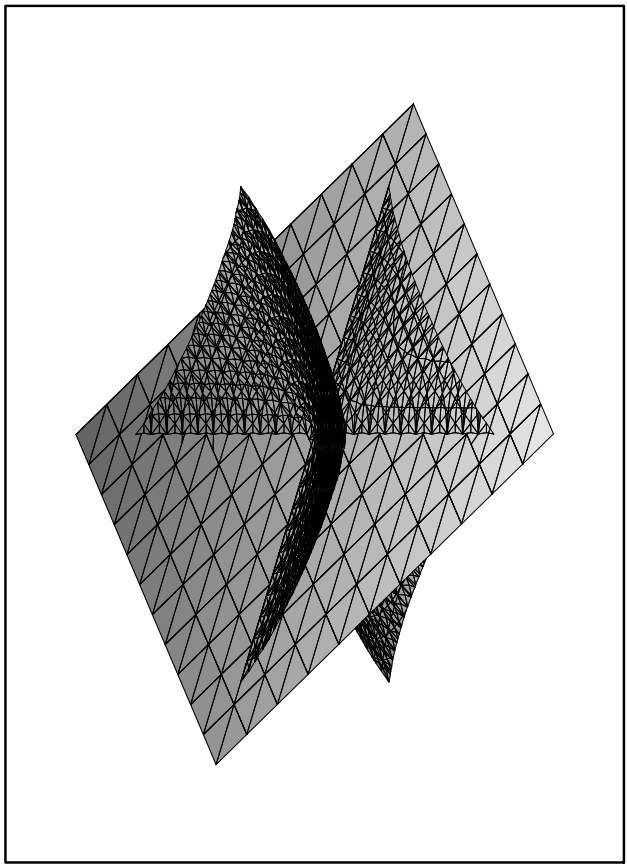}
\caption{The homogenized zero variety near $(1,1,1)$}
\label{fig:cone graph}
\end{figure}

Observe that $r$ and $s$ can be both positive
and negative, while $t$ is always positive and at least $|r| + |s|$.
The denominator $Q (\ZZ)$ is a Laurent polynomial whose zero
set in $(\C^*)^d$ has two isolated singularities at $\pm (1,1,1)$.
Near each of these, the zero sets of the two factors in the
denominator look like a cone and a plane respectively: letting
$Q = J H$, where $J$ is the quadratic and $H = 1-YZ$ is the log-linear
term, we have 
$$\homog (Q,(1,1,1)) = (Z^2 - \frac{1}{2} (U^2 + V^2)) 
   (Y-Z) \, .$$
The zero set of this homogeneous polynomial is shown in 
Figure~\ref{fig:cone graph}.

\begin{figure}[ht]
\centering
\includegraphics[scale=0.5]{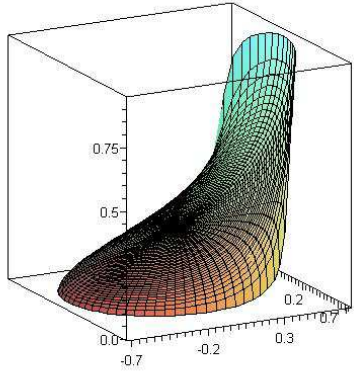}
\caption{Scaling limit of $a_{r,s,t}$}
\label{fig:aztec plot}
\end{figure}

Let $T'$ denote the flat torus $\xx + \R^d / (2 \pi \Z)^d$ 
where $\xx$ is any element of $B$, the component of 
$\amoeba (Q)^c$ over which the series converges.  Changing 
variables via $\ZZ = \exp (\zz) = (e^{z_1}, e^{z_2}, e^{z_3})$ and 
$d\zz = d\ZZ / \ZZ$, then writing $\zz = \xx + i \yy$ 
and $f (\cdot) = F \circ \exp (\xx + i \cdot)$, yields
\begin{equation} \label{eq:log cauchy}
a_\rr = (2 \pi)^{-d} e^{- \rr \cdot \xx} 
   \int_{T'} \exp (- i \rr \cdot \yy) f (\yy) \, d\yy \, .
\end{equation}
Up to this point the expression for $a_\rr$ is exact.  
One can show that approximating $f$ by $\homog (f)$ and 
$T'$ by $\R^d$ does not change the leading asymptotic term.  
This is a somewhat lengthy verification and it is here that
the conical deformation of Corollary~\ref{cor:VF 2} is required.
The upshot is that $a_\rr = E(\rr)$ where $E$ is the inverse
Fourier transform of $\homog (f)$.  For quadratics, the
inverse Fourier transform is the dual quadratic.  For a
product of a quadratic and a log-linear function, an
explicit computation yields~\cite[Theorem~3.9]{BP-cones} 
$$C \arctan \left ( 
   \frac{ \sqrt{ A^* (\rr , \rr)} \sqrt{-A^* (\ell , \ell) }}
   {A^* (\rr , \ell)} \right )$$
where $A^*$ is the quadratic form dual to $\homog (f)$ and $\ell$ is
the the log-linear factor viewed as an element of $(\R^d)^*$ in
the logarithmic space.  For example, when $F$ is the Aztec 
generating function one obtains the following result, a pictorial
version of which is given in Figure~\ref{fig:aztec plot}.
\begin{thm}[\protect{\cite[Theorem~3.7]{BP-cones}}] \label{th:a}
The Northgoing placement probabilities $\{ a_{rst} \}$ for the Aztec 
Diamond are asymptotically given by
$$a_{rst} \sim \frac{1}{\pi} \arctan \left ( 
   \frac{\sqrt{t^2 - 2r^2 - 2s^2}}{t-2s} \right )$$
when $(r,s,t)$ is in the cone $\{ (r,s,t) : t^2 > 2r^2 + 2s^2 \}$
and $r+s+t$ is odd.  
\end{thm}

\begin{unremark}
The way the generating function is indexed, $a_{rst} = 0$ when 
$r+s+t$ is even.  Asymptotics in the corners of the diamond outside
the inscribed circle are given by $a_{rst} \to 1$ in one corner
and $a_{rst} \to 0$ in the other three, where convergence occurs
at an exponential rate as $t \to \infty$. 
\end{unremark}

\clearpage

\begin{center}
{\Huge Part II: Stability}
\end{center}
\vspace{1in}

We now turn to the notion of stability.  As we will see in 
Sections~\ref{ss:LP} and~\ref{sec:stability 2}, the most
useful properties of stability are closure properties:
stability is preserved by a wealth of operations that
are natural from an algebraic or probabilistic point
of view.  Stability is defined for general multivariate
complex polynomials.  Some of these properties or their
proofs simplify when restricted to certain subclasses,
such as polynomials whose coefficients are real or positive,
or \Em{multi-affine} polynomials, whose degrees in each 
variable never exceed~1.  Even in these cases, however, 
some properties are proved only by going through the more
general setting of complex coefficients.

\setcounter{equation}{0}
\section{Stability theory in one variable} \label{sec:stability 1} 

As usual, the univariate theory is older and simpler.  As with
hyperbolicity, the origins of stability theory are in differential
equations and control theory.  After discussing this, we turn
to combinatorial uses of the univariate theory.  Stable generating 
polynomials produce coefficient sequences satisfying Newton's 
inequalities, implying, among other things, log-concavity.
We end the section on univariate stability with a discussion of
stability preserving operations and the so-called Laguerre--P{\'o}lya
class.

\subsection{Stability over general regions} \label{ss:hurwitz}

In this section we will begin by thinking more generally
of polynomials that avoid roots in some region $\Omega$.
We will call these $\Omega$-stable.  Throughout the 
remainder of this paper we will let
$$\halfplane := \{ z : \Im \{ z \} > 0 \} \, .$$
denote the open upper half-plane.  I will always
use the term ``stable'' to mean $\halfplane$-stable, 
as in Definition~\ref{def:stable}, and will use
``Hurwitz stable'' or ``$\Omega$-stable'' for other
notions of stability.
Recalling G.-C.~Rota's philosophical observation, one
might keep in mind that the Riemann Hypothesis is
equivalent to $\Omega$-stability where $\Omega = \{ 
z : \Re \{ z \} \neq 1/2 \}$.  This may seem like a 
strained connection, but in fact some of the literature
we review in Section~\ref{sec:determinants} is explicitly 
motivated by the desire to understand the zeta function.

Statistical physicists have a different motive for understanding
regions free of zeros.  To explain this, we examine some graph 
theoretic models, paraphrasing the exposition in~\cite{sokal-zeros}.  
Let $G = (V,E)$ be a finite graph and let $q$ be a positive integer.  
Define a polynomial in variables $\{ x_e : e \in E \}$ 
$$Z(q , \{ x_e : e \in E \}) := \sum_{\sigma} 
   \prod_{e = \overline{vw} \in E} 
   \left ( 1 + x_e \delta_{\sigma (v) , \sigma (w)} \right )$$
where $\sigma$ ranges over $q$-colorings of the vertices, that is,
all maps from $V$ to $\{1 , \ldots , q \}$.  This is elementarily 
seen to be equivalent to the alternative definition
$$Z(q , \{ x_e : e \in E \}) := \sum_{E' \subseteq E} q^{k(E')}
   \prod_{e \in E'} x_e$$
where $k(E')$ denotes the number of connected components in the
subgraph $(V,E')$.  The last formulation makes it obvious that
$Z$ is a polynomial in $q$ as well as in $\{ x_e \}$.  When
$x_e = -1$ for all $e$ this specializes to the chromatic polynomial
$\chi (q)$; more generally, taking $x_e = x$ for all $e$ gives the
bivariate \Em{Tutte polynomial}.  

The polynomial $Z$ is also the \Em{partition function} for the 
$q$-state Potts model.  Statistical physicists are interested
in \Em{phase transitions} where the behavior of the system 
depends non-analytically on its parameters.  Being a polynomial, 
$Z$ is of course analytic in all its parameters.  Typically though,
one is interested in infinite-volume limits such as the free energy
$$f (q,x) := \lim_{n \to \infty} n^{-1} \log Z_n (q , x)$$
where $Z_n$ is taken on a graph $G_n$ in a family of graphs
whose increasing limit is $\Z^d$.   The function $f$ may
fail to be analytic at $(q,x)$ if the zeros of $Z_n (q,x)$
have a limit point at $(q,x)$.  Therefore, there is physical
significance in keeping the zeros of the polynomial $Z_n$ 
out of specified regions.  Such results are often known as
\Em{Lee-Yang} theorems, after~\cite{lee-yang1,lee-yang2}
(see also~\cite{lieb-sokal}).  

One such theorem was proved by D.\ Wagner in 2000.  The 
\Em{reliability polynomial} $R_G$ of the graph $G = (V,E)$ is
the probability that the graph is connected when each edge
is kept or deleted independently with probability $x$, which
is a polynomial in the parameter $x$.
Brown and Colbourn conjectured~\cite{brown-colbourn}
that all zeros of $R_G (x)$ are in the closed disk $\{ z :
|z-1| \leq 1 \}$.  In other words, $R_G$ is ${\cal D}^c$-stable
where ${\cal D} := \{ z : |z-1| \leq 1 \}$.  This was proved for the class 
of series-parallel graphs in~\cite[Theorem~0.2]{wagner-reliability}.  A 
simpler proof was found by Sokal~\cite[Section~4.1,~Remark~3]{sokal-zeros},
who actually proved the stronger multivariate stability result.
On the other hand, Royle and Sokal~\cite{sokal-royle}
showed that the original univariate conjecture is false for 
general graphs; they also showed that multivariate ${\cal D}^c$-stability
of $R_G$ holds {\em if and only if} $G$ is series-parallel.
Sokal~\cite{sokal-lectures} has recently conjectured that
$R_G$ is (univariate) ${\cal D}^c$-stable for the complete 
graphs, $K_n$.

\subsubsection*{Hurwitz stability}

The association of zero-free regions with the term ``stability'' 
originated in ODE's and control theory.  Its use is attributed 
to Hurwitz.  In fact, $\Omega$-stability with $\Omega = 
\{ z : \Re \{ z \} \geq 0 \}$ is sometimes called 
``Hurwitz stability''.  In ODE's, it is easy to see 
the physical significance of Hurwitz stability.  Let $M$
be a matrix and consider the linear system $\yy' = M \yy$.  
Hurwitz stability is equivalent to all eigenvalues of $M$ having 
negative real parts, which is equivalent to all homogeneous
solutions decaying, and hence good control over the system
$\yy' - M \yy = \vv (t)$.
One may also consider discrete-time analogues, such as 
the system $\yy^{(n+1)} - Q \yy^{(n)} = \vv^{(n)}$
where $Q$ would correspond to the exponential of $M$ 
in the previous system.
Now decay of homogeneous solutions is equivalent to 
$\Omega$-stability when $\Omega$ is the complement 
of the open unit disk; polynomials whose zeros
are all in the open unit disk are said to be \Em{Schur-stable}.

For a less trivial example we turn to control theory.
Following~\cite[Section~10.3]{henrici-2}, we consider
systems which turn an input signal $f(t)$ into an output
$u(t)$.  Many systems, and in particular those built from 
networks of impedances, are not only linear but also time
homogeneous, that is, the map $f \mapsto u$ commutes
with time translation.  In terms of Laplace transforms, this
means that the system acts multiplicatively, meaning that
if $\lap$ denotes the Laplace transform, then such systems 
obey the law
\begin{equation} \label{eq:transfer}
\lap u = g \cdot \lap f
\end{equation}
for some function $g : [0,\infty) \to \R$ that is called the
\Em{transfer function}.  For example an L-C-R circuit, with
an input voltage function $f(t)$ between the inductor and the 
resistor producing an output $u(t)$ of current from the inductor 
to the capacitor, satisfies~\eqref{eq:transfer} with $g(s) = 
(R + L s + C^{-1} / s)^{-1}$.  In impedance networks the transfer
function is in fact always rational, $g(s) = p(s) / q(s)$ for
some polynomials $p$ and $q$.  Generally, a system is said to
be stable if bounded input produces bounded output.  

To check whether the system with transfer function $g$ is stable,
it suffices to test it on inputs of the form $f(t) = e^{i \omega t}$,
the Laplace transform of which is $(\lap f) (s) = 1/(s - i \omega)$.
Let the poles of the rational function $g$ be $a_1 , \ldots , a_r$.
The $\lap u = g \cdot \lap f$ has poles $a_1 , \ldots , a_r , i \omega$.  
We may write $\lap u$ in a partial fraction expansion resulting
in $\sum_{j=0}^r p_j / (s - a_j)$ for some polynomials $\{ p_j \}$;
here we have set $a_0 := i \omega$.  Inverting the Laplace transform
gives a sum $\sum_{j=1}^r P_j (t) e^{a_j t}$ for some collection
$\{ P_j \}$ of polynomials.  This is bounded if and only if 
each $a_j$ is either in the open left half-plane or is 
imaginary and a simple zero.  But when $g$ has an imaginary
zero $a_j$, then setting $a_0 = a_j$ (that is, taking the input 
to be $e^{i a_j t}$) produces a doubled root and an unbounded
output.  Hence, bounded inputs produce bounded outputs if and only
if all poles of $g$ lie strictly in the left half-plane.  

We conclude that stable behavior of the system corresponds to 
Hurwitz stability of the denominator of the transfer function $g$.
For example, in the L-C-R circuit above, the denominator of
$g$ is a quadratic with positive coefficients.  The real parts
of the roots are always negative, whence such a system always 
behaves stably.  

Differential equations of order $d$ may be transformed into
first order systems in $d$ variables by the well known trick 
of representing the first $d-1$ derivatives as new variables.
It is not surprising, therefore, that Hurwitz stability also
arises in the stability analysis of order-$d$ linear differential 
equation with constant coefficients.  Given $d$ initial conditions 
and an inhomogeneous term, we may write such an equation as
\begin{equation} \label{eq:diffeq2}
f^{(d)} + a_{d-1} f^{(d-1)} + \cdots + a_0 f = H
\end{equation}
with initial conditions
$f^{(j)} (0) = b_j$ for $0 \leq j \leq d-1$.
Assuming $H$ to grow at most exponentially, a Laplace transform 
exists near the origin.  We may take the Laplace transfrom of 
both sides of~\eqref{eq:diffeq2}.  Using linearity and the
rule $\lap f' (s) = s \lap f (s) - f(0^+)$ and denoting
$g := \lap f$ gives inductively
$$\lap (f^{(d)}) (s) = s^d g(s) - s^{d-1} f(0^+) - s^{d-2} f' (0^+) - 
   \cdots - f^{(d-1)} (0^+) \, .$$
Plugging this and the boundary conditions into~\eqref{eq:diffeq2}
yields the following equation for $g$:
\begin{eqnarray*}
s^d g(s) - s^{d-1} b_0 - s^{d-2} b_1 - \cdots - b_{d-1} && \\
+ \; a_{d-1} \left [ s^{d-1} g(s) - s^{d-2} b_0 - \cdots - b_{d-2} 
   \right ] && \\[1ex]
\vdots \hspace{0.5in} && \\[1ex]
+ \; a_1 \left [ s^1 g(s) - b_0 \right ] && \\
+ \; a_0 g(s) & = & \lap H (s) \, .
\end{eqnarray*}
Letting $p(s) := s^d + a_{d-1} s^{d-1} + \cdots + a_0$ denote 
the characteristic polynomial of the equation~\eqref{eq:diffeq2},
and $p_k (s)$ denote the shifted polynomial $s^k + a_{d-1} s^{k-1}
+ \cdots + a_{d-k}$, we may rewrite the equation for $g$ as
$$ p \cdot g = \lap H + \sum_{j=0}^{d-1} b_j p_{d-1-j} $$
and hence
\begin{equation} \label{eq:LT}
g(s) = \frac{\lap H (s) + \sum_{j=0}^{d-1} b_j p_{d-1-j}(s)} {p(s)} \, .
\end{equation}
Again, poles of $g$ with positive real part produce unbounded
output, as do purely imaginary poles $i \omega$ when the driving
term is taken to be $e^{i \omega t}$, whose Laplace transform 
$\lap H$ has a pole at $i \omega$.  Thus, again, Hurwitz stability
is equivalent to bounded output on bounded input.

\subsection{Real roots and Newton's inequalities}

I will now return, permanently, to upper half-plane stability, 
which will be called, simply, ``stability''.  In one variable, 
when the coefficients of $f$ are real, the zeros come in 
conjugate pairs.  Stability, therefore, is equivalent to having 
only real zeros.  When the coefficients are nonnegative, a strictly
positive zero is impossible, whence stability is further equivalent 
to having all zeros on the negative half line.  

Combinatorialists have long sought to exploit the properties 
of stable generating functions.
The universal problem in combinatorial enumeration is to count
a family of structures indexed by one or more positive integer parameters.
In the case of one parameter, say $n$, the sequence $\{ a_n \}$ of counts
and the corresponding generating function $f(z) := \sum_n a_n z^n$ are of
fundamental interest.  In the setting of probability generating
functions, the coefficients $a_n$ are nonnegative and sum to one.
In this case, if $f$ has only real roots, then the distribution
which gives probability $a_n$ to the value $n$ is representable
as the sum of independent random variables each taking the value~0
or~1 (Bernoulli random variables).  There may or may not be a 
natural interpretation for these Bernoulli variables; 
see~\cite[Example~23]{HKPV} for an example in which there
can be no natural interpretation.

Often it is intuitively
plausible that the sequence $\{ a_n \}$ is \Em{unimodal}, meaning
that for some $k$, we have $a_0 \leq \cdots \leq a_k \geq a_{k+1}
\geq a_{k+2} \geq \cdots$.  The enumeration literature is littered
with examples in which unimodality is conjectured (see for example
the survey~\cite{stanley-unimodal} and the follow-up to
this~\cite{brenti-unimodal}), but a proof is often elusive 
for the reason that there is no obvious theoretical framework 
within which to prove unimodality.  There are, however some
stronger properties for which natural avenues of proof exist.

\begin{defn}[log-concavity]
A finite or infinite sequence $\{ a_k \}$ of nonnegative numbers
is said to have no internal zeros if the indices of the nonzero
terms form an interval $[r,s]$.  The sequence is said to be 
\Em{log-concave} if it has no internal zeros and if 
$$a_{k-1} a_{k+1} \leq a_k^2$$
for $k \in [r+1,s-1]$.  The sequence $\{ a_r , \ldots , a_{r+k} \}$
with no internal zeros is said to be \Em{ultra log-concave} if 
$$\left ( \frac{a_j}{\binom{k}{j}} \right )^2 \geq 
   \frac{a_{j+1}}{\binom{k}{j+1}} \; 
   \frac{a_{j-1}}{\binom{k}{j-1}} \; .$$
\end{defn}

It is immediate that ultra log-concavity implies log-concavity which
implies unimodality.  While ultra log-concavity appears to be the 
least natural of these properties, it was shown three centuries
ago by Newton to follow from stability of the generating function.

\begin{thm}[Newton's inequalities] \label{th:newton}
Suppose $f(z) := \sum_{k=0}^n a_k z^k$ is a real stable polynomial.
Then for $1 \leq k \leq n-1$, 
\begin{equation} \label{eq:newton}
\left ( \frac{a_k}{\binom{n}{k}} \right )^2 \geq 
   \frac{a_{k+1}}{\binom{n}{k+1}} \; 
   \frac{a_{k-1}}{\binom{n}{k-1}} \; .
\end{equation}
If furthermore the coefficients $a_k$ are all nonnegative then
the roots of $f$ are all nonpositive and the sequence is 
ultra log-concave.
\end{thm}

\noindent{\sc Proof:} By Rolle's Theorem, if a univarite polynomial 
$f$ has only real zeros then so does its derivative.  This observation
will be useful many times below.  We use it now to deduce that 
$$Q (z) := \left ( \frac{d}{dz} \right )^{k-1} f(z)$$
has only real zeros.  Reversing a sequence of coefficients
via $R(z) := z^{n-k+1} Q(1/z)$ also preserves the property of
having all real roots.  By Rolle's Theorem again, $S(z) := 
(d/dz)^{n-k+1} R(z)$ has all real roots.  But $S(z)$ is the
trinomial
$$\frac{n!}{2} \left ( \frac{a_{k-1}}{\binom{n}{k-1}} \; z^2 + 
   2 \frac{a_k}{\binom{n}{k}} \; z + \frac{a_{k+1}}{\binom{n}{k+1}}
   \right ) \, ,$$
and the theorem follows from the discriminant test for quadratics.
$\Cox$

A curious result in a similar vein was proved by 
Gurvits~\cite[Lemma~3.2]{gurvits-vdW}; the proof, which
is a few lines of calculus, is omitted.
\begin{pr} \label{pr:1 bound}
Let $f$ be a probability generating polynomial of degree $d$ and let
$C := \inf f(t)/t$.  If $f$ is stable then
$$a_1 = f'(0) \geq \left ( \frac{d-1}{d} \right )^{d-1} \; C \, .$$
Equality holds if and only if $f = (1 - q + q t)^d$ generates
a binomial distribution.  The value $[(d-1)/d]^{d-1}$ increases
to $e^{-1}$ as $d \to \infty$; for infinite series, the bound 
$a_1 \geq e^{-1} C$ holds with equality if and only if
$f(t) = e^{\lambda (t-1)}$ is the generating function for
a Poisson distribution.
$\Cox$
\end{pr}

There are a number of techniques that can be used to establish
stability (also known, in the univariate case, as the real root 
property).  In some cases
one has a sequence $\{ f_n \}$ where $f_n$ has degree $n$ and
the roots of $f_{n-1}$ interlace the roots of $f_n$, meaning
that each interval between roots of $f_n$ contains a root of $f_{n-1}$.
If this is true, then often it is possible to prove it by induction.
The following example from~\cite{stanley-unimodal} illustrates this.

\begin{example}[Hermite polynomials] \label{eg:hermite}
The Hermite polynomials are a sequence of polynomials defined by 
$$H_n (z) := \sum_{k=0}^{\lfloor n/2 \rfloor} \frac{(-1)^k n! (2z)^{n-2k}}
   {k! (n-2k)!} \, .$$
The are orthogonal with respect to the Gaussian measure $e^{-x^2} \, dx$.
They satisfy the recursion
$$H_n (z) = - e^{z^2} \frac{d}{dz} 
   \left ( e^{-z^2} H_{n-1} (z) \right ) \, .$$
Assume for induction that $H_{n-1} (z)$ has $n-1$ real zeros.  
It is clear from this that $H_n$ has $n-2$ zeros interlacing the
zeros of $H_{n-1}$.  To see that $H_n$ also has a zero less than
all the zeros of $H_{n-1}$, observe that $e^{-z^2} H_{n-1}$ tends
to zero as $z \to -\infty$, therefore $e^{-z^2} H_{n-1}$ has 
an extreme value to the left of its leftmost zero, which is 
a zero of $H_n$.  Similarly, $H_n$ has a zero greater than
the greatest zero of $H_{n-1}$, and the induction is established.
\end{example}

There is a more methodical way that often works to prove interlacing 
by induction.  If $f$ and $g$ are real polynomials of respective 
degrees $n$ and $n-1$ having all real roots and the zeros of $g$ 
interlace the zeros of $f$, it may be seen elementarily that  
$f + \lambda g$ has $n$ real zeros for any real $\lambda$.   
A simple application of this idea is to the sequence
$\{ g_n \}$ defined by the three-term recurrence
$$g_{n+1} (x) = a x g_n (x) + b g_{n-1} (x) \, ,$$
which include the Chebyshev polynomials and Laguerre
polynomials ($a$ and $b$ depending on $n$ in the latter case).
The same idea is the basis for a theorem independently 
proved by Heilmann and Lieb~\cite{heilmann-lieb}, by Gruber 
and Kunz~\cite{gruber-kunz} and, in part, by Nijenhuis~\cite{nijenhuis}.
A matching on a finite weighted graph or multi-graph is a subset
$M$ of the edge set such that any two edges in $M$ are
disjoint (share no vertex).  Fix a graph $G = (V,E)$ and
let $\{ W(e) : e \in E \}$ be a set of nonnegative weights.
Let $a_j$ denote the number of matchings with $j$ edges, counted
by weight, where the weight of a matching is the product
of the weights of the edges in the matching; by convention
we take $a_0 :- 1$.  In the simplest
case, $W(e) \equiv 1$ and $a_j$ simply counts matchings of size $j$,
size being the number of edges.

If we are allowed to set $W(e) = 0$ for some edges, then the
complete graph is universal.  We therefore assume that $G = K_n$
for some $n$.

\begin{thm}[matchings]
The number of weighted matchings of $K_n$ enumerated by size
is ultra log-concave.
\end{thm}

\noindent{\sc Proof:} We show stability of a certain generating
function, the most convenient being 
\begin{equation} \label{eq:Q}
Q_n (z) := \sum_{j=0}^{\lfloor n/2 \rfloor} (-1)^j a_j z^{n-2j} 
   \, .
\end{equation}
This is monic, of degree $n$, and is odd or even depending on $n$.
If we show $Q_n$ is stable and has $n$ distinct roots, then letting
$Q_n(z) = R_n(z^2)$ in the even case and $Q_n(z) = z R_n(z^2)$ in the
odd case, it follows that $R_n(-z)$ has $\lfloor n/2 \rfloor$
negative real roots, hence is stable; it is also the generating
function for matchings by size, hence the proof will be complete.

Taking limits at the end, we may assume without loss of generality 
that the weights $W(e)$ are all strictly positive.  To see that 
$Q_n$ is stable, let $Q_n^S$ denote the generating function 
defined by~\eqref{eq:Q} on the graph $K_n \setminus S$.  We 
use the recursion
\begin{equation} \label{eq:match recur}
Q_n (z) = z Q_n^u (z) 
   - \sum_{v \neq u} W(u,v) Q_n^{u,v} (z)
\end{equation}
which is obvious from the definition.  The induction hypothesis
is that $Q_n^u$ is stable with distinct roots, which are
interlaced by the roots of $Q_n^{u,v}$ for any $v \neq u$.

To verify the induction, we consider the sign of $Q_n$ at
the $n-1$ zeros of $Q_n^u$ in decreasing order.  By the
interlacing property, the sign of each $Q_n^{u,v}$ alternates,
starting out positive because $Q_n^{u,v} (+\infty) = +\infty$.
Therefore, $z Q_n^u - \sum_v W(u,v) Q_n^{u,v}$ alternates,
starting out negative.  This implies the existence of 
$n-2$ roots of $Q_n$ interlaced by the roots of $Q_n^u$.
But also there is one root of $Q_n$ to the right of every root 
of $Q_n^u$ because $Q_n$ is negative at the rightmost root
of $Q_n^u$ and positive at $+\infty$.  Similarly there is
one root of $Q_n$ to the left of every root of $Q_n^u$.  This
completes the inductive proof.
$\Cox$

A second method by which one can establish univariate stability
is via closure properties of the class of univariate stable
polynomials.  This topic will be expanded in the next section,
but for now we mention a result from~\cite{brenti-LC},
whose proof we will omit.
Let $(x)_i$ denote the falling product $x (x-1) \cdots (x-i+1)$
and let $(x)^i$ denote the rising product $x (x+1) \cdots (x+i-1)$.
\begin{thm}[\protect{\cite[Theorems~2.4.2--2.4.3]{brenti-LC}}]
\label{th:falling}
Let $f(x) := \sum_{k=0}^n a_k x^k$ be real stable with nonnegative
coefficients.  Then the polynomials $\sum_{k=0}^n a_k (x)_k$
and $\sum_{k=0}^n a_k (x)^k$ are real stable as well.
$\Cox$
\end{thm}

In the following example, a function from the set $\{ 1 , \ldots , n \}$ 
to itself is represented as a directed graph with an edge from
$j$ to $f(j)$ for each $j$.

\begin{example}[functions enumerated by components]
Let $b(n,k)$ be the number of functions from the set 
$\{ 1 , \ldots , n \}$ to itself whose directed graph has 
precisely $k$ components.  It is elementary, 
cf.~\cite[Example~3.3.28]{goulden-jackson}, that
$$b(n,k) = \sum_{i=1}^n \binom{n-1}{i-1} n^{n-i} c(i,k)$$
where $c(i,k)$ is a signless Stirling number of the first kind.
Multiplying by $x^k$ and summing over $k$ gives
$$f_n (x) := \sum_{k=1}^n b(n,k) x^k 
   = \sum_{i=1}^n  \binom{n-1}{i-1} n^{n-i} (x)^i \, .$$
The corresponding sum replacing $(x)^i$ by $x^i$ has the short
closed form
$$\sum_{i=1}^n \binom{n-1}{i-1} n^{n-i} x^i = x (x+n)^{n-1} \, .$$
Evidently this has only real zeros, so applying the conclusion of
Theorem~\ref{th:falling}, we see that $f_n$ is stable as well.
\end{example}

A third method is an equivalent characterization for ultra log-concave
that is due to Edrei~\cite{edrei}.  This relies on the notion of
\Em{total positivity}, a theme developed at length by 
Karlin~\cite{karlin-TP} and for which a number of combinatorial 
applications are given in~\cite{brenti-TP}.  

\begin{defn}[Total positivity; P\'olya frequency sequence] 
\label{def:TP} ~~

$(i)$ An infinite matrix is said to be totally positive (TP)
if all of its minors are nonnegative. 

$(ii)$ A sequence $\{ a_n \}$ is said to be 
a P\'olya frequency sequence (PF-sequence) if the matrix
$A_{nk} := a_{n-k}$ is totally positive.  Here we
take $a_i := 0$ if $i < 0$ or the sequence $\{ a_n \}$
is finite and has length less than $i$.
\end{defn}

A proof of the following equivalence may be found in
in~\cite[Theorem~5.3]{karlin-TP}.
\begin{thm}[Edrei's equivalence theorem] \label{th:edrei}
The nonnegative sequence $(a_0 , \ldots , a_d)$ is a PF-sequence
if and only if the generating polynomial $\sum_{k=0}^d a_k z^k$
is stable.
$\Cox$
\end{thm}

Brenti~\cite[Section~3]{brenti-unimodal} points out that only
the $2 \times 2$ minors are required for log-concavity, hence
unimodality.  Nevertheless, nonnegativity of all minors has
a combinatorial interpretation which is given 
in~\cite[Theorem~3.5]{brenti-TP}.  This interpretation is
somewhat abstract, but in special cases the interpretation
can be more concrete.  

\begin{example}[$r$-derangements] \label{eg:derange}
A permutation of $\{ 1 , \ldots , n\}$ is said to be an 
$r$-derangement if all its cycles have length at least~$r$.  
Thus a 1-derangement is any permutation, a 2-derangement is
a classical derangement, and so forth.  Let $b(n,r)$
count the number of $r$-derangements of $\{ 1 , \ldots , n\}$.  
Brenti~\cite[Corollary~5.9]{brenti-TP} shows that
$f_n := \sum_{r=1}^n b(n,r) z^r$ is stable.  To do
so, he shows that a related sequence $\{ c_r (n) \}$
is a PF-sequence, hence stable, then applies 
Theorem~\ref{th:falling} to recover the result 
with $b(n,r)$ in place of $c_r (n)$.
\end{example}

We end this section with one more example of a stable
generating polynomial, this one taken from~\cite{stanley-unimodal}.
Recall from Example~\ref{eg:hermitian} that if $A$ is
a real symmetric matrix and $B$ is real positive semi-definite
then $f(z) := \Det (A + z B)$ is a real stable polynomial.

\begin{example}[spanning forests enumerated by component]
\label{eg:forests}
Let $G$ be a finite graph with possibly multiple edges and let 
$$a_k (G) := \sum_F \gamma (F)$$ 
where the sum is over all spanning forests $F$ of $G$ such that
$F$ has precisely $k$ components, where $\gamma (F)$ is the product 
of the cardinalities
(numer of vertices) of the $k$ components.  The factor of $\gamma$
changes the eumeration from forests to rooted forests.
Following~\cite[Proposition~4]{stanley-unimodal}, let us see that 
the polynomial $f_G := \sum_k a_k (G) z^k$ is stable.  Let $A$ be the
matrix whose rows and columns are indexed by the vertices of $G$
with entry $A_{ij} := \deg (u)$ if $u = v$ and otherwise equal to
$-N(i,j)$ where $N(i,j)$ is the number of edges between $i$ and $j$.
It is known that $f_G = \Det (A + z I)$.  It follows that $f_G$
is stable.  Therefore, the number of rooted forests enumerated by 
components is ultra log-concave, hence unimodal.  Stanley attributes
the identity $f_G = \Det (A + z I)$ to Kelmans.
\end{example}

\subsection{The Laguerre--P\'olya class}
\label{ss:LP}

Suppose we transform a polynomial $f(z) = \sum_{k=1}^n a_k z^k$
by multiplying each coefficient $a_k$ by a specified constant 
$\lambda_k$.  We may ask which sequences $\blambda := \{ \lambda_k \}$ 
are \Em{multiplier sequences}, meaning that the resulting operator 
$T_\blambda$ preserves the class of real stable polynomials.  
A classical theorem due to P\'olya and Schur~\cite{polya-schur14} 
gives a complete and beautiful characterization of multiplier sequences.

\begin{thm}[P\'olya-Schur 1914] \label{th:PS}
Let $\blambda = \{ \lambda_n : n \geq 0 \}$ be a sequence
of real numbers and let $T_\blambda$ denote the linear operator
defined by 
$$T_\blambda \left ( \sum_{k=0}^n a_k z^k \right ) = 
   \sum_{k=0}^n \lambda_k a_k z^k \, .$$
Denote by $\Phi$ the formal power series
$$\Phi (z) := \sum_{k=0}^\infty \frac{\lambda_k}{k!} z^k \, .$$
Then the following are equivalent.
\begin{enumerate}[(i)]
\item $\blambda$ is a multiplier sequence;
\item $\Phi$ is an entire function and is the limit, uniformly on
compact sets, of the polynomials with all zeros real and of the 
same sign;
\item $\Phi$ is entire and either $\Phi (z)$ or $\Phi (-z)$ 
has a representation 
$$C z^n e^{\alpha_0 z} \prod_{k=1}^\infty (1 + \alpha_k z)$$
where $n$ is a nonnegative integer, $C$ is real, and $\alpha_k$ 
are real, nonnegative and summable.
\item For all nonnegative integers $n$, the polynomial
$T_\blambda [ (1+z)^n ]$ is real stable with all roots 
of the same sign.
\end{enumerate}
\end{thm}

We will not prove the general P\'olya-Schur theorem here, 
proving only some special cases later as we need them.
Some examples and remarks will clarify its meaning and
possible uses.  Condition~$(iv)$ may be thought of as
saying that the polynomials $(1+z)^n$ are universal
test cases for stability preserving: a multiplier sequences
preserving stability of these will preserve stability of
all real stable polynomials.

\begin{example}[dilation] \label{eg:PS1}
If $b$ is a nonnegative integer, setting $\lambda_k := b^k$
produces the operator $T_\blambda f (z) = f(b z)$.  Clearly
this preserves stability.  The representation in~$(iii)$ is 
obvious because $\Phi (z) = e^{bz}$.  
\end{example}

\begin{example}[factorials] \label{eg:PS2}
If $n \geq 1$ is an integer then the sequence $\lambda_k 
:= (n)_k := n! / (n-k)!$, defined to vanish when $k > n$,
produces the exponential generating function $\Phi (z) =
(1+z)^n$.  By criterion~$(iii)$ this is a multiplier sequence.
Dividing by the constant $n!$ we see that $\lambda_k = 1/(n-k)!$
for $k \leq n$ and zero for $k > n$ defines a multiplier sequence.
For real polynomials of degree $n$, stability of $f$ is equivalent 
to stability of the inversion $z^n f(1/z)$, hence $\lambda_k = 1/k!$
defines a multiplier sequence on polynomials of degree $n$.  This
is true for every $n$, whence $\{ 1/k! \}_{k=0}^{\infty}$ 
is a multiplier sequence.  This result is due to Laguerre.
\end{example}

\begin{example}[coefficientwise multiplication] \label{eg:PS3}
Applying the previous example to a polynomial $g(z) = 
\sum_{k=0}^n a_k z^k$ with negative real roots shows that 
$\sum_{k=0}^n (a_k / k!) z^k$ also has negative real roots.
Setting $\lambda_k := a_k$ for $k \leq n$ and zero for $k > n$
we see that its exponential generating function $\Phi (z) = 
\sum_{k=0}^n a_k z_k / k!$ is of the form in~$(iii)$ with
$\alpha_0 = 0$.  We conclude that $\{ a_k \}$ is a multiplier
sequence.  In other words, if we multiply term by term 
the coefficient sequences of two real rooted polynomials,
at least one of which has nonnegative coefficients, we get
another real rooted coefficient sequence.  Thus the real root
property for polynomials with nonnegative coefficients is
closed under \Em{Hadamard products}.
\end{example}

Multiplier sequences are related to the so-called 
\Em{Laguerre--P\'olya} class.  An entire function 
is said to be in the Laguerre--P{\'o}lya class, 
denoted $\LP$, if it is the limit, uniformly on 
compact subsets of $\C$, of polynomials with only 
real zeros; the terminology goes back at least 
to~\cite{schoenberg}.  It is said to belong to the 
subclass $\LPI$ if it is the limit of polynomials with only
real zeros all of the same sign.  Thus Theorem~\ref{th:PS}
asserts that generating functions of multiplier sequences
are precisely the functions of class $\LPI$.

A related notion to that of a multiplier sequence is the 
notion of a \Em{complex zero decreasing sequence} (CZDS).
Let $Z_c (f)$ denote the number of non-real zeros of $f$.
Say that a finite or infinite sequence $\{ a_k \}$ is
a CZDS if for any real polynomial $f(z) = \sum_{k=0}^n b_k z^k$,
\begin{equation} \label{eq:CZDS}
Z_c \left ( \sum_{k=0}^n a_k b_k z^k \right ) \leq Z_c (f) \, .
\end{equation}
In particular, in order for $\{ a_k \}$ to be a CZDS, 
a value of zero on the right of~\eqref{eq:CZDS} (no non-real zeros)
implies no non-real zeros on the left,
so any CZDS is a multiplier sequence.  The converse, however, 
is not true.  

To see that there is any nontrivial CZDS, we observe that 
the operator $z (d/dz)$ is represented by the sequence
$a_k := k$ and can never increase the number of non-real 
zeros.  Although not every multiplier sequences is a CZDS,
each multiplier sequence leads to a CZDS via the following 
result going back to Laguerre.

\begin{pr}
Let $\Phi \in \LP$ have zeros only in $(-\infty , 0]$.
Then the sequence $\{ \Phi (k) : k \geq 0 \}$ is a CZDS.
$\Cox$
\end{pr}

Related to the notions of multiplier sequences and CZDS 
is the notion of a multiplier sequence for the property
of being nonnegative on real inputs.  Say that $\{ c_k \}$
is a $\Lambda$-sequence if multiplication of coefficients
term by term preserves the property of being everywhere nonnegative:
$$\sum_{k=0}^n a_k x^k > 0 \mbox{ for all } x \in \R \Longrightarrow 
   \sum_{k=0}^n c_k a_k x^k > 0 \mbox{ for all } x \in \R \, .$$ 
Note that negative values are allowed for the numbers $c_k$.
The following classical relationship is proved in~\cite{widder}.
\begin{pr}
If $\{ a_k \}$ is a CZDS then $\{ a_k^{-1} \}$ is a $\Lambda$-sequence.
$\Cox$
\end{pr}

Connections abound between these notions: stability, multiplier
sequences, the Laguerre--P\'olya class, $\Lambda$-sequences, etc.,
which we will not have time to survey here; the reader is referred
to~\cite{craven-csordas} for an introduction.  $\Lambda$-sequences,
for example, are neatly characterized by a determinant
condition and also by the so-called Hamburger moment problem:
$\{ c_n \}$ is a $\Lambda$-sequence if and only if $c_n =
\int_{-\infty}^\infty t^n \, d\mu (t)$ for some measure 
$\mu$ on $\R$ not supported on finitely many points.
At the root of much of the present interest about these 
properties is their connection to the Riemann hypothesis.
An account is given in~\cite{csordas-varga,csordas03}.
Finally, a small related literature on nonlinear transformations
that preserve real stable polynomials may be of interest.
It was conjectured independently by S.~Fisk and R.~Stanley
that if the real polynomial $\sum_k a_k z^k$ has all real roots
then so does the polynomial whose $z^k$ coefficient is
$a_k^2 - a_{k-1} a_{k+1}$.  Some progress was made by McNamara
and Sagan in~\cite{sagan-log-concave}; a proof of the
conjecture and more was recently found by 
Br{\"a}nd{\'e}n~\cite{branden-iterated}.

\setcounter{equation}{0}
\section{Multivariate stability} \label{sec:stability 2} 

Let $\halfplane$ denote the open upper half-plane; thus
stability is equivalent to having no zeros in the
region $\halfplane^d$.  We remark on some connections
to other notions of stability.  For homogeneous polynomials,
the zero set is invariant under multiplication 
by $e^{i \theta}$ in each coordinate, therefore 
all notions of half-plane stability coincide.  Any
\Em{circular} region (the interior or exterior of a circle)
is transformed into $\halfplane$ by a M{\"o}bius transformation,
which allows the theory of $\Omega$-stable functions to be
mapped to the ordinary theory of stable functions via a 
bi-rational change of variables whenever $\Omega$ is
the product of circular regions.  This mapping acts nicely
with respect to some aspects such as closure properties
but not as nicely with respect to coefficient sequences.
Thus, for example, it may shown that \Em{Schur-stability},
where $\Omega = {\cal D}^d$ and ${\cal D}$ is the open unit disk, 
has an unexpected closure property: if $f$ and $g$ are
multi-affine Schur-stable polynomials then the \Em{Hadamard
product} $f \hadamard g$ is Schur-stable as well; here the
Hadamard product of $\sum_S a(S) \zz^S$ and $\sum_S b(S) \zz^S$
is defined to be $\sum_S a(S) b(S) \zz^S$.  Our principal
interest is in $\halfplane$-stable polynomials and their
coefficients.  We will not henceforth consider notions of 
stability other than upper half-plane stability.

A number of the closure properties of the class of stable
polynomials extend immediately from the univariate to the 
multivariate setting or are otherwise elementary.  In
particular, if $f \in \C[z_1 , \ldots , z_d]$ is stable 
then so are the polynomials obtained by the following
operations:
\begin{enumerate}[(a)]
\item dilation: replacing $f$ by $f(b_1 z_1 , \ldots , b_d z_d)$, 
where $\{ b_j \}$ are nonnegative constants;
\item permuting the names of the variables;
\item specialization: setting $z_j$ equal to a constant in $\halfplane$;
\item diagonalization: setting $z_j$ equal to $z_i$;
\item inversion: replacing $f$ by $f(-1/z_1 , z_2 , \ldots , z_d)$;
\item differentiation: replacing $f$ by $\partial f / \partial z_j$;
\item limits: $f_n \to f$ uniformly on compact sets and $f_n$ stable 
implies $f$ is stable;
\end{enumerate}
see for example~\cite[Lemma~2.4]{wagner-BAMS}; here, and after,
in results such as the this one, we will take ``stable'' to
include the zero polynomial in order not to have to make 
exceptions.

We also recall that Proposition~\ref{pr:orthant} equates
stability of a real homogeneous polynomial to hyperbolicity
in all directions in the positive orthant.  Although the
definition of hyperbolicity is more complicated for
non-homogeneous polynomials, the corresponding fact characterization
of real stable polynomials is not.
\begin{pr} \label{pr:inhom orthant}
A (non identically zero) real $d$-variable polynomial $f$ is 
stable if and only if for all $\vv \in \R^d$ and $\uu \in \R^d_+$, 
the polynomial $t \mapsto f(\vv + t \uu)$ has only real zeros.
\end{pr}

\noindent{\sc Proof:} If $f$ has a zero $(z_1 , \ldots , z_d)$
in $\halfplane^d$, then writing $z_j = v_j + i u_j$, we see
that $\uu \in \R^d_+$ and $i$ is a zero of $f(\vv + t \uu)$.
Conversely, if $\uu \in \R^d_+$ and $\vv \in \R^d$ such that
$f(\vv + t \uu)$ has a root with nonvanishing imaginary part, 
we may conjugate if necessary to find a root $t = a+bi$ with $b > 0$;  
then $\vv + t \uu$ is a zero of $f$ in $\halfplane^d$.
$\Cox$

There is one geometric feature of functions of two or more
complex variables, absent from the univariate theory, that 
is worth noting because it strengthens closure property~(c).
Let $f$ be any function in $d$ variables, vanishing at 
$\bb = (b_1 , \ldots , b_d)$ and analytic in a neighborhood of $\bb$.
Suppose that $f (\cdot , b_2 , \ldots , b_d)$ is not the
zero function.  Then there is a continuous multi-valued function 
$\phi$ on a neighborhood $U$ of $(b_2 , \ldots , b_d)$ in $\C^{d-1}$
such that $\phi (\zz) \to b_1$ as $\zz \to (b_2 , \ldots , b_d)$
and such that $f (\phi (\zz) , \zz) = 0$.  This perturbation
property follows from the fact that at least one multivariate 
Puiseux series solution exists for $z_1$ as a function of
$z_2 , \ldots , z_d$ which converges to $b_1$.  This leads
to the following lemma.

\begin{lem} \label{lem:closure}
Let $B_j, 1 \leq j \leq d$ be open subsets of $\C$.  Suppose
a polynomial $f$ is nonvanishing on $\prod_{j=1}^d B_j$.
Then $f$ is nonvanishing on $\prod_{j=1}^d \overline{B_j}$
except possibly in two cases: $f$ can vanish on the product
of boundaries $\prod_{j=1}^d \partial B_j$, and $f$ can vanish at
a point $(b_1 , \ldots , b_d)$ with coordinates $\{ b_j : j \in S \}$ 
in the interior of $B_j$ for some set $S \subseteq \{ 1 , \ldots , d \}$ 
provided that $f$ vanishes identically for all values of those 
coordinates, that is, $f(b_1' , \ldots , b_d') = 0$
whenever $b_j' = b_j$ for all $j \notin S$.

In particular, taking $B_j = \halfplane$ for all $j$, 
we see that stability implies no zeros in $\overline{\halfplane}^d$
except for real zeros or degenerate cases.

A further consequence of this lemma is that~(c) can be strengthened 
to include setting $z_j$ equal to a real constant. 
\end{lem}

\noindent{\sc Proof:} Let $f$ vanish at $\bb = (b_1 , \ldots , b_d)$
where $b_j \in \partial B_j$ for $j \leq r$ and $b_j \in B_j$
for $j > r$.  We need to show that the hypotheses are contradicted
if $k > r$ is a coordinate such that $f$ does not vanish identically 
when the remaining coordinates are fixed.  This follows if we 
perturb each coordinates $b_j, j \neq k$ so as to lie in the
open region $B_j$, by sufficiently small amounts so that 
the perturbation solution for $z_k$ lies in $B_k$.
$\Cox$

The above facts are generally quite elementary.  Before we get to 
the fun stuff, there is some serious overhead in deriving further
closure properties.  One route to developing properties of multivariate 
stable functions is via the notion of \Em{proper position}.  This is 
in some sense a multivariate version of interlacing of roots.  The
development here specifically avoids this because I do not find it 
intuitive and because we can do everything we need without it.  
For a development incorporating the notion of proper position, 
see~\cite{wagner-BAMS}.

\subsection{Equivalences}

Recall from Proposition~\ref{pr:orthant} the relation between
stability and hyperbolicity: a real homogeneous polynomial
is stable if and only if it has a cone of hyperbolicity 
containing the (strictly) positive orthant.  A number of
equivalent formulations of stability will prove useful.
If $f \in \C [z_1 , \ldots , z_d]$ is any polynomial, let 
$m$ denote the maximum total degree of any monomial in $f$ 
and define the homogenization $f_H$ of $f$ to be the unique 
degree-$m$ homogeneous polynomial in the variables 
$z_1 , \ldots , z_{d+1}$ such that $f_H (z_1 , \ldots , z_d , 1)
= f(z_1 , \ldots , z_d)$.  We may write explicitly 
$f_H (z_1 , \ldots , z_{d+1}) = z_{d+1}^m \, f(z_1/z_{d+1} , 
\ldots , z_d / z_{d+1})$.  The following proposition is proved
in~\cite[Section~4]{borcea-branden-liggett}.

\begin{pr}[homogenization] \label{pr:homogenization}
The real polynomial $f \in \R[z_1 , \ldots , z_d]$ 
is stable if and only if the homogeneous
polynomial $f_H$ is hyperbolic in all directions $\vv$ such that 
$v_j > 0$ for $1 \leq j \leq d$ and $v_{d+1} = 0$.

If $f$ also has nonnegative coefficients then the four properties
are all equivalent:
\begin{enumerate}[(i)]
\item $f$ is stable;
\item $f_H$ is stable;
\item $f_H$ is hyperbolic with respect to some vector in
the nonnegative orthant;
\item $f_H$ is hyperbolic with respect to every vector in the
positive orthant.
\end{enumerate}
\end{pr}

\noindent{\sc Proof:} For one direction, assume $f_H$ to be 
hyperbolic in every direction $\yy$ with $y_j$ positive when 
$j \leq d$ and $y_{d+1} = 0$.  For any $\xx \in \R^d$ and 
$\yy \in \R_+^d$, the polynomial 
$$t \mapsto f(\xx + t \yy) = f_H ((\xx , 1) + t (\yy , 0))$$ 
is not identically zero because $\lim_{t \to \infty} t^{-m} 
f(\xx + t \yy)$ is equal to $f_H (\yy , 0) \neq 0$.  It
has only real zeros, by hyperbolicity of $f_H$ in direction 
$(\yy , 0)$, hence $i$ is not a root, hence $f(\xx + i \yy) 
\neq 0$; but $\xx \in \R^d$ and $\yy \in \R_+^d$ are arbitrary, 
so this is stability of $f$.

For the other direction, suppose $f_H$ fails to be hyperbolic 
in a direction $\vv$ for which $v_j > 0, 1 \leq j \leq d$, and 
$v_{d+1} = 0$.  Then $f_H (\xx + t \vv)$ has a non-real zero
$t = u$ for some $\xx \in \R^{d+1}$.  The complex conjugate
of $u$ is also a zero so we may assume without loss of generality
that the imaginary part of $u$ is positive.  The vector 
$\yy := \xx + u \vv$ lies in $\halfplane^d \times \R$.  There are
three cases.  First, if $y_{d+1} > 0$ then dividing though
by this gives a zero of $f$ in $\halfplane^d$, implying
$f$ is not stable.  Second, if $y_{d+1} < 0$, dividing yields
a zero in $- \halfplane^d$ and conjugating, we see again that
$f$ is not stable.  Finally, if $y_{d+1} = 0$, we may use
perturbation to reduce to the case where $y_{d+1} \neq 0$.

For the equivalence of the four properties in the case of 
nonnegative coefficients, recall that the cone of hyperbolicity 
of $f_H$ containing some vector $\vv$ is the connected component 
in $\R^{d+1}$ that contains $\vv$ of the set where $f_H$ is nonzero.
By nonnegativity of the coefficients, if this cone contains
any vector in the closed nonnegative orthant then it contains
all vectors in the positive orthant.  This establishes the equivalence
of~(iii) and~(iv).  Equivalence between (iii) and~(ii) is 
Proposition~\ref{pr:orthant}.
Finally, we note that~(ii) implies~(i) by setting $z_{d+1} = 1$
(as is allowed by the last consequence of Lemma~\ref{lem:closure}),
while~(i) implies~(iii) by the first part of this proposition.
$\Cox$

\begin{unremark}
The perturbation argument also proves that stability of $f$
implies stability of $\homog (f)$, where $\homog (f) 
(z_1 , \ldots , z_d) := f_H (z_1 , \ldots , z_d , 0)$
is the leading homogeneous part of $f$; see, 
e.g.,~\cite[Proposition~2.2]{COSW}.
\end{unremark}

\subsubsection*{The special case of multi-affine polynomials}

Suppose $f$ is multi-affine, that is, no variable appears in any
monomial with power two or higher.  The following equivalence is
taken by some to be the definition of stability in the multi-affine
case.  Its significance will become evident in the next section.

\begin{thm}[\protect{\cite[Theorem~5.6]{branden07}}]  
\label{th:mixed partials}
Let $f$ be real and multi-affine.  Then $f$ is stable if and only
if the inequality 
\begin{equation} \label{eq:rayleigh}
\frac{\partial f}{\partial x_i} (\xx) 
   \frac{\partial f}{\partial x_j} (\xx)
   \geq f(\xx) \frac{\partial^2 f}{\partial x_i \partial x_j} (\xx) 
\end{equation}
holds for all real $\xx$ and distinct $i,j \leq n$.
\end{thm}

The proof given here is considerably simpler than the published
proof, which is a ``proper position'' argument.  It comes from 
the same source (P. Br{\"a}nd{\'en}, personal communication) and
begins with the following lemma.

\begin{lem} \label{lem:branden}
Let $Q , R \in \C [z_1 , \ldots , z_n]$ and define $P = Q + z_{n+1} R
\in \C [z_1 , \ldots , z_{n+1}]$.  Let $\Omega \subseteq \C^n$ be any
connected set and let $D_1 , D_2$ be closed subsets of $\C$ with
disjoint interiors, such that $D_1 \cup D_2 = \C$ and $J := D_1 \cap D_2$
a simple curve separating $D_1$ and $D_2$.  

If $P$ has no roots in $\Omega \times J$ and $R$ has no roots in $\Omega$
then $P$ either has no roots in $\Omega \times D_1$ or has no roots in
$\Omega \times D_2$.
\end{lem}

\noindent{\sc Proof:} $P$ having no roots in $\Omega \times J$ is 
equivalent to $- Q(\zz) / R(\zz) \notin J$ for all $\zz \in \Omega$.
Hence $-Q/R$ maps $\Omega$ either to the interior of $D_1$ or
the interior of $D_2$, hence to $D_2^c$ or $D_1^c$, and the
result follows.
$\Cox$

\noindent{\sc Proof of theorem:} For the forward implication, suppose
$f$ is stable.  Let $e_i$ denote the $i^{th}$ standard basis
vector and define $G(s,t) := f(\xx + t e_i + s e_j)$.  For any
$\xx \in \R^d$, the function $G$ is stable in the variables $(s,t)$.
To see this, note that this is equivalent to $f$ having no
zeros in $\overline{\halfplane}^d$ with coordinates $i$ and $j$
in $\halfplane$; this is ruled out by Lemma~\ref{lem:closure}.
By multi-affinity, we may express $G$ as $a + bs + ct + dst$ where 
$a,b,c,d$ are given by partial derivatives:
$$G(s,t) = f(\xx) + \frac{\partial f}{\partial z_i} (\xx) s 
   + \frac{\partial f}{\partial z_j} (\xx) t 
   + \frac{\partial^2 f}{\partial z_i \partial z_j} (\xx) st \, .$$
For bivariate multi-affine polynomials $a+bs+ct+dst$, stability is
equivalent to $ad \leq bc$, and the forward direction follows.

For the reverse direction, use induction on $n$.  When $n=1$ the
result is true because all nonzero real polynomials $a+bx$ are stable.
Assume the reverse result for $n$ and suppose $P = Q + z_{n+1} R$ 
satisfies~\eqref{eq:rayleigh}.  If $Q$ or $R$ is identically
zero then we are done by induction, so assume not.  For any 
real $x$, it is clear that $P(z_1 , \ldots , z_n , x)$ 
satisfies~\eqref{eq:rayleigh}, hence, by induction, the polynomials
$Q, R$ and $Q + xR$ are all stable or identically zero.  If 
$Q + \alpha R \equiv 0$ for some real $\alpha$ then 
$P = (z_{n+1} - \alpha) R$ which is stable.  We may assume, therefore, 
that there is no $\alpha$ for which $Q + \alpha R$ is identically zero.
It follows that $P$ has no zeros in $H^n \times \R$ where $H$ is
the open upper half-plane.  By Lemma~\ref{lem:branden} this
implies that either $P$ or $P(z_1 , \ldots , z_n , -z_{n+1})$ 
is stable.  In the latter case, by the forward direction, 
equation~\eqref{eq:rayleigh} holds with the signs reversed 
for $j = n+1$ and $i \leq n$, and hence with equality for such
$i$ and $j$.  This works out to $(\partial R / \partial x_i) Q 
= (\partial Q / \partial x_i) R$ for all $i \leq n$,
that is, $Q$ and $R$ are multiples of each other, which again
implies that $P$ is stable.
$\Cox$

\subsection{Operations preserving stability}

In a brilliant series of papers~\cite{borcea-branden-LYPS1,borcea-branden-LYPS2,borcea-branden-classification},
Borcea and Br\"and\'en finished off the problem of extending
the P\'olya-Schur Theorem to multivariate stability preserving
maps.  This work not only generalizes from one variable to
several, but also from multiplier sequences to arbitrary
$\C$-linear maps.  I will present here two of their results,
one for multiplier sequences and one giving universal test 
functions for general linear maps, proving only those parts 
that will be of use later.  

As before, if $\blambda = \{ \lambda_\rr : \rr \in \Z_+^d \}$ is 
an array of real numbers, we denote by $T_\blambda$ the 
operator for which 
$$T_\blambda \left ( \sum_\rr a_\rr \zz^\rr \right )
   = \sum_\rr \lambda_\rr a_\rr \zz^\rr \, .$$
Also generalizing from before, we say that $\blambda$ is a 
\Em{multivariate multiplier sequence} if $T_\blambda$ preserves
the class of real stable polynomials.  The following theorem
is proved in~\cite[Theorem~1.8]{borcea-branden-classification}.

\begin{thm}[multivariate multiplier sequences] \label{th:mv multiplier}
The array $\blambda$ is a $d$-variate multiplier sequence
if and only if there are $d$ univariate multiplier sequences
$\blambda^{(1)} , \ldots , \blambda^{(d)}$ such that 
$$\lambda_\rr = \lambda^{(1)}_{r_1} \cdots \lambda^{(d)}_{r_d}$$
and satisfying a further sign condition: either every $\lambda_\rr$
is nonnegative, or every $\lambda_\rr$ is nonpositive, or the 
same holds for $(-1)^{|\rr|} \blambda$.
\end{thm}

\noindent{\sc Proof in one direction:} The easy direction, and
all we will need below, is to see that the product of univariate 
nonnegative multiplier sequences is a multiplier sequence (also
that nonnegative multiplier sequences preserve complex stability,
not just real stability).  
Let $\lambda_\rr = \prod_{j=1}^d \lambda^{(j)}_{r_j}$ where
each sequence $\blambda^{(j)}$ is a univariate multiplier 
sequence.  First assume that $\lambda^{(j)}$ is identically 1
for $j \geq 2$.  Let $f$ be a real stable polynomial of $d$ variables.
Fixing $z_2 , \ldots , z_d$ in the upper half-plane, the
polynomial $f_1 := f(\cdot, z_2 , \ldots , z_d)$ is univariate
stable, and by the stability preserving assumption on 
$\blambda^{(i)}$, we see that $T_{\blambda_1} (f_1)$ is
univariate stable.  This polynomial having no zeros in the
upper half-plane is equivalent to $T_\blambda f$ having no
zeros with $z_1$ in the upper half-plane and the specified
values of $z_2 , \ldots , z_d$.  Because these values were
arbitrary values in the upper half-plane, this finishes the
proof in this special case.  For the general case, write 
$T_\blambda$ as the composition of operators of this
form in each coordinate.
$\Cox$

One consequence of this is that the multi-affine part of
a stable polynomial is stable.
\begin{cor} \label{cor:multi-affine}
Let $f$ be a stable polynomial in $d$ variables and let
$f_{\rm ma}$ denote the multi-affine part of $f$, that is,
the sum of the square-free monomials in $f$.  Then
$f_{\rm ma}$ is stable.
\end{cor}

\noindent{\sc Proof:} The sequence $a_0 = 1, a_1 = 1 , 
a_n = 0$ for $n \geq 2$ is a multiplier sequence.  
By Theorem~\ref{th:mv multiplier}, the multivariate
sequence defined by $\blambda_\rr = 1$ if $\rr_j = 0$
or~1 for all $j \leq d$ and $\blambda_\rr = 0$ otherwise
is a multiplier sequence.  This multiplier sequence
extracts the multi-affine part.
$\Cox$

A complete characterization of linear operators preserving 
stability of complex polynomials was proved 
in~\cite[Theorem~2.3]{borcea-branden-LYPS1}.
We will discuss here only the version of this result 
concerning operators preserving stability of multi-affine 
functions.  We need to define what Borcea and Br\"and\'en 
call the \Em{algebraic symbol} $G_T$ of a linear operator $T$.
Their definition requires a degree bound on each 
variable.  We consider here only the multi-affine case; 
if you compare to the original text, you should set the degree 
bound vector $\kappa$ to $(1, \ldots , 1)$.  Let 
$\C_{\rm ma}[z_1 , \ldots , z_d]$ denote the space of
complex multiaffine polynomials in $d$ variables.  For
a linear operator $T : \C_{\rm ma} [z_1 , \ldots , z_d] 
\to \C [z_1 , \ldots , z_d]$ define its algebraic symbol
to be the element $G_T \in \C [z_1 , \ldots , z_d , w_1 , \ldots , w_d]$
satisfying 
\begin{equation} \label{eq:GT}
G_T (z_1 , \ldots , z_d , w_1 , \ldots , w_d) = 
   \sum_{S \subseteq \{ 1 , \ldots , d \}} \prod_{j \notin S} w_j 
   T \left ( \prod_{j \in S} z_j \right ) \, .
\end{equation}
In other words, $G$ applies $T$ to $\prod_{j=1}^d (z_j + w_j)$
by treating $\{ w_j \}$ as constants and applying $T$ to the 
resulting monomials in $\{ z_j \}$.
\begin{thm}[\protect{\cite[Lemma~3.2]{borcea-branden-LYPS1}}]
\label{th:multi-affine PS}
Let $T : \C_{\rm ma} [z_1 , \ldots , z_d] \to C [z_1 , \ldots , z_d]$ 
be any $\C$-linear operator.  Then $T$ preserves stability if and
only one of the two following conditions hold:
\begin{enumerate}[(a)]
\item There is a stable polynomial $P$ and a linear functional $\alpha :
\C_{\rm ma} [z_1 , \ldots , z_d] \to \C$ such that $T (f) = \alpha (f) P$;
\item The polynomial $G_T (z_1 , \ldots , z_d , w_1 , \ldots , w_d)$
is stable as a complex polynomial in $2d$ variables.
\end{enumerate}
\end{thm}

\begin{unremark}
This is a very powerful theorem and it is worth asking 
what part of it one needs to understand in order to 
understand the Borcea--Br{\"a}nd{\'e}n--Liggett theory
of negatively dependent random variables surveyed in
Section~\ref{sec:negdep}.  The answer is that we only
use one direction (sufficiency) and this is only to
prove Proposition~\ref{pr:xy}, which could be proved
elementarily as it concerns only bivariate functions.  
Nevertheless, because this would be messy and unenlightening,
we will prove the sufficiency direction of Theorem~\ref{th:multi-affine PS}
here, which we will then use to prove Proposition~\ref{pr:xy}.
The proof relies on a lemma proved thirty years ago 
by Lieb and Sokal~\cite{lieb-sokal}, which we will not prove 
here (the proof is half a page in~\cite{wagner-BAMS}). 
\end{unremark}
\begin{lem}[\protect{\cite[Lemma~2.3]{lieb-sokal}}] \label{lem:lieb-sokal}
Let $P(\zz) + w Q(\zz)$ be stable in $\C [z_1 , \ldots , z_d , w]$.
If the degree of the variable $z_j$ is at most~1 then the polynomial
$$P(\zz) - \frac{\partial Q(\zz)}{\partial z_j}$$
is either identically zero or is stable.
$\Cox$
\end{lem}

\noindent{\sc Proof of sufficiency:} If $T$ satisfies~$(a)$ the
result is immediate, so assume 
$T$ satisfies~(b).  Because $w \mapsto - w^{-1}$ preserves $\halfplane$,
we see that $G_T (\zz,\ww)$ is stable if and only if $w_1 \cdots w_d
G_T (\zz , - \ww^{-1})$ is stable; here $- \ww^{-1}$ denotes
$(- w_1^{-1} , \ldots , - w_d^{-1})$.  From this, we see that if
$f \in \C [v_1 , \ldots , v_d]$ is stable and multi-affine, then
\begin{equation} \label{eq:lieb-sokal}
( w_1 \cdots w_d) \cdot G_T (\zz , - \ww^{-1}) \cdot f(\vv) = 
   \sum_{S \subseteq \{ 1 , \ldots , n \}} T[\zz^S] (-\ww)^S f(\vv)
\end{equation}
is stable in $\C [\zz, \ww , \vv]$.  We now apply the Lieb-Sokal lemma
$d$ times: the $j^{th}$ time we take $w = w_j$ and $z_j = v_j$.
Each time, the Lieb Sokal lemma replaces $w_j$ in~\eqref{eq:lieb-sokal}
by $- (\partial / \partial v_j)$.  The final result is that
\begin{equation} \label{eq:spec}
\sum_{S \subseteq \{ 1 , \ldots , n \}} T[\zz^S] f^{(S)} (\vv)
\end{equation}
is stable in $\C [ \zz , \vv]$.  Specializing each $v_j$ to
zero preserves stability.  But this specialization yields 
$\sum_S f^{(S)} (0) T[z^S]$, which is another name for $T f(\zz)$.  
Thus stability of $f$ implies stability of $T f$, proving the theorem.
$\Cox$

\subsection{More closure properties} \label{ss:further}

The multivariate P\'olya-Schur theorems of Borcea and Br\"and\'en are 
powerful but there are yet more closure properties that are useful
in probabilistic settings because they have specific probabilistic
meanings.  We state these now but defer the proofs to the next section.

Let $n_j , 1 \leq j \leq d$ be the maximum degree of the variable
$z_j$ in a given polynomial $f \in \C [z_1 , \ldots , z_d]$.
One natural way to create a multi-affine polynomials out of $f$
is to replace each occurrence of each power $z_j^r$ by 
$\binom{n_j}{r}^{-1} e_r (z_{j,1} , \ldots , z_{j,n_j})$,
where $e_r$ is the $r^{th}$ elementary symmetric function.
Here $\{ z_{j,s} : 1 \leq j \leq d , 1 \leq s \leq n_j \}$
are a collection of \Em{clone} variables.  Formally, 
define the \Em{polarization} $f_p$ of $f$ to be the 
unique polynomial symmetric separately in each set of 
clone variables $z_{j,1}, \ldots , z_{j,n_j}$ such that
substituting the value $z_j$ for all the clones $z_{j,s}$
yields the original function $f$.  In Section~\ref{ss:SR proofs}
we will show that stability is closed under polarization.

One final useful closure result is the following.  The proof
will be given after a probabilistic interpretation is given
in Section~\ref{sss:sym}.
\begin{thm}[partial symmetrization] \label{th:partial sym}
Let $f$ be a $d$-variate complex poynomial.  Fix $1 \leq i < j \leq d$  
and define $\tau f$ to be $f$ with the roles of $z_i$ and $z_j$
swapped: 
$$\tau f (z_1 , \ldots , z_d) = f (z_1 , \ldots z_{i-1} , z_j , 
   z_{i+1}, \ldots , z_{j-1} , z_i , z_{j+1} , \ldots , z_d) \, .$$
Fix $\theta \in [0,1]$ and define $f_{\theta; i,j} := (1 - \theta) f 
+ \theta \tau (f)$.  If $f$ is stable, so is $f_{\theta; i , j}$.
\end{thm}

\setcounter{equation}{0}
\section{Negative dependence} \label{sec:negdep} 

\subsection{A brief history of negative dependence}

When random variables in a collection influence each other
in a consistent direction, this can be very useful for obtaining
one-sided bounds on moments and probabilities of natural events.
For example, when all pairwise correlations are negative, the
variance of the sum is bounded above by the sum of the variances.
For this reason, the literature is sprinkled with definitions
of various properties of positive and negative dependence, that
hold in examples of interest and which imply useful consequences. 
Some of these apply to random variables taking values in
$\R, \R^+, \Z$ or $\Z^+$.  To keep matters simple and coherent,
we will discuss only the case where the variables are
\Em{binary valued}, that is, take values in the set $\{ 0 , 1 \}$.
The joint law of $d$ binary variables is a probability measure
on $\B_d := \{ 0 , 1 \}^d$.

The simplest and one of the weakest conditions 
is pairwise positive (respectively negative) correlation: 
$\E X_i X_j \geq (\E X_i) (\E X_j)$ (respectively $\leq$), 
which may also be written as $\P (X_i = X_j = 1) \geq 
\P (X_i = 1) \P (X_j = 1)$ (respectively $\leq$).  As we
have just seen, this condition is strong enough to imply
one-sided bounds on the variance of $\sum_j X_j$.  It
is not strong enough to give bounds on the moments.  For
this, one requires at least positive or negative 
\Em{cylinder dependence}.  Positive cylinder dependence 
is said to hold if 
\begin{equation} \label{eq:cyl 1}
\E \prod_{j \in S} X_j \geq \prod_{j \in S} \E X_j
\end{equation}
for all subsets $S \subseteq \{ 1 , \ldots, n \}$.
Negative cylinder dependence is the reverse inequality.
Sometimes this is strengthened to require as well that
\begin{equation} \label{eq:cyl 2}
\E \prod_{j \in S} (1-X_j) \geq \prod_{j \in S} (1- \E X_j) \, .
\end{equation}
The equation~\eqref{eq:cyl 1} immediately implies that the
moments of $Z := \sum_j X_j$ are at least what they would be
for independent Bernoullis with the same marginals, while
negative cylinder dependence implies the moments are at most
what they would be in the independent case.  Expanding the
moment generating function $\phi (\lambda) := \E e^{\lambda Z}$ 
in powers of $Z$, and then into monomials in the variables $X_j$ 
shows that positive (respectively negative) cylinder dependence
implies that for $\lambda \geq 0$, the moment generating function 
$\phi (\lambda)$ is at least (respectively at most)
what it would be for independent Bernoullis of the same marginal;
including the inequality for complementary cylinders~\eqref{eq:cyl 2} 
extends the inequality for moment generating functions to $\lambda < 0$.
The bounds on moment generating functions transfer to concentration
inequalities (see, e.g.,~\cite[Theorem3.4]{panconesi-srinivasan}).

A yet stronger property is association.  Say that the 
collection of random variables $\{ X_1 , \ldots , X_d \}$ 
is \Em{positively associated} if 
$$\E f(X_1 , \ldots , X_d) g (X_1 , \ldots , X_d) \geq 
   \E f(X_1 , \ldots , X_d) (\E g(X_1 , \ldots , X_d)$$
whenever $f,g : \B_d \to \R$ are nondecreasing functions
(with respect to the coordinatewise partial order on $\B_d$).
Negative association cannot be defined by simply reversing the
inequality because any function $f$ is always nonnegatively
correlated with itself.  Therefore, we amend the definition 
to say that $\{ X_1 , \ldots , X_d \}$ are \Em{negatively
associated} if $\E f g \leq \E f \; \E g$ whenever $f,g : \B_d
\to \R$ are nondecreasing and there is a set $S \subseteq 
\{ 1 , \ldots , n \}$ such that $f$ depends only on 
$\{ X_j : j \in S \}$ while $g$ depends only on 
$\{ X_j : j \notin S \}$.  

One of the most useful results on positive association is
a result of Fortuin, Kasteleyn and Ginibre~\cite{FKG}.
Say that $y$ covers $z$ in the Boolean lattice $\B_d$
if $y > z$ in the coordinatewise partial order 
but there is no element strictly between $y$ and $z$;
we denote this relation by   
$$y \covers z \, ,$$
If $x$ and $y$ both cover $z$ and $w$ covers
$x$ and $y$, we call this configuration a face of the lattice 
$\B_d$.  Say that the probability measure $\mu$ satisfies 
the \Em{positive lattice condition} if whenever $x,y,z,w$ 
form a face of $\B_d$ (with $z$ and $w$ at the bottom and top), 
$\mu (x) \mu (y) \leq \mu (z) \mu (w)$.  
It is immediate to see that this property for faces implies 
$\mu (x \wedge y) \mu (x \vee y) \geq \mu (x) \mu (y)$
for every pair $x,y$, hence we also call in (multiplicative)
submodularity.  What is not so obvious is the theorem now known
as the FKG Theorem.

\begin{thm}[FKG] \label{th:FKG}
If $\mu$ satisfies the positive lattice condition then $\mu$ is
positively associated.
$\Cox$
\end{thm}

The power of this theorem is that its hypotheses are easy to
check while its conclusions are quite strong.  In many statistical
mechanical systems, ratios of probabilities over a face of $\B_d$
are easy to compute even if the absolute probabilities are not.
The FKG theorem immediately implies positive association for
ferro-magnetic Ising models, the random cluster model with
$q \geq 1$, and a number of other models from physics, many 
of which are surveyed in~\cite{liggett-book}.  

For negative dependence the situation is not as nice.  The negative
lattice condition (reverse the inequality on each face) does not
imply negative association, or seemingly anything of value.  Unlike
the positive lattice condition it is not closed under the most
natural of operations, namely integrating out one variable
(probabilistically this means ignoring this variable and viewing 
$\mu$ as a measure on $\B_{n-1}$).  As a result, it has been 
difficult to identify negatively dependent laws and to establish
properties such as negative association.  

\subsection{Search for a theory}

For some time now, a more satisfying theory of negative dependence 
has been sought.  The goal was to find a property of laws $\mu$
on $\B_d$ which would be: 
\begin{itemize}
\item checkable for a handful of known or conjectured \Em{examples};
\item shown to imply \Em{consequences} such as negative association;
\item closed under natural \Em{operations}.
\end{itemize}
This goal was popularized in the article~\cite{pemantle-NA} which 
posed several challenges and stated many conjectures but contained 
few concrete results.  \Em{Spoiler alert:} for binary valued 
random variables, the property is that the probability generating
function is stable; such a probability distribution is called
{\em strong Rayleigh}; see Theorem~\ref{th:SR} below.  In the 
forthcoming lists of examples, consequences and closure properties, 
only in two cases are they not known to hold for the class of
strong Rayleigh measures: Example~\ref{eg:RC} (the random cluster 
measure) is conjectured to be Rayleigh but known not to be 
strong Rayleigh, and the proposed closure property~(\ref{prop:interval}) 
from Section~\ref{sss:closure} (restriction to an interval)
is known to fail.  

One open question remaining from~\cite{pemantle-NA} is whether the 
property known there as h-NLC$^+$ implies negative association.
To define h-NLC$^+$, let us first define an \Em{external field}.
Let $\mu$ be a probability measure on $\B_n$ and let $\blambda_j > 0$ 
be positive real numbers, for $1 \leq j \leq d$.  Define
a measure  $\mu'$ by 
$$\mu' (x_1 , \ldots , x_d) = \frac{\prod_{j=1}^d \lambda_j^{x_j}}{Z} 
   \mu (x_1 , \ldots , x_d)$$
where $Z := \sum_{\xx} \blambda^\xx \mu (\xx)$ is the normalizing
constant.  We call $\mu'$ the perturbation of $\mu$ by the 
external field $\blambda$.  The negative lattice condition 
is not closed under integrating out a variable; we say that the 
hereditary negative lattice condition plus external fields holds 
(h-NLC$^+$) if the negative lattice condition holds for all measures
obtained from $\mu$ by imposing an external field and integrating
out some of the variables.  This was later shown to be equivalent
to Wagner's \Em{Rayleigh condition}~\cite{wagner-rayleigh}, 
which will discuss further in Section~\ref{sub:grail}.  
We now survey the proposed examples, consequences and closure 
properties for the desired class of negatively dependent laws.

\subsubsection{Proposed examples}

\begin{example}[conditioned Bernoullis] \label{eg:bern}
Let $\mu$ be the law of independent Bernoulli variables
$\{ X_1 , \ldots , X_d \}$
with means $p_1 , \ldots , p_d$, let $k$ be an integer
between 1 and $d$ and let $\mu' := (\mu | \sum_{j=1}^d X_j = k)$
be the conditional law of $\{ X_j \}$ given that they 
sum to $k$.  This should be negatively dependent.  Note 
that this is equal to the uniform measure on $k$-subsets
of $\{ 1 , \ldots , n \}$ under the external field
$\lambda_j = p_j / (1-p_j)$.
\end{example}

\begin{example}[spanning trees] \label{eg:UST}
Let $G = (V,E)$ be a finite connected simple graph and 
let $\{ \lambda_e : e \in E \}$ be nonnegative weights.
Let $W$ be the set of subsets $T \subseteq E$ such that the 
resulting $(V,T)$ is connected and has no cycles.  Each $T \in W$
will have cardinality exactly $|V| - 1$.  Such subgraphs are
called \Em{spanning trees} of $G$.  Let $w(T) := \prod_{e \in T}
\lambda_e$ denote the weight of $T$ under the multiplicative
weighting scheme determined by $\blambda$.  We define the
weighted spanning measure $\mu$ (depending on $G$ and $\blambda$) by 
$$\mu (T) = \frac{w(T)}{\sum_{T \in W} w(T)} \, .$$
This was proven in~\cite{pemantle-spanning} to have negative 
correlations and in~\cite{feder-mihail} to have negative
association, and we would hope it to included under any
reasonable definition of ``negatively dependent''.
\end{example}

\begin{example}[random cluster measure] \label{eg:RC}
The random cluster model on a graph $G = (V,E)$ with 
parameters $q > 0$ and $\{ \lambda_e : e \in E \}$ 
is the probability measure $\mu$ on $\{ 0 , 1 \}^E$ 
obtained by normalizing the weights
$$w(T) = \prod_e \lambda_e^{T(e)} q^{N(T)}$$
where $T \subseteq E$ and $N(T)$ is the number of connected 
components of the subgraph $(V,T)$ of $G$.  When $q=2$ this 
is equivalent to the Ising model for ferromagnetism and for 
integers $q \geq 2$ it is equivalent to the Potts model, 
also from statistical physics (see, e.g.~\cite{wu-potts}). 
When $q \geq 1$, the FKG Theorem immediately shows $\mu$ to
be positively associated.  When $0 < q < 1$, the measure is
conjectured but not known to be negatively associated.  We
would hope that the right definition of negative dependence 
would settle this conjecture.
\end{example}

\begin{example}[exclusion measures] \label{eg:exclusion}
The \Em{exclusion process} is a continuous time Markov chain on 
$\B_d$ as follows.  The state $\xx \in B_d$ is interpreted as a 
collection of particles at the sites $j \in \{ 1 , \ldots , d \}$ 
for which $x_j = 1$; the sites $j$ for which $x_j = 0$ are
considered vacant.  Let $\lambda_{\{i,j\}}$ be nonnegative
real numbers, and independently at rate $\lambda_{\{i,j\}}$,
let the values $x_i$ and $x_j$ swap.  This is interpreted
as a particle at one site $i$ or $j$ jumping to the other site.  
If both sites or neither site is occupied, then nothing happens.
The exclusion process is a special case of the more general
\Em{exchange process}, where the labels of the sites are arbitrary,
and in particular might be distinct rather than being drawn from
the two element set $\{ 0 , 1 \}$.  Let $\mu_t$ be the law of the
exclusion process at time $t$, starting from some deterministic
state $\eta$ at time~0.  It was conjectured that $\mu_t$ is always
negatively associated.  We would like the definition of negatively
dependent measures to include $\mu$.
\end{example}

\begin{example}[determinantal measures with Hermitian kernels] 
\label{eg:determinantal}
Let $K$ be a $d \times d$ Hermitian matrix with spectrum in $[0,1]$.
For $S \subseteq \{ 1 , \ldots , d \}$, let $w(S)$ denote the 
determinant of the submatrix $K_S$ of $K$ when restricted to rows and 
columns in the set $S$.  There is a unique probability measure
$\mu$ on subsets of $\{ 1 , \ldots , d \}$ such that
$$\mu \{ T : T \supseteq S \} = |K_S|$$
for each $S$.  This is called the \Em{determinantal} measure
with kernel $K$.  Such measures were proved in~\cite{lyons03}
to be negatively associated.  We would hope these measures
satifsy our definition of negative dependence.
\end{example}

\subsubsection{Proposed consequences} \label{sub:consequences}

It was proposed in~\cite{pemantle-NA} that the right definition
of negative dependence for binary variables would imply the
following properties, with definitions immediately to follow.
\begin{enumerate}[(a)]
\item negative association
\item stochastically increasing levels
\item stochastic covering property
\item log-concave rank sequence
\end{enumerate}

Recall that a probability measure $\mu$ on a lattice $W$
is said to \Em{stochastically dominate} a law $\nu$, written
$\mu \succeq \nu$, if $\mu (A) \geq \nu (A)$ for every
upwardly closed set $A$ (the set $A$ is upwardly closed if
$x \in A$ and $y > x$ implies $y \in A$).  
An equivalent condition is stated in terms of coupling: 
$\mu \succeq \nu$ if and
only if there is a measure $Q$ on $W \times W$ whose
first marginal is $\mu$, whose second marginal is $\nu$,
and which is supported on $\{ (x,y) \in W^2 : x \geq y \}$.
In other words, we may simultaneously sample from $\mu$ and
$\nu$ in such a way that the sample from $\mu$ is always
greater than or equal to the sample from $\nu$.

Given a measure $\mu$ on $\B_d$ and an event $G \subseteq \B_d$
with nonzero measure, denote by $(\mu | G)$ the
measure $\mu (\cdot) / \mu (\cdot \cap G)$ obtained by 
conditioning on $G$.  The conditional measure $(\mu | X_j = 1)$
may be identified with a measure on $\B_{d-1}$ by ignoring
the $j^{th}$ coordinate.  

\begin{property}[stochastically increasing levels]
Say that a measure $\mu$ on $\B_d$ has stochastically increasing
levels if $(\mu | \sum_j x_j = k) \preceq (\mu | \sum_j x_j = k+1)$
for each $k$ such that $\mu (\sum_j x_j = k)$ and 
$\mu (\sum_j x_j = k+1)$ are both nonzero.  
\end{property}

The property of having stochastically increasing levels is not 
implied by negative association, as was hoped in~\cite{pemantle-NA},
but is a desired consequence of the right definition of
negative dependence.

\begin{property}[stochastic covering] \label{prop:SC}
Say that a measure $\nu$ on a lattice $W$ stochastically covers 
the measure $\mu$, and denote this by $\nu \mcovers \mu$, 
if there is a coupling measure $Q$ on $W^2$
with marginals $\nu$ and $\mu$, such that $Q$ is supported
on the set $\{ (x,y) \in W^2 : x = y \mbox{ or } x \covers y \}$.
In other words, $\nu$ stochastically dominates $\mu$ but
``only by one'', in the sense that you can sample from $\mu$
by sampling from $\nu$ and then flipping at most one one to a
zero.
\end{property}

Negative association of a measure $\mu$ implies that 
$(\mu | X_j = 0) \succeq (\mu | X_j = 1)$ when viewed
as a measure on $\B_{d-1}$.  The \Em{stochastic covering
property} is defined in~\cite{PP-rayleigh} to hold when
$(\mu | X_j = 0) \mcovers (\mu | X_j = 1)$, viewed as
measures on $\B_{d-1}$.  This property was wrongly conjectured
in~\cite{pemantle-NA} to follow from negative association.
We would hope for it to follow from the right definition of
negative dependence.

\begin{property}[(ultra) log-concave rank sequence]
The \Em{rank sequence} for the measure $\mu$ on $\B_d$ is the
sequence $(a_k := \mu \{ \sum_j X_j = k \})_{0 \leq k \leq d}$.
It was wrongly conjectured to follow from negative association
that the rank sequence is log-concave.  The right definition of
negative dependence turns out to imply not only log-concavity
but \Em{ultra log-concavity}, namely that for $0 \leq k \leq d$,
equation~\eqref{eq:newton} holds, which I repeat here for convenience.
\begin{equation} \label{eq:ULC}
\frac{a_k^2}{\disp{\binom{d}{k}^2 }} \geq 
   \frac{a_{k-1}}{\disp{\binom{d}{k-1}}} \, 
   \frac{a_{k+1}}{\disp{\binom{d}{k+1}}} \, .
\end{equation}
\end{property}

\subsubsection{Proposed closure properties} \label{sss:closure}

In addition to negative association, it was thought that the
right negative dependence property would imply the following
closure properties.
\begin{enumerate}[(a)]
\item closure under products, projections (integrating out a variable),
and limits;
\item closure under external fields;
\item conditioning on the event $\{ \sum_j X_j = k \}$;  \label{clo}
\item more generally, conditioning on the event $k_1 \leq
\sum_j X_j \leq k_2$; \label{prop:interval} 
\item total symmetrization: replacing $\mu (\xx)$ by $(1/n!)
\sum_{\pi \in S_d} \mu (\pi \xx)$ where $\pi$ acts by permuting
the coordinates;
\item more generally, evolution via the symmetric exclusion process
for any finite time.
\end{enumerate}

It turns out that~\eqref{clo} is too strong, but that the class of
strong Rayleigh measures is indeed closed under \Em{rank rescaling},
(see the upcoming definition) when the rank sequence ${\bf b}$ is 
ultra log-concave.
\begin{defn}[rank rescaling]
Let $\mu$ be a measure on $\B_d$ and let $\bb := (b_0 , \ldots , b_d)$ 
be a vector of nonnegative real numbers.  For $\xx = (x_1 , \ldots , x_d) 
\in \B_d$ let $N (\xx) := \sum_{j=1}^d x_j$ denote the
sum of the coordinates.  For nontriviality assume that
$Z := \sum_{i=0}^d b_i \mu \{ N=i \} > 0$.  The rank rescaling 
of $\mu$ by $\bb$ is the measure $\mu_{\bf b}$ defined by
$$\mu_{\bf b} (\xx) = \frac{N(\xx)}{Z} \, \mu (\xx) \, .$$
\end{defn}

\subsection{The grail is found: application of stability theory 
to joint laws of binary random variables} \label{sub:grail}

It is now acknowledged that the ``correct'' definition of 
negative dependence for binary random variables is the
following definition due to~\cite{borcea-branden-liggett}.
\begin{defn}[strong Rayleigh]
Say that the probability measure $\mu$ on $\B_d$ is 
strongly Rayleigh (also ``strong Rayleigh'') if 
its generating polynomial $f \in \C [ z_1 , \ldots , z_d]$
defined by 
$$f(\zz) := \sum_{\omega \in \B_d} \mu (\omega) \zz^\omega$$
is a real stable polynomial.
\end{defn}

As a generating polynomial, $f$ is automatically multi-affine,
with real nonnegative coefficients.  Recalling 
Theorem~\ref{th:mixed partials}, we see that strong
Rayleigh is equivalent to the inequality~\eqref{eq:rayleigh}
on the mixed partial derivatives of $f$ holding for all
real $\xx$.  This is sometimes taken as the definition
of strong Rayleigh.  

The origin of the name comes from Wagner's Rayleigh property.  
This is the mixed partial inequality~\eqref{eq:rayleigh} but
only required to hold for nonnegative real arguments.
As was mentioned before, this is equivalent to h-NLC$^+$, which
is not hard to believe because varying $\xx$ over $\R_+^d$
corresponds to varying the external field arbitrarily.
Strangely, when this is required to hold for negative values
of the external field as well, the property becomes more
natural and robust.  Because the Rayleigh property is 
equivalent to h-NLC$^+$ and is obviously a consequence
of the strong Rayleigh property, we see immediately that
strong Rayleigh implies pairwise negative correlations.
This is, of course, just the tip of the iceberg.
\begin{thm}[properties of strong Rayleigh measures] \label{th:SR}
~\\ \vspace{-0.2in}

$(i)$ Strong Rayleigh measures are negatively associated.

$(ii)$ Strong Rayleigh measures have stochastically 
increasing levels, satisfy the stochastic covering 
property~\ref{prop:SC}, and have ultra log-concave 
rank sequences as in~\eqref{eq:ULC}.  The class of 
strong Rayleigh measures is also closed under 
rank rescaling by sequences $(b_0 , \ldots , b_d)$ 
when the nonzero values $b_i$ are the coefficient
sequence of a stable univariate polynomial $\sum_{i=0}^d
b_i z^i$.  In particular, this implies 
closure of the class of strong Rayleigh measures under conditioning 
on $\{ \sum_j X_j = k \}$ or $\{ k \leq \sum_j X_j \leq  k+1 \}$.

$(iii)$ The class of strong Rayleigh
measures is closed under products, projections, external fields, 
total symmetrization and the symmetric exclusion dynamics.  
\end{thm}

As an example of the utility of this theorem, here is one
of the two applications that brought the theory to my
attention.
\begin{example}[concentration inequalities for Lipschitz functions]
Let $\mu$ be a measure on $\B_d$ and let $f : \B_d \to \R$ 
be Lipschitz with constant~1.  When $\mu$ is a product measure, 
the martingale $W_k := \E (f \| \omega_1 , \ldots , \omega_k)$ 
has bounded increments $|W_k - W_{k-1}| \leq 1$.  
It follows from Azuma's inequality (see, e.g,~\cite{alon-spencer})
that $f$ satisfies the Gaussian concentration inequality 
\begin{equation} \label{eq:concentration}
\P (f - \E f \geq a) \leq e^{- a^2 / (2n)} \, .
\end{equation}
In fact independence is not needed.  The stochastic covering 
property was defined in~\cite{PP-rayleigh} precisely to 
imply $|W_k - W_{k-1}| \leq 1$.  It follows 
that~\eqref{eq:concentration} holds for any strong Rayleigh 
measure.
\end{example}

\subsection{Proof of Theorem~\protect{\ref{th:SR}}}
\label{ss:SR proofs}

\subsubsection{Proof of part~$(i)$~: negative association}

``Strong Rayleigh implies negative association''
was the result of Borcea--Br\"and\'en--Liggett
that pointed the way to application of
stable polynomial theory to the probabilistic setting.
We follow their proof, which begins by paving the way for 
a reduction to the case of homogeneous measures.

\begin{lem}[homogenization] \label{lem:homogenization}
Let $\{ X_1 , \ldots , X_d \}$ be a collection
of random variables whose law is strong Rayleigh 
and define $X_{d+1} := d - \sum_{j=1}^d X_j$.
Then the collection $\{ X_1 , \ldots , X_{d+1} \}$ 
is strong Rayleigh.
\end{lem}

\noindent{\sc Proof:} The generating function for the collection
$\{ X_1 , \ldots , X_{d+1} \}$ is just the homogenization of the
generating function for $\{ X_1 , \ldots , X_d \}$.  Therefore, 
what we need to check is that the homogenization of a 
stable generating function is stable.  This is not true 
for all polynomials\footnote{Indeed $f(x) = x-1$ is a trivial example
of a stable polynomial whose homogenization $f_H (x,y) = x-y$
fails to be stable.} but for probability generating functions,
or more generally any stable polynomial with nonnegative
coefficients, this follows from~$(i) \Rightarrow (ii)$ of 
Proposition~\ref{pr:homogenization}.  
$\Cox$

Next we recall the definition of polarization via clone variables
from Section~\ref{ss:further}: it is the unique polynomial 
symmetric separately in each set of clone variables 
$\{z_{j,1}, \ldots , z_{j,n_j} \}$ such that
substituting in $z_j$ for all the clones $z_{j,s}$
yields the original function $f$.  The following
result does not rely on $f$ having nonnegative, or
even real coefficients.
\begin{thm}[polarization] \label{th:polarization}
The complex polynomial $f$ is stable if and only if
its polarization $f_p$ is stable.
\end{thm}

To prove this we require the Grace--Walsh--Szeg{\H{o}} theorem.
This result may be derived from the theory of stable polynomials, 
as shown in~\cite{borcea-branden-LYPS2} and streamlined
in~\cite[Section~4]{wagner-BAMS}.  We will be content here
to quote this century old result; for a proof, 
see~\cite[Theorem~2.12]{COSW} or Section~4 of~\cite{wagner-BAMS}.
\begin{lem}[GraceWalsh-Szeg{\H{o}} Theorem] \label{lem:GWS}
Let $f \in \C [z_1 , \ldots , z_d]$ be symmetric and multi-affine
and let $\CC$ be a convex circular region containing the points 
$\zeta_1 , \ldots , \zeta_d$, for instance the upper half-plane 
$\halfplane$.  Then there exists at least one point $\zeta \in \CC$ 
such that $f(\zeta_1 , \ldots , \zeta_d) = f(\zeta , \ldots , \zeta)$.
$\Cox$
\end{lem}

\noindent{\sc Proof of Theorem}~\ref{th:polarization}: 
One direction is elementary: if $f_p$ is stable then the diagonalization
property allows us to set $z_{j,s} = z_j$ for all $j$ and we deduce
stability of $f$.  For the nontrivial direction, suppose $f_p$ is
not stable and let ${\bf a} := \{ a_{j,s} \}$ be numbers in the 
upper half-plane such that $f_p ({\bf a}) = 0$.  Fixing the numbers
$\{ a_{j,s} : j > 1 \}$, the Grace-Walsh-Szeg{\H{o}} Theorem 
implies the existence of a number $a_1$ in the upper half-plane
such that $f(a_1 , \ldots , a_1 , {\bf a}') = 0$, where
${\bf a'}$ are the numbers $a_{j,s}$ for $j > 1$.  Iterating,
we arrive at numbers $a_2, \ldots , a_n \in \halfplane$ for
which $f(a_1 , \ldots , a_n) = 0$, showing that $f$ is not stable.
$\Cox$

A measure $\mu$ on $\B_d$ is called \Em{Projected Homogeneous Rayleigh}
(PHR) if there is a measure $\nu$ on $\B_m$ for some $m \geq d$
such that $\nu$ is a homogeneous measure with the Rayleigh property
and the projection of $\nu$ to the law of the first $d$ variables 
gives $\mu$.  If, furthermore, the law $\nu$ is strong Rayleigh, 
we say that $\mu$ is PHSR.   

\begin{cor} \label{cor:PHR}
Strong Rayleigh measures are PHSR.
\end{cor}

\noindent{\sc Proof:} First extend the law $\mu$ of the binary 
variables $\{ X_1 , \ldots , X_d \}$ by homogenization to a law 
$\mu_H$ on variables $\{ X_1 , \ldots , X_{d+1} \}$, the last of 
which takes integer values.  By Lemma~\ref{lem:homogenization}
this is still in the class of strong Rayleigh measures.  Next replace
the integer variable $X_{d+1}$ by clones $\{ X_{d+1,1} , \ldots , 
X_{d+1,d} \}$ by polarization; the resulting measure $\nu$ on
$\B_{2d}$ is strong Rayleigh by Theorem~\ref{th:polarization}
and its projection onto the first $d$ coordinates is $\mu$.
$\Cox$

\noindent{\sc Proof of negative association:} 
To prove negative association, we use an argument due to 
Feder and Mihail~\cite{feder-mihail}.  This argument shows
that homogeneous strong Rayleigh measures (thus by 
Corollary~\ref{cor:PHR}, all strong Rayleigh measures) are 
negatively associated.  All that is needed in the Feder--Mihail
lemma is for this class to be closed under conditioning on 
$\{ X_j = x \}$ (which is obvious) and have pairwise negative 
correlations, which we have seen in~\eqref{eq:rayleigh}.  In
short, negative association for strong Rayleigh measures is 
reduced to the following lemma.

\begin{lem}[Feder--Mihail] \label{lem:FM}
Let $\SS$ be a class of homogeneous measures on finite 
Boolean algebras (of differing sizes) which is closed
uner conditioning and each of which has pairwise
negative correlations.  Then all measures in $\SS$
are negatively associated.
\end{lem}

Assuming this lemma, the proof of negative association
finishes as follows.  The class of homogeneous strong
Rayleigh measures is closed under conditioning on the
value of any variable.  We have noted that strong Rayleigh
implies negative pairwise correlations, so we conclude 
from the Feder--Mihail lemma that homogeneous strong
Rayleigh measures are negatively associated.  Any strong
Rayleigh measure is a projection of a homogeneous strong
Rayleigh measure, hence inherits the negative association
property.  It remains only to prove the lemma.

\noindent{\sc Proof of Feder--Mihail Lemma:} 
Any increasing function on $\B_d$
is a positive linear combination of increasing indicator
functions, that is, indicator functions $\one_A$ of 
upwardly closed events $A$.  It therefore suffices to
prove negative correlation for upwardly closed events 
$A$ and $B$ depending on disjoint sets of variables.  
(Note: throughout the proof, the phrase ``negatively
correlated'' includes the case of zero correlation.)

We first prove this in the special case where $A$ is the
event $\{ X_j = 1 \}$.  We use induction on the number 
of variables, $m$.  When $m=2$, the only nontrivial
case is $\one_A = X_1$ and $\one_B = X_2$, which are
negatively correlated by hypothesis.  Now assume for 
induction that for all measures $\nu \in \SS$ on Boolean
lattices $\B_d$ for $d \leq m$, for all $j$, and for 
all upwardly closed events $B$ not depending on $X_j$, 
the functions $X_j$ and $\one_B$ are negatively correlated.  
Let $\mu$ be a law on $\B_{m+1}$, and let $j \leq m+1$
and $B \subseteq \B_{m+1}$ not depending on $X_j$ be given.
We will use the inequality 
$$p a + (1-p) b \leq p' c + (1-p') d$$
holding when $p, p', a,b,c,d \in [0,1]$ with
$p \leq p', a \leq c, b \leq d$ and $a \geq b$.
To use this inequality, fix $i \leq m+1$ to be
determined later.  Let $p := \mu (X_i = 1 \| X_j = 1)$ 
and $p' := \mu (X_i = 1 \| X_j = 0)$, the inequality
$p \leq p'$ holds because it is equivalent to pairwise 
negative correlation of $X_i$ and $X_j$.  We let
\begin{eqnarray*}
a & := & \mu (B \| X_i = 1 , X_j = 1) \, ; \\
b & := & \mu (B \| X_i = 0 , X_j = 1) \, ; \\
c & := & \mu (B \| X_i = 1 , X_j = 0) \, ; \\
d & := & \mu (B \| X_i = 0 , X_j = 0) \, .
\end{eqnarray*}
The inequalities $a \leq c$ and $b \leq d$ follow
because the conditional measures $(\mu \| X_j = 1)$,
$(\mu \| X_i = 1)$ and $(\mu \| X_i = 0)$ are all
subject to the induction hypothesis.  Finally, to 
ensure that $a \geq b$, we now choose $i$ judiciously.
Equivalent conditions for the inequality $a \geq b$ are
\begin{eqnarray*}
   \frac{\mu (B, X_i = 1, X_j = 1)}{\mu (B^c , X_i = 1 , X_j = 1)} & \geq &
   \frac{\mu (B, X_i = 0, X_j = 1)}{\mu (B^c , X_i = 0 , X_j = 1)}  \\
   & \Updownarrow & \\ 
   \frac{\mu (B, X_i = 1, X_j = 1)}{\mu (B, X_i = 0, X_j = 1)} & \geq &
   \frac{\mu (B^c , X_i = 1 , X_j = 1)}{\mu (B^c , X_i = 0 , X_j = 1)}  \\
   & \Updownarrow & \\ 
\E (X_i \| B \, , X_j = 1) & \geq & \E (X_i \| B^c \, , X_j = 1) \, .
\end{eqnarray*}
Homogeneity of the measure $\mu$ implies that the sum over $i$ of
the left and right-hand sides of the last inequality are both 
equal to the deterministic value $\sum_{i=1}^{m+1} X_i$.  Therefore
the inequality must hold for at least one $i \leq m+1$.  
We conclude that
$$\mu (B \| X_j = 1) = p a + (1-p) b \leq p' c + (1-p') d 
   = \mu (B \| X_j = 0) \, .$$

Having established that single variables are negatively correlated
with upwardly closed events for all measures in $\SS$, we consider 
the general case where $A$ and $B$ are upwardly closed events 
depending on disjoint sets of variables.  Again we use induction,
assuming the result for all measures on fewer variables and also
already possessing the result in the special case.  We run the 
argument with $\one_A$ in place of $X_j$, thus $p := \mu (X_i = 1 \| A)$,
$p' := \mu (X_i = 1 \| A^c)$, and so forth.  The inequality $p \leq p'$
is a consequence of the special case.  The inequalities $a \leq c$
and $b \leq d$ come from the induction hypothesis applied to
the conditional measures $(\mu \| X_i = 1)$ and $(\mu \| X_i = 0)$
respectively.  Again, constancy of $\sum_{i=1}^{m+1} X_i$ 
forces the inequality $\E (X_i \| B , A) \geq \E (X_i \| B^c , A)$
to hold for some $i$, finishing the proof.
$\Cox$

\subsubsection{Proof of part~$(ii)$~: stochastic inequalities}

We begin with the stochastic covering property, which follows
directly from Corollary~\ref{cor:PHR} and the following lemma,
first proved in~\cite{PP-rayleigh}.
\begin{lem}[\protect{\cite[Proposition~2.2]{PP-rayleigh}}] 
\label{lem:PHR -> SC}
If a measure is projected homogeneous strong Rayleigh then it has the 
stochastic covering property.
\end{lem}

\noindent{\sc Proof:} Let $\mu$ be a projection of a homogeneous
strong Rayleigh measure $\nu$.  Without loss of generality, we assume
that the projection is onto the last $d$ out of $m$ coordinates.
Let $\nu_0$ denote the probability measure on $\B_{m-1}$ which is
the $\nu$-law of $(X_1 , \ldots , X_{m-1})$ conditioned on $X_m = 0$;
similarly, let $\nu_1$ denote the probability measure on $\B_{m-1}$ 
which is the $\nu$-law of $(X_1 , \ldots , X_{m-1})$ conditioned on 
$X_m = 1$.  Negative association of $\nu$ implies that $\nu_0$
stochastically dominates $\nu_1$: let $B \subseteq \B_m$ be any
event not depending on the $m^{th}$ coordinate then $B$, that is
$B = B_* \times \{ 0 , 1 \}$ for some $B_* \in \B_{m-1}$; 
then negative correlation of $B$ and
$X_m$ implies that $\nu_0 (B_*) \geq \nu_1 (B_*)$.  The equivalent
coupling formulation of stochastic domination (see the beginning of
Subsection~\ref{sub:consequences}) shows there is a random variable
$(X,Y)$ such that $X$ has law $\nu_0$, $Y$ has law $\nu_1$ and 
$X \geq Y$.  By homogeneity of $\nu$, we see in fact that $X$ always
covers $Y$.  Projecting back to $\B_d$, this yields a random pair
$(\overline{X}, \overline{Y})$ such that $\overline{X}$ has law
$(\mu \| X_m = 0)$, $\overline{Y}$ has law $(\mu \| X_m = 1)$, 
and $\overline{X}$ is always equal to or covering $\overline{Y}$.
$\Cox$

\noindent{\sc Proof of ultra log-concavity of the rank sequence:}
Next we examine the rank function.  Let $\mu$ be a strong
Rayleigh measure on $\B_d$.  The rank sequence $a_k := \mu (N = k)$
has generating function $g(z) := \sum_{k=0}^d a_k z^k =  
f_\mu (x, \ldots , x)$.  This is a stable polynomial
because it is a diagonalization of the stable polynomial $f_\mu$,
recalling property~(d) at the beginning of Section~\ref{sec:stability 2}.
Newton's inequalities (Theorem~\ref{th:newton}) for univariate
real stable polynomials then yield ultra log-concavity of the 
coefficients of the rank sequence. 
$\Cox$

\noindent{\sc Proof of rescaling by a stable coefficient sequence:}
Let $b_0 , \ldots , b_d$ be a sequence of nonnegative real
numbers whose nonzero elements are an interval $b_r , b_{r+1},
\ldots , b_s$ and are the coefficients of a stable polynoimal.
We have seen in Example~\ref{eg:PS3} that the sequence 
$\{ b_k \}$ is a multiplier sequence.  It follows from
Theorem~\ref{th:mv multiplier} that if $f (\zz) = \sum_\rr 
a_\rr \zz^\rr$ is a real stable function of $d+1$ variables 
then $g(z) := \sum_\rr a_\rr b_{r_{d+1}} \zz^\rr$ is also stable.

Now let $\mu$ be a strong Rayleigh measure on $\B_d$ 
with probability generating function $f_\mu$.  By 
Lemma~\ref{lem:homogenization}, the measure $\nu$ 
whose generating function is the homogenization
$f_\nu = (f_\mu)_H$ is also stable (recall that $\nu$
is the law of $(X_1 , \ldots , X_d , d - \sum_{j=1}^d X_j)$ 
when $(X_1 , \ldots , X_d)$ has law $\mu$).  With $\{ b_k :
0 \leq k \leq d \}$ any ultra log-concave sequence as above, 
we have seen that $g(\zz) := \sum_\rr a_\rr b_{r_{d+1}} \zz^\rr$ 
is stable.  Setting $z_{d+1}$ equal to~1, we see that $g(\zz) /  
g(1 , \ldots , 1)$ is a real stable probability generating function, 
whence the measure described at the end of part~$(ii)$ of
Theorem~\ref{th:SR} is strong Rayleigh.
$\Cox$

\noindent{\sc Proof of stochastically increasing levels:}
If the interval of support $[r,s]$ of the ultra log-concave
sequence ${\bf b}$ satisfies $s-r = 0$ or~1, then it is automatically
ultra log-concave.  Taking $s=r$, we see that if the law
$\mu$ of $(X_1 , \ldots , X_d)$ is strong Rayleigh then 
so is $(\mu \| N = k)$ as long as $\mu (N=k) > 0$; here 
again $N := \sum_{j=1}^d X_j$.  Taking $s = r+1$ we see 
that $(\mu \| k \leq N \leq k+1)$ is strong Rayleigh, 
provided that the event conditioned on is not null.  

Fix $k$ such that $\mu (N=k)$ and $\mu (N=k+1)$ are both 
nonzero, and let $\mu_k := (\mu \| N=k)$ and $\mu_{k+1} := 
(\mu \| N = k+1)$.  The measure $(\mu \| k \leq N \leq k+1)$
is strong Rayleigh.  Homogenizing by adding the check bit
$x_{d+1} := k + 1 - \sum_{j=1}^d x_j$ gives a measure $\nu$ 
on $\B_{d+1}$ that is also strong Rayleigh.  Let $A \subseteq \B_d$ 
be any upwardly closed set and let $A_* := A \times \{ 0,1 \}$.  By
negative association, the event $A_*$ is negatively correlated
with $X_{d+1}$.  Thus $\nu (A_* \| X_{d+1} = 1) \leq 
\nu (A_* \| X_{d+1} = 0)$, which is equivalent to 
$\mu_k (A) \leq \mu_{k+1} (A)$.  Because $A$ was an
arbitrary upwardly closed event, this proves that
$\mu_k \preceq \mu_{k+1}$.
$\Cox$

\subsubsection{Proof of part~$(iii)$~: exclusion and other 
closure properties}  \label{sss:sym}

Most of the closure properties in part~$(iii)$ of Theorem~\ref{th:SR}
are easy, following from the closure properties of the class of
stable polynomials listed at the beginning of Section~\ref{sec:stability 2}.
The generating function for the product measure $f_{\mu \times \nu}$
is the product of functions $f_\mu \cdot f_\nu$, whence the
strong Rayleigh property for products of strong Rayleigh measures 
follows from stability of the product of stable polynomials. 
Projection corresponds to setting some variables equal to 1, whence
closure under projections follows from specialization to real values
(Lemma~\ref{lem:closure}).  External fields correspond to 
replacing $f$ by $c f(a_1 z_1, \ldots , a_d z_d)$ for positive
real values of the $a_j$, which follows from the closure of 
stable polynomials under dilations.  To prove closure of the
class of strong Rayleigh measures under total symmetrization,
there are two avenues.  One is to deduce this from the more
difficult result for partial symmetrization (see below), 
observing that total symmetrization may be achieved by 
repeated partial symmetrization.  
As pointed out in~\cite[Remark~4.5]{borcea-branden-liggett}, 
a more direct argument is as follows.  Let $\mu$ be a
strong Rayleigh measure on $\B_d$.  The sum $N := \sum_{j=1}^d X_j$ 
has generating function $f_N (z) := f_\mu (z , \ldots , z)$, 
which is stable because it is a diagonalization of $f_\mu$.  The 
generating function of the total symmetrization $\overline{\mu}$ 
is the polarization of $f_N (z)$, hence is stable by
Theorem~\ref{th:polarization}.  

The key to proving closure of the class of strong Rayleigh 
measures under exclusion dynamics is to prove that a single 
step of \Em{partial symmetrization} preserves the strong Rayleigh 
property.  Let $\mu$ be a probability measure on $\B_d$, let
$1 \leq i < j \leq d$ be indices and fix $0 < \theta < 1$.
The $\theta$-partial symmetrization of $\mu$ with respect to 
indices $i$ and $j$ is defined to be the measure $\mu_{\theta; i,j}$
whose verbal description is as follows.
\begin{quote}
To sample from $\mu_{\theta; i,j}$, first sample from $\mu$, 
then flip an independent $\theta$-coin to decide whether to 
transpose the $i$ and $j$ coordinates.  
\end{quote}

At the level of generating functions, this is 
Theorem~\ref{th:partial sym}, which we now recall
and prove.  The generating polynomial $f_{\mu, \theta ; i , j}$ is
the one described in Theorem~\ref{th:partial sym} by
$$f_{\mu, \theta; i,j} = (1 - \theta) f_\mu + \theta \tau (f_\mu)$$
where $\tau$ operates by switching arguments $i$ and $j$.
The theorem states that if $f$ is stable then so is 
$(1 - \theta) f + \theta \tau (f)$, which when applied to
$f = f_\mu$ yields closure of strong Rayleigh measures under
partial symmetrization.  We now prove the theorem.

\noindent{\sc Proof of Theorem~\ref{th:partial sym}:}
Let $\mu$ be a strong Rayleigh measure.
Assume without loss of generality that $i=1$ and $j=2$.  We need
to show that $f_{\mu, \theta , 1 , 2}$ is stable.  For every
choice of complex numbers $a_3 , a_4 , \ldots , a_d \in \halfplane$,
we know that $f_\mu (\cdot , \cdot , a_3 , \ldots , a_d)$ is 
multi-affine and stable and we need to show that $f_{\mu,\theta,1,2} 
(\cdot , \cdot , a_3 , \ldots , a_d)$ is stable.  Therefore it 
suffices to show: 
\begin{pr} \label{pr:xy}
For all multi-affine stable functions $f : \C^2 \to \C$,
the function $(1-\theta) f (x,y) + \theta f(y,x)$ is stable.
\end{pr}
\begin{unremark}
Because of the restriction to the multi-affine case, we know
that $f(x,y) = a + bx + cy + dxy$, so it cannot be too hard to
prove this!  We observe that if we needed this only for real 
stable functions, it would follow immediately from the characterization
that $f$ is stable if and only if $ad \leq bc$, which follows from
the mixed partial derivative criterion~\eqref{eq:rayleigh}.  Indeed,
the proof of Theorem~\ref{th:partial sym} appearing 
in~\cite[Section~4.4]{borcea-branden-liggett} begins with
this and extends via the multivariate Obreschkoff theorem.
To keep things self-contained, we will instead derive this
from the portion of Theorem~\ref{th:multi-affine PS} that we have
proved; another self-contained proof is given 
in~\cite[Theorem~7]{liggett-exclusion}; and exhaustive description 
of this four parameter family is yet another way to finish this.
\end{unremark}

\noindent{\sc Proof of Proposition}~\ref{pr:xy}:
The operator $T$ on $\C_{\rm ma} [x,y]$ 
corresponding to partial symmetrization is defined by its action
on the four basis elements $1, x, y$ and $xy$:
$$T(1) = 1; \; T(x) = (1 - \theta) x + \theta y; \; 
   T(y) = (1 - \theta) y + \theta x; \; T(xy) = xy \, .$$
Plugging into the definition~\eqref{eq:GT} of $G_T$ gives
$$G_T (x,y,u,v) = xy + [(1 - \theta) x + \theta y ] u
   + [(1 - \theta) y + \theta x ] v + uv \, .$$
We need to show this is stable.  It is multi-affine and
has real coefficients, so we may apply the mixed partial
derivative test~\eqref{eq:rayleigh}.  There are six unordered
pairs of variables.  In each case we easily verify that the
left side of~\eqref{eq:rayleigh} minus the right side is
equal to a nonnegative multiple of a square when 
$0 \leq \theta \leq 1$.  For example, for the pair $(x,u)$,
we evaluate $( \partial^2 f / \partial x \partial u) \, f  
- (\partial f / \partial x) (\partial f / \partial u)$ 
to obtain $\theta (y+v)^2$, while for the pair $(x,y)$ we obtain 
$(\partial^2 f / \partial x \partial y) \, f 
- (\partial f / \partial x) (\partial f / \partial y)
= \theta (1-\theta) (u-v)^2$; up to the symmetries of $f$, 
these are the only two cases, therefore the proposition
is proved.
$\Cox$

The final step, from partial symmetrization to exclusion dynamics,
is a small one.  The operators $T_{\theta ; i , j}$ defined by
$\mu \mapsto \mu_{\theta; i,j}$ generated a semigroup $S$ of 
operators on the space of probability measures on $\B_d$.
For any set $\{ \lambda_{\{i,j\}} \}$ of swap rates and
any $t \geq 0$, the operator mapping $\mu$ to the time $t$
law of the exclusion process with swap rates $\lambda$ started
from the measure $\mu$ is in the closure of $S$.  This is 
more or less self-evident; for a more formal proof, note that
when only one swap rate is nonzero, the exclusion evolution
operator is already equal to an operator $T_{\theta; i,j}$,
after which one can use the Trotter product formula 
(see~\cite[Proposition~5.1]{borcea-branden-liggett}) to 
write the general exclusion evolution operator as a limit of
products of those with only one nonvanishing swap rate.  This
completes the proof of Theorem~\ref{th:SR}.
$\Cox$

To finish the discussion of probabilistic applications, we 
discuss examples.  The theorem explicitly addresses exclusion 
measures (Example~\ref{eg:exclusion}).  Conditioned
Bernoullis (Example~\ref{eg:bern}) are strong Rayleigh due
stability of $1-p+pz$ for $0 \leq p \leq 1$, closure under 
products, (yielding all measures with independent coordinates), 
and closure under conditioning on $N=k$.  Random cluster
measures with $q < 1$ are conjectured to be Rayleigh and known
not to be strong Rayleigh and spanning trees are a special case of 
determinantal measures, so among the
examples~\ref{eg:bern}--\ref{eg:determinantal} it remains
only to show that determinantal measures are strong Rayleigh.
We will require a result which is the analogue of
Example~\ref{eg:hermitian} but for stability instead of
hyperbolicity.  The exposition is taken 
from~\cite[Proposition~3.2]{borcea-branden-liggett} though
the result itself has been know for a long time.

\begin{pr} \label{pr:hermitian}
Let $A_1 , \ldots A_d$ be (complex) positive semi-definite
$n \times n$ matrices and let $B$ be a Hermitian matrix,
also $n \times n$.  
\begin{enumerate}[1.]
\item The polynomial
$$(z_1 , \ldots , z_d) \mapsto f(z_1 , \ldots , z_d) 
   := \Det \left ( z_1 A_1 + \cdots + z_d A_d + B \right )$$
is either identically zero or it is real stable.
\item If $B$ is also positive semi-definite then the $f$ has
all nonnegative coefficients.  
\item It follows that if $Z := \diag (z_1 , \ldots , z_d)$
and $B$ is any positive semi-definite $d \times d$ matrix, 
then $\Det (B + Z)$ is a multi-affine real stable polynomial
with all nonnegative coefficients, hence equal to $c f_\mu$
for some strong Rayleigh measure $\mu$.
\end{enumerate}
\end{pr}

\noindent{\sc Proof:} Stability is closed under limits,
so it suffices to prove the proposition in the case where
the matrices $A_1 , \ldots , A_d$ are positive definite.
The polynomial $f$ is real on real inputs because the 
determinant of a complex Hermitian matrix is real.
Pick a vector $\uu \in \R^d_+$ and $\vv \in \R^d$.  Define
a function $\zz$ by $\zz (t) := \vv + t \uu$.  The matrix
$P := \sum_{j=1}^d v_j A_j$ is positive definite, hence 
invertible with a positive definite square root, call it $Q$.
We may write
$$f(\zz (t)) = \Det (P) \; \Det \left ( t I + Q H Q^* \right )$$
where $H := B + \sum_{j=1}^d u_j A_j$ is Hermitian.  Thus
$f(\zz(t))$ is a polynomial in $t$ that is a constant multiple
of the characteristic polynomial of a Hermitian matrix,
hence has all real zeros.  Because $\vv \in \R^d$ and 
$\uu \in \R^+_d$ are arbitrary, the criterion for stability
in Proposition~\ref{pr:inhom orthant} is satisfied, proving
the first statement of the theorem.

Expanding the determinant, we find that the coefficients of $f$ 
are products of principal minors of the matrices $A_1 , \ldots , A_d$
and $B$.  When these are all positive semi-definite, the principal
minors are positive, proving the second statement.  The last statement
follows immediately from setting $A_i$ to the matrix with a~1 in the
$(i,i)$-entry and zeros elsewhere, noting that this matrix is 
positive semi-definite.
$\Cox$

\noindent{\sc Proof that determinantal measures are strong Rayleigh:}
Assume first that the Hermitian kernel $K$ is invertible.
It is easily seen that the generating polynomial $f_\mu$ is given by
$$f_\mu (z_1 , \ldots , z_d) = \Det (I - K + K Z) 
   = \Det (K) \; \Det (K^{-1} - I + Z) \, ,$$
where $Z := \diag (z_1 , \ldots , z_d)$.  Because $K$ is also a 
contraction, $K^{-1} - I$ is positive semi-definite.  The constant
$\Det (K)$ is positive, so it follows from Example~\ref{eg:hermitian}
that $f_\mu$ is stable.  The general case may be obtained by
taking limits because positive definite kernels are dense 
in the space of positive semi-definite kernels.
$\Cox$

\setcounter{equation}{0}
\section{Further applications of stability: determinants, permanents 
and moments} \label{sec:determinants} 

From the outset, determinants have been prominent in the theory of
hyperbolic polynomials.  Already in~\cite[Example~2]{garding-hyperbolic}, 
the determinant function on the space of Hermitian matrices was
given as an example of a hyperbolic polynomial, with the 
nonnegative definite matrices being a cone of hyperbolicity
(see Example~\ref{eg:hermitian} above).  Related to the determinant
but more enigmatic is the \Em{permanent}.  The definition of the
permanent,
$$\Per (A) := \sum_{\sigma \in S_n} \prod_{i=1}^n
   a_{i, \sigma (i)} \, ,$$
differs from that of the determinant only in that there is
no alternating sign factor $(-1)^{{\rm parity}\, (\sigma)}$
in the summand.  This makes the permanent much less tractable
than the determinant.  For example, it is \#{P}-hard to compute
the permanent of a zero-one matrix, while determinants may be evaluated
in polynomial time.  In this section we discuss two results on
permanents and one on determinants.  The first of these is 
a lower bound on the permanent, conjectured by van der Waerden
in 1926, proved independently by Egorychev and Falikman in 1981,
and re-proved in a simpler and more general way by Gurvits in
2009 using the theory of stable functions.  The second of
these results, the Monotone Column Permanent Conjecture
asserts the stability of a certain polynomial obtained as
a permanent.  It was conjectured in 1999 and proved in 2009.
The third result, the so-called BMV conjecture, is only tangentially
related to stability theory, but has garnered enough attention
to mandate its inclusion here.  It was conjectured in 1975.
A proof posted recently to the arXiv is believed to be correct.

\subsection{The van der Waerden conjecture}

One prolific area of research in understanding the permanent has
been to identify extremal cases for various families of matrices.
Some motivation for understanding permanents comes from 
graph theory.  If $A$ is the incidence matrix of a bipartite
graph $G$, then $\Per (A)$ is the number of perfect matchings
of $G$.  This interpretation has led to an emphasis on extrema
over zero-one matrices.  For example, one might consider matrices of
zeros and ones with prescribed row sums.  An upper bound, conjectured
by Minc and proved by Br{\'e}gman~\cite{bregman} is as follows.
\begin{thm}[Br{\'e}gman's Theorem] \label{th:bregman}
Let $A$ be a nonnegative $n \times n$ matrix of zeros and ones
and let $r_1 , \ldots , r_n$ denote the row sums of $A$.  Then
$$\Per (A) \leq \prod_{j=1}^n (r_j !)^{1/r_j} \, .$$
$\Cox$
\end{thm}

The lower bound is of course zero because there could be a column
of zeros.  If we prescribe the column sums as well as the row sums,
can we find a nontrivial lower bound?  Removing the restriction 
to zero-one matrices, this question was posed by van der Waerden in
1926 for \Em{doubly stochastic} matrices.
A matrix $A = (a_{ij})_{1 \leq i,j \leq n}$ is said to be
stochastic if it has nonnegative entries and all row sums 
$\sum_j a_{ij}$ are equal to~1.  The terminology comes from 
the fact that these matrices are precisely the transition kernels 
for Markov chains.  The matrix $A$ is said to be doubly stochastic
if all column sums $\sum_i a_{ij}$ are equal to~1 as well.  One
might not see at first that the permanent of a doubly stochastic
matrix must be nonzero, but this follows from the well known
fact that doubly stochastic matrices are positive linear 
combinations of permutation matrices (and such matrices have
permanent equal to~1).  

Intuition may suggest that the minimum occurs when all entries
are equal to $1/n$.  Indeed, van der Waerden conjectured in 1926 
that if $A$ is doubly stochastic then
\begin{equation} \label{eq:vdW}
\Per (A) \geq \frac{n!}{n^n}
\end{equation}
with equality if and only if $a_{ij} = 1/n$ for all $i,j$.
This theorem was proved thirty years ago, independently by
Egorychev~\cite{egorychev} and Falikman~\cite{falikman}.
Recently, Gurvits~\cite{gurvits-vdW} gave a different proof 
using stability.  Not only does this represent a considerable 
simplification, but the result is general enough to imply 
several other well known results.  
Gurvits begins by identifying the permanent as a coefficient of 
a polynomial, as follows.  If $A = (a_{ij})$ is any $n \times n$ 
matrix, we may define a homogeneous polynomial $p$ by 
$$p(\xx) = p_A (\xx) := \prod_{i=1}^n \sum_{j=1}^n a_{ij} x_j \, .$$
The permanent of $A$ is then the $(1, \ldots , 1)$ coefficient of $p$.
Thus,
\begin{equation} \label{eq:per}
\Per (A) = \frac{\partial^n}{\partial x_1 \cdots \partial x_n}
   p (0 , \ldots , 0) \, .
\end{equation} 
When $A$ is stochastic, each factor in the product evaluates 
to~1 at $(1, \ldots , 1)$.  Taking the derivative with respect 
to $x_j$ and evaluating at $(1,\ldots , 1)$ gives the $j^{th}$ 
column sum, hence if $A$ is doubly stochastic then 
\begin{equation} \label{eq:dbl stoch}
\frac{\partial}{\partial x_j} p_A (1 , \ldots , 1) = 1, \;\; 
   1 \leq j \leq n \, .
\end{equation}
Let $\CC_n$ be the class of homogeneous polynomials of 
degree $n$ in $n$ variables with all coefficients nonnegative.  
Following Gurvits, we extend the definition of the term 
``doubly stochastic'' from matrices to $\CC_n$ by saying
that $p$ is doubly stochastic if~\eqref{eq:dbl stoch} holds.

To prove the van der Waerden result, we need to derive
$\Per (A) \geq n! / n^n$
from~\eqref{eq:dbl stoch}.  This will not be true for every
$p \in \CC_n$ satisfying~\eqref{eq:dbl stoch} but Gurvits' idea
was that it should hold for all stable $p \in \CC_n$.  Because
stability is closed under product and each linear polynomial
with positive coefficients is stable, we see immediately that
$p_A$ is always stable.  The other ingredient in Gurvits' proof
is to strengthen the induction, replacing the lower bound of 
$n! / n^n$ by $(n! / n^n) \Cap (p)$, where $\Cap (p)$ is a
constant that evaluates to~1 when $p = p_A$ for a doubly
stochastic matrix $A$.  

\begin{defn} \label{def:cap}
Let $\CC_n$ be the class of homogeneous polynomials of 
degree $n$ in $n$ variables with all coefficients nonnegative.  
For $p \in \CC_n$, define the capacity of $p$ by
$$\Cap (p) := \inf \frac{p(x_1 , \ldots , x_n)}{\prod_{i=1}^n x_i}$$
where the infimum is over nonnegative values of the variables $x_i$.
\end{defn}

\begin{pr}[\protect{\cite[Fact~2.2]{gurvits-vdW}}]
If $p$ is doubly stochastic then $\Cap (p) = 1$.
\end{pr}

\noindent{\sc Proof:}  
We may view $p \in \CC_n$ as a generating function for a 
probability distribution $\mu$ supported on $\Delta_n := 
\{ \omega \in \Z_+^n : |\omega| = n \}$.  Letting $\E$ denote 
expectation with respect to $\mu$ we see that~\eqref{eq:dbl stoch} 
is equivalent to $\E \omega_j = 1$ for all $j$.

For $\xx \in \R_+^n$, we may evaluate
$p(\xx) = \log \E \prod_j x_j^{\omega_j}$.  Applying Jensen's 
inequality to the concave function $\log (x)$ gives
\begin{eqnarray*}
\log p (x_1 , \ldots , x_n) & = & \log \E \prod_{j=1}^n x_j^{\omega_j} \\
& \geq & \E \log \prod_{j=1}^n x_j^{\omega_j} \\
& = & \E \sum_j \omega_j \log x_j \, .
\end{eqnarray*}
Because $\E \omega_j = 1$ for all $j$, this last quantity is just
$\log \prod_j x_j$.  Thus, on the positive orthant, 
$p(\xx) \geq \prod_{j=1}^n x_j$ showing that $\Cap (p) \geq 1$.
Setting $x_j = 1$ for all $j$ shows that $\Cap (p) \leq 1$.
$\Cox$

Recalling the lower bound on $f' (0)$ for univariate 
stable polynomials with nonnegative coefficients from
Proposition~\ref{pr:1 bound} it is not hard to envision
some kind of induction.  The engine of Gurvits' proof,
used in the induction step, is the following inequality.

\begin{lem}[\protect{\cite[Theorem~4.10]{gurvits-vdW}}] \label{lem:gurvits}
Let $p \in \CC_n$ be stable and define 
$$q (x_1 , \ldots , x_{n-1}) := \frac{\partial}{\partial x_n} 
   p(x_1 , \ldots , x_{n-1} , 0) \, .$$
Then
\begin{equation} \label{eq:gurv}
\Cap (q) \geq \Cap (p) \cdot \left ( \frac{d-1}{d} \right )^{d-1}
\end{equation}
where $d$ is the maximum degree of $x_n$ in $p$.
\end{lem}

\noindent{\sc Proof:} Letting $x_1 , \ldots , x_{n-1}$ range over
positive numbers whose product is~1, we need to show that
$$q(x_1 , \ldots , x_{n-1}) \geq \Cap (p) G(d)$$
where $d \leq n$ is the maximum degree of $x_n$ in $p$ and $G(m) :=
[(m-1)/m]^{m-1}$.  Fix $x_1 , \ldots , x_{n-1}$.  The specialization
property implies that the univariate polynomial
$R(t) := p(x_1 , \ldots , x_{n-1}, t)$ is stable.
By definition of capacity, $\Cap (p) \leq \inf_{t > 0} R(t) / t$. 
The degree of $R$ is equal to $d$, whence it follows from
Proposition~\ref{pr:1 bound} that
$$q(x_1 , \ldots , x_{n-1}) = R' (0) \geq \Cap (p) G(d) \, .$$
$\Cox$

We may now give Gurvits' result implying the van der Waerden
conjecture.
\begin{thm} \label{th:gurvits}
Let $p$ be stable with nonnegative coefficients.  Then
$$\frac{\partial^n}{\partial x_1 \cdots \partial x_n}
   p (0 , \ldots , 0) \geq \Cap (p) \frac{n!}{n^n} \, .$$
\end{thm}

\noindent{\sc Proof:}  Let $q_n = p$ and in general define
$$q_i (x_1 , \ldots , x_i) := \frac{\partial^{n-i}}{\partial x_{i+1}
   \cdots \partial x_n} p(x_1 , \ldots , x_i , 0 \ldots , 0) \, .$$
Stability is closed under differentiation and specializing to
real values, hence by induction on $n-i$, each $q_i$ is stable.
Also, $q_i \in \CC_i$ for each $i$.  Applying Lemma~\ref{lem:gurvits}
with $q_i$ in place of $p$ and $i$ in place of $n$ shows that 
$$\Cap (q_{i+1}) \geq \Cap (q_i) G(i)$$
because $i$ is an upper bound for the degree of $x_i$ in $q_i$.
Inductively, we see that
\begin{equation} \label{eq:both}
\Cap (q_1) \geq \Cap (p) \prod_{i=2}^n G(i)  \, .
\end{equation}
The right-hand side may be identified as
\begin{equation} \label{eq:RHS}
\prod_{i=2}^n G(i) = n! / n^n \, . 
\end{equation}
To identify the left-hand side, observe that $q_1 (x_1)$ 
is homogeneous of degree~1 in $x_1$, that is, $q_1 = c x_1$.  
Thus 
\begin{equation} \label{eq:LHS}
\Cap (q_1) = c = q_1' (0) = 
   \frac{\partial^n}{\partial x_1 \cdots \partial x_n}
   p (0 , \ldots , 0)  \, .
\end{equation}
Plugging~\eqref{eq:LHS} and~\eqref{eq:RHS} into~\eqref{eq:both}
proves the theorem.
$\Cox$

Another result that succumbs to this method is the Schrijver--Valiant
conjecture, proved in 1998.  This concerns an integer version of
doubly stochastic matrices.  Let $\Lambda (k,n)$ denote the collection
of $n \times n$ matrices whose entries are nonnegative integers
and whose rows and columns all sum to $k$.  Define
$$\lambda (k,n) = k^{-n} \min \{ \Per (A) : A \in \Lambda (k,n) \}$$
and 
$$\theta (k) = \lim_{n \to \infty} \lambda(k,n)^{1/n} \, .$$
In 1980 it was proved~\cite{schrijver-valiant} that
$\theta (k) \leq G(k)$ and equality was conjectured. 
The proof by Schrijver~\cite{schrijver} was difficult.
This result is a corollary of Theorem~\ref{th:gurvits}.
For details on this, as well as a separate application
of Theorem~\ref{th:gurvits} to prove a lower bound on 
the mixed discriminant, see~\cite{gurvits-vdW}.

\subsection{The Monotone Column Permanent Conjecture}

Say that the $n \times n$ matrix $A$ is a \Em{monotone column matrix}
if its entries are real and weakly decreasing down each column, 
that is, $a_{i,j} \geq a_{i+1,j}$ for $1 \leq i \leq n-1$ and 
$1 \leq j \leq n$.  Let $J_n$ denote the $n \times n$ 
matrix of all ones.  It was conjectured in~\cite{haglund-ono-wagner}
that whenever $A$ is a monotone column matrix, the univariate 
polynomial $\Per (z J_n + A)$ has only real roots.  A proof
was given there for the case where $A$ is a zero-one matrix.

Although there was strong intuition already present, the proof
had to wait for the development of multivariate stable function
theory.  In 2009, Br{\"a}nd{\'e}n, Haglund, Visontai and 
Wagner~\cite{monotone-permanent} were able to prove this 
conjecture by proving something stronger, namely multivariate 
stability.
\begin{thm}[Multivariate Monotone Column Permanent Theorem (MMCPT)]
\label{th:MMCPT}
Let $Z_n$ be the diagonal matrix whose entries are the $n$
indeterminates $z_1 , \ldots , z_n$.  Let $A$ by an $n \times n$
monotone column matrix.  Then $\Per (J_n Z_n + A)$ is
a stable polynomial in the variables $z_1 , \ldots , z_n$.
Specializing to $z_j \equiv z$ for all $j$ preserves stability,
hence the original conjecture follows.  
\end{thm}

I will give only a brief sketch of the proof, which
relies on two results proved by Borcea and Br{\"a}nd{\'e}n.  
One of these is the criterion, Theorem~\ref{th:multi-affine PS} 
for a linear operator on multi-affine polynomials to be stability 
preserving.  They use the sufficiency direction, namely that 
stability of the $(2n)$-variable ``algebraic symbol'' $G_T (\zz , \ww)$ 
implies that $T$ is a stability preserving operator on 
$\C_{\rm ma} [z_1 , \ldots , z_d]$; this is the direction for
which a proof was included above.  The other lemma they require
is the following result.  
\begin{lem}[\protect{\cite[Proposition~2.5]{monotone-permanent}}]
\label{lem:multilinear}
Let $V$ be a real vector space and let $\phi : V^n \to \R$
be a multilinear form.  Let $e_1 , \ldots , e_n , v_2 , \ldots , v_n$
be fixed vectors in $V$ and suppose that the polynomial
$$\phi (e_1 , v_2 + z_2 e_2 , \ldots , v_n + z_n e_n)$$
is not identically zero in $\R [z_2 , \ldots , z_n]$.  Then
the set of all $v_1 \in V$ for which the polynomial
$$\phi (v_1 + z_1 e_1 , \ldots , v_n + z_n e_n)$$
is stable is either empty or a convex cone over the origin
containing $\pm e_1$.
$\Cox$
\end{lem}

By means of this lemma, the MMCPT may be reduced to the special
case where the entries of $A$ are all zero or one.  I will not
reproduce this reduction here, but the main idea is an induction.
If we have proved the result in the case where $k$ columns are 
real monotone and $n-k$ are zero-one monotone, then applying the 
lemma to one of the $n-k$ columns shows that stability is 
achieved over a convex cone containing all zero-one monotone
columns, and such a cone necessarily contains all real
monotone columns.  

After a change of variables to $y_j = 1 + 1/z_j$ and the 
introduction of new variables $x_1 , \ldots , x_n$, 
stability of $\Per (z_j + a_{ij})$ will follow if we show
stability of $\Per (a_{ij} y_j + (1-a_{ij} x_i)$ for
monotone zero-one matrices.  Denote this last permanent
as $\Per (B(A))$.  Expanding the permanent
on the last row gives a differential recurrence relation
which may be written in form $f_j = T f_{j-1}$ where
$f_j$ is the permanent of the upper $j \times j$ 
submatrix of $B(A)$.  Here $T$ is an operator of the
form $k + z_j \sum_{i=1}^{j-1} \partial / \partial z_j$.
The sufficiency criterion reduces the task to checking
stability of the algebraic symbol, $G_T$, which may be 
accomplished by a simple computation.



\subsection{The BMV conjecture}

In 1975, Bessis, Moussa and Villani~\cite{BMV} formulated what
is now know as the BMV conjecture.
\begin{conj}[BMV] \label{conj:BMV}
Let $A$ and $B$ be Hermitian $n \times n$ matrices with $B$
positive semi-definite.  Then the function
$$\lambda \mapsto \trace \exp \left ( A - \lambda B \right )$$
is the Laplace transform of a positive measure on $[0,\infty)$;
here, $\trace ( \cdot )$ denotes the trace.
\end{conj}
This seems vaguely tied to a number of the themes we have been
discussing, but not necessarily related in a direct way to the
theory of stable functions.  In fact, both the conclusion and 
the hypotheses will reveal more upon further scrutiny.  Beginning
with the obvious, we recall that the linear cone of 
positive semi-definite matrices 
is a cone of hyperbolicity for the determinant function on the 
space of Hermitian matrices.  Thus if we take the determinant 
instead of the trace of the exponential, the resulting function 
$\lambda \mapsto \Det (A - \lambda B)$ is stable.  

Recall from Edrei's Equivalence Theorem (Theorem~\ref{th:edrei})
that stability of the polynomial $f := \sum_{j=0}^k a_j z^j$ is
equivalent to a sequence of inequalities, which can be summed
up by saying that $\{ a_k \}$ is a P{\'o}lya frequency sequence.
Similarly, a function $\phi$ on $[0,\infty)$ is a Laplace
transform of a positive measure on $\R^+$ if and only $\phi$
is a \Em{completely monotone} function, meaning that the derivatives
of $\phi$ do not change sign on $[0,\infty)$ and alternate in sign:
for each nonnegative integer $n$ and positive real $t$,
$$(-1)^n \frac{d^n}{dt^n} \phi (t) \geq 0 \, .$$
This result is the Bernstein--Widder Theorem~\cite{bernstein28}. 
Evidently the conclusion is a property whose definition has
features in common with stability, a property we know to be
true for $\lambda \mapsto \Det (A - \lambda B)$.

In 2004, Lieb and Seiringer~\cite{lieb-seiringer} found an 
equivalent formulation of the BMV conjecture that looks even 
more similar to stability theory.
\begin{thm} \label{th:lieb}
The BMV conjecture is equivalent to the polynomials
$$\trace \left ( A + \lambda B \right )^n$$
having nonnegative coefficients for all integers $n \geq 1$
and all pairs of matrices $A,B$ that are {\em both} positive
semi-definite.
\end{thm}
The exponential is gone, the coefficient on $\lambda$ is positive, 
only integer powers are involved, and what is needed is a 
countable sequence of inequalities: all coefficients in a
sequence of polynomials must be nonnegative.  Another
equivalence brings the BMV conjecture directly into the
realm of stable function theory.  Let $p(x,y)$ be any stable
bivariate polynomial with nonnegative coefficients and
nonzero constant term.  Let $a_{mn}$ be the Taylor coefficients
of $(\partial p / \partial x) / p$.  The BMV conjecture is 
equivalent to the proposition that $(-1)^{m+n} a_{mn} \geq 0$
for all $m,n$.

In July, 2011,
a proof of the BMV conjecture by H. Stahl~\cite{stahl} was posted 
on arXiv.  At the time of writing, the refereed publication has 
not appeared, but word on the street is that the proof should be 
correct.  The proof is sufficiently complicated that understanding 
its specialization to the $3 \times 3$ case is considered
a good project, perhaps for a dissertation.   There is a sense 
that there may be a simpler proof out there, and that the 
``book proof'' of the BMV-Stahl Theorem is likely to involve 
stable polynomials.

\section*{Acknowledgements}

I owe a huge debt to Yuliy Baryshnikov for helping me to understand the
material in Part~I and to Petter Br{\"a}nd{\'e}n for helping me to
understand the material in Part~II.  Without their help, this survey
could not have been written.  Thanks are also due to Mirk{\'o} Visontai 
for his careful reading and comments, to Alan Sokal for crucial help
and corrections on the development of multivariate stability, and to 
David Wagner for helpful discussions.

\bibliographystyle{alpha}
\bibliography{RP}

\end{document}